\documentclass[12pt]{amsart}
\usepackage[hyperindex]{hyperref} % Only for ArXiv
\RequirePackage{amsmath}
\RequirePackage{amsopn}
\RequirePackage{amsfonts}
\RequirePackage{amssymb}
\RequirePackage{amsthm}
\RequirePackage{amsrefs}
\evensidemargin=\oddsidemargin
\makeindex

\newtheorem{thm}{Theorem}[section]

\newtheorem{cor}[thm]{Corollary}
\newtheorem{lemma}[thm]{Lemma}
\newtheorem{prop}[thm]{Proposition}

\theoremstyle{definition}
\newtheorem{definition}[thm]{Definition}

\theoremstyle{remark}
\newtheorem{remark}[thm]{Remark}
\newtheorem{example}[thm]{Example}

\numberwithin{equation}{section}
\def\mathcs{C^{*}}
\newcommand{\cs}{\ensuremath{\mathcs}}

\DeclareMathSymbol{\rtimes}{\mathbin}{AMSb}{"6F}
\def\R{\mathbf{R}}
\def\C{\mathbf{C}}
\def\T{\mathbf{T}}
\def\Z{\mathbf{Z}}
\def\K{\mathcal{K}}
\DeclareMathOperator*{\supp}{supp}

\DeclareMathOperator{\id}{id}
\def\set#1{\{\,#1\,\}}
\let\tensor=\otimes
\def\restr#1{|_{{#1}}}

\newbox\hidebox
\def\spechide#1{\setbox\hidebox=\hbox{$#1$}
\hbox to\wd\hidebox{$\box\hidebox^\wedge$\hss}}
\makeatletter
\def\labelenumi{\textnormal{(\@alph\c@enumi)}}
\def\theenumi{\@alph \c@enumi}
\def\labelenumii{\textnormal{(\@roman\c@enumii)}}
\def\theenumii{\@roman \c@enumii}
\newcount\charno
\def\alphapart#1{\charno=96
\advance\charno by#1\char\charno}
\def\partref#1{\textup{(}\textnormal{\alphapart{#1}}\textup{)}}
\makeatother
\def\<{\langle}
\def\>{\rangle}
\let\ipscriptstyle=\scriptscriptstyle
\def\lipsqueeze{{\mskip -3.0mu}}
\def\ripsqueeze{{\mskip -3.0mu}}
\def\ipcomma{\nobreak\mathrel{,}\nobreak}
\newbox\ipstrutbox
\setbox\ipstrutbox=\hbox{\vrule height8.5pt% depth 3.5pt
width 0pt}
\def\ipstrut{\copy\ipstrutbox}
\def\lip#1<#2,#3>{\mathopen{\relax_{\ipstrut\ipscriptstyle{
#1}}\lipsqueeze
\langle} #2\ipcomma #3 \rangle}
\def\blip#1<#2,#3>{\mathopen{\relax_{\ipstrut
\ipscriptstyle{ #1}}\lipsqueeze\bigl\langle} #2\ipcomma #3 \bigr\rangle}
\def\rip#1<#2,#3>{\langle #2\ipcomma #3
\rangle_{\ripsqueeze\ipstrut\ipscriptstyle{#1}}}
\def\brip#1<#2,#3>{\bigl\langle #2\ipcomma #3
\bigr\rangle_{\ripsqueeze\ipstrut\ipscriptstyle{#1}}}
\def\angsqueeze{\mskip -6mu}
\def\smangsqueeze{\mskip -3.7mu}
\def\trip#1<#2,#3>{\langle\smangsqueeze\langle #2\ipcomma #3
\rangle\smangsqueeze\rangle_{\ripsqueeze\ipstrut\ipscriptstyle{#1}}}
\def\btrip#1<#2,#3>{\bigl\langle\angsqueeze\bigl\langle #2\ipcomma
#3
\bigr\rangle
\angsqueeze\bigr\rangle_{\ripsqueeze\ipstrut\ipscriptstyle{#1}}}
\def\tlip#1<#2,#3>{\mathopen{\relax_{\ipstrut\ipscriptstyle{
#1}}\lipsqueeze \langle\smangsqueeze\langle} #2\ipcomma #3
\rangle\smangsqueeze\rangle}
\def\btlip#1<#2,#3>{\mathopen{\relax_{\ipstrut\ipscriptstyle{
#1}}\lipsqueeze
\bigl\langle\angsqueeze\bigl\langle} #2\ipcomma #3
\bigr\rangle\angsqueeze\bigr\rangle}

\def\ip(#1|#2){(#1\mid #2)}
\def\bip(#1|#2){\bigl(#1 \mid #2\bigr)}
\def\Bip(#1|#2){\Bigl( #1 \bigm| #2 \Bigr)}
\expandafter\ifx\csname BibSpec\endcsname\relax\else
\BibSpec{collection.article}{%
    +{}  {\PrintAuthors}                {author}
    +{,} { \textit}                     {title}
    +{.} { }                            {part}
    +{:} { \textit}                     {subtitle}
    +{,} { \PrintContributions}         {contribution}
    +{,} { \PrintConference}            {conference}
    +{}  {\PrintBook}                   {book}
    +{,} { }                            {booktitle}
    +{,} { }                            {series}
    +{,} { \voltext}                    {volume}
    +{,} { }                            {publisher}
    +{,} { }                            {organization}
    +{,} { }                            {address}
    +{,} { \PrintDateB}                 {date}
    +{,} { pp.~}                        {pages}
    +{,} { }                            {status}
    +{,} { \PrintDOI}                   {doi}
    +{,} { available at \eprint}        {eprint}
    +{}  { \parenthesize}               {language}
    +{}  { \PrintTranslation}           {translation}
    +{;} { \PrintReprint}               {reprint}
    +{.} { }                            {note}
    +{.} {}                             {transition}
}
\BibSpec{article}{%
    +{}  {\PrintAuthors}                {author}
    +{,} { \textit}                     {title}
    +{.} { }                            {part}
    +{:} { \textit}                     {subtitle}
    +{,} { \PrintContributions}         {contribution}
    +{.} { \PrintPartials}              {partial}
    +{,} { }                            {journal}
    +{}  { \textbf}                     {volume}
    +{}  { \PrintDatePV}                {date}
    +{,} { \eprintpages}                {pages}
    +{,} { }                            {status}
    +{,} { \PrintDOI}                   {doi}
    +{,} { available at \eprint}        {eprint}
    +{}  { \parenthesize}               {language}
    +{}  { \PrintTranslation}           {translation}
    +{;} { \PrintReprint}               {reprint}
    +{.} { }                            {note}
    +{.} {}                             {transition}
}
\BibSpec{book}{%
    +{}  {\PrintPrimary}                {transition}
    +{,} { \textit}                     {title}
    +{.} { }                            {part}
    +{:} { \textit}                     {subtitle}
    +{,} { \PrintEdition}               {edition}
    +{}  { \PrintEditorsB}              {editor}
    +{,} { \PrintTranslatorsC}          {translator}
    +{,} { \PrintContributions}         {contribution}
    +{,} { }                            {series}
    +{,} { \voltext}                    {volume}
    +{,} { }                            {publisher}
    +{,} { }                            {organization}
    +{,} { }                            {address}
    +{,} { pp.~}                        {pages}
    +{,} { \PrintDateB}                 {date}
    +{,} { }                            {status}
    +{}  { \parenthesize}               {language}
    +{}  { \PrintTranslation}           {translation}
    +{;} { \PrintReprint}               {reprint}
    +{.} { }                            {note}
    +{.} {}                             {transition}
}
\fi

\newcommand{\cox}{\ensuremath{C_{0}(X)}}

\IfFileExists{mathrsfs.sty}{\usepackage{mathrsfs}}
{\let\mathscr\mathcal} 
\newcommand{\ib}{im\-prim\-i\-tiv\-ity bi\-mod\-u\-le}

\newcommand{\sme}{\,\mathord{\mathop{\text{--}}\nolimits_{\relax}}\,}
\newcommand{\modulefont}[1]{\mathsf{#1}}
\newcommand{\X}{\modulefont{X}}

\newcommand{\bundlefont}[1]{\mathscr{#1}}
\newcommand{\A}{\bundlefont A}
\newcommand{\B}{\bundlefont B}
\newcommand{\E}{\bundlefont E}
\newcommand\CC{\bundlefont C}
\newcommand{\HH}{\bundlefont H}
\newcommand\VV{\bundlefont V}

\def\cp(#1,#2,#3){#1\rtimes_{#3}#2}

\newcommand{\go}{G^{(0)}}
\newcommand{\ho}{H^{(0)}}
\def\sa_#1(#2;#3){\Gamma_{#1}(#2;#3)}

\def\bchip<#1,#2>{\tlip\cp(\B,H,\beta)<{#1},{#2}>}
\def\bbchip<#1,#2>{\btlip\cp(\B,H,\beta)<{#1},{#2}>}
\def\acgip<#1,#2>{\trip\cp(\A,G,\alpha)<{#1},{#2}>}
\def\bacgip<#1,#2>{\btrip\cp(\A,G,\alpha)<{#1},{#2}>}

\newcommand{\op}[1]{{#1}^{\text{\normalfont op}}}

\renewcommand{\H}{\mathcal{H}}
\newcommand{\half}{\frac12}
\newcommand{\hoo}{\H_{00}}
\let\phi\varphi
\usepackage{bbold}
\def\charfcn#1{\mathbb{1}_{#1}}

\newcommand{\usc}{upper semicontinuous}
\newcommand{\stars}{*_{s}}
\newcommand{\starr}{*_{r}}

\def\sacc(#1;#2){\mathscr{G}(#1;#2)}

\newtheorem*{claim}{Claim}

\newcommand{\bb}{\mathcal{B}^{b}}

\newcommand\atensor{\odot}

\def\pip<#1,#2>{\brip<[#1],[#2]>}

\newcommand\gcb{\sa_{c}(G;\B)}
\newcommand\grm{generalized Radon measure}
\newcommand\ccg{C_{c}(G)}
\newcommand\ccpg{C_{c}^{+}(G)}
\newcommand\nnu{|\nu|} 

\theoremstyle{plain}

\theoremstyle{definition}

\theoremstyle{remark}
\newtheorem*{conventions}{Conventions}

\newcommand\prerep{pre-representation}
\newcommand\ilt{inductive limit topology}
\newcommand\End{\operatorname{End}}

\newcommand\gcbatho{\gcb\atensor\H_{0}}
\newcommand\sgcbatho{\sgcb\atensor\H_{0}}
\newcommand\N{\mathscr{N}}

\newcommand\sgcb{\Sigma_{c}^{1}(G;\B)}
\newcommand\bboc{\mathcal{B}^{1}_{c}}

\renewcommand\AA{\mathcal{A}_{0}}

\newcommand\D{\mathscr{D}}

\newcommand\nsigma{|\sigma|}
\newcommand\hoop{\hoo'}
\newcommand\eij{\epsilon^{ij}}
\newcommand\gcbathoop{\gcb\atensor\hoop}

\newcommand\strutornot{\relax}
\def\ipu<#1,#2>{\langle #1 \ipcomma #2 \rangle_{\strutornot u}}
\def\ipv<#1,#2>{\langle #1 \ipcomma #2 \rangle_{\strutornot v}}
\def\bipu<#1,#2>{\bigl\langle #1 \ipcomma #2 \bigr\rangle_{\strutornot u}}
\def\bipv<#1,#2>{\bigl\langle #1 \ipcomma #2 \bigr\rangle_{\strutornot
    v}}
\def\bipb #1<#2,#3>{\bigl\langle #2 \ipcomma #3 \bigr\rangle_{\strutornot #1}}

\newcommand\hatpi{\check\pi}

\newcommand\Lin{\operatorname{Lin}}

\newcommand\spxi{\zeta}

\newcommand\gucb{\sa_{c}(G^{u};\B\restr{G^{u}})}

\newcommand\Lambdasubc{\Lambda_{c}}

\newcommand\hg{\lambda_{G}}
\newcommand\hh{\lambda_{H}}

\def\fcs(#1,#2){\cs(#1,#2)}

\newcommand\tlips{\tlip{\scriptstyle\star}}
\newcommand\trips{\trip{\scriptstyle\star}}
\newcommand\myf{\Upsilon}

\renewcommand\op[1]{{#1}^{\text{op}}}
\newcommand\bq{\bar q}
\newcommand\bE{\overline{\E}}
\newcommand\dmap{\flat}

\newcommand\tlipsd{\tlip{\bar{\scriptstyle\star}}}
\newcommand\tripsd{\trip{\bar{\scriptstyle\star}}}

\usepackage[all]{xy} 
\SelectTips{cm}{} 
\allowdisplaybreaks[3]

\date{3 June 2008}


\begin{document}

\title[Equivalence of Fell Bundles]{\boldmath Equivalence and
  Disintegration Theorems for Fell Bundles and their \cs-algebras}

\author[Muhly]{Paul S. Muhly}
\address{Department of Mathematics\\
  University of Iowa\\
  Iowa City, IA 52242} \email{pmuhly@math.uiowa.edu}

\author[Williams]{Dana P. Williams}
\address{Department of Mathematics\\
  Dartmouth College\\
  Hanover, NH 03755-3551} \email{dana.williams@dartmouth.edu}

\begin{abstract}
We study the \cs-algebras of Fell bundles.  In particular, we 
prove the analogue of Renault's disintegration theorem for groupoids.
As in the groupoid case, this result is the key step in proving a deep
equivalence theorem for the \cs-algebras of Fell bundles.
\end{abstract}

\maketitle

\tableofcontents

\section*{Introduction}
\label{sec:introduction}

Induced representations and the corresponding imprimitivity theorems
constitute a substantial part of the representation theory of locally
compact groups and play a critical role in harmonic analysis.  These
constructs extend naturally to the crossed products of \cs-algebras by
locally compact groups as illustrated in \cite{wil:crossed}.  One can
push the envelope considerably further to include twisted crossed
products of various flavors, and eventually arrive at the Banach
$*$-algebraic bundles of Fell.  (The latter are discussed in
considerable detail in the two volume treatise of Fell and Doran
\citelist{\cite{fd:representations1}\cite{fd:representations2}}.)
There are also extensions of these concepts to locally compact
groupoids, groupoid crossed products, and as we will describe here, to
the groupoid analogue of Fell's Banach $*$-algebraic bundles over
groups, which we call just \emph{Fell bundles}.  These bundles, their
representation theory and their role in harmonic analysis will be our
focus here.

In our approach, induction and imprimitivity are formalized using the
Rieffel machinery as described in \cite{rw:morita}.  Therefore our
results are stated in the form of Morita equivalences, and induction
is defined algebraically as in \cite{rw:morita}*{\S2.4}.  In this
setting, the central result is Raeburn's Symmetric Imprimitivity
Theorem
\citelist{\cite{rae:ma88}\cite{wil:crossed}*{Chap.~4}}.\footnote{The
  imprimitivity theorem for Fell's Banach $*$-algebraic bundles,
  \cite{fd:representations2}*{Theorem~XI.14.17}, is not, strictly
  speaking, a Morita equivalence result.  Of course, such a result
  will be a consequence of our equivalence theorem:
  Theorem~\ref{thm-yam-equi}.}

However, the there are daunting technical difficulties to surmount.
For example, the correspondence between representations of a groupoid
\cs-algebra and ``unitary'' representations of the underlying groupoid
is, unlike in the group case, a deep result which is known as
Renault's Disintegration Theorem \cite{ren:jot87}*{Proposition~4.2}
(see also, \cite{muhwil:nyjm08}*{Theorem~7.8} for another proof).  The
formulation of Renault's Disintegration Theorem is very general, and
is designed to facilitate the proof of the Equivalence Theorem
\cite{mrw:jot87}*{Theorem~2.8} which states, in essence, that
equivalent groupoids have Morita equivalent \cs-algebras.  The
equivalence theorem can be extended to groupoid crossed products
\cite{ren:jot87}*{Corollaire~5.4}.  (Groupoid crossed products and
Renault's equivalence theorem are discussed at length in
\cite{muhwil:nyjm08}.)  The equivalence theorem for groupoid crossed
products is a very powerful imprimitivity type theorem with far
reaching consequences.  For example, it subsumes Raeburn's Symmetric
Imprimitivity Theorem \cite{muhwil:nyjm08}*{Example~5.12}.

In his 1987 preprint, Yamagami proposed a natural generalization of
Fell's Banach $*$-algebraic bundles over a locally compact group to a
Fell bundle over a locally compact groupoid.  He also suggested that
there should be a disintegration theorem \cite{yam:xx87}*{Theorem~2.1}
and a corresponding equivalence theorem \cite{yam:xx87}*{Theorem~2.3}.
However he gave only bare outlines of how proofs might be constructed.
The object of this paper is to work out carefully the details of these
results in a slightly more general context (as indicated in the next
paragraph).  As will become clear, there are significant technical
hurdles to clear.  In the first part of this paper, we want to
formalize the notion of a Fell bundle and its corresponding
\cs-algebra, and then to prove the disintegration theorem.  This will
require considerable new technology which we develop here.  The
remainder of the paper is devoted to describing the appropriate notion
of equivalence for Fell bundles, and then to stating and proving the
equivalence theorem.

Since a groupoid can only act on objects which are fibred over its
unit space, it has been clear from their inception that groupoid
crossed products must involve \cs-bundles of some sort.  Moving to
equivalence theorems and Fell bundles means that we will have to widen
the scope to include Banach bundles.  Originally, Renault and
Yamagami, for example, worked with continuous Banach bundles of the
sort studied by Fell (cf.,
\cite{fd:representations1,fd:representations2}).  However, more
recently (for example, see \cite{muhwil:nyjm08}), it has become clear
that it is unnecessarily restrictive to insist on continuous Banach
bundles.  Rather, the appropriate notion is what we call an
\emph{\usc-Banach bundle}.\index{upper semicontinuous-Banach
  bundle@\usc-Banach bundle} Fell called such bundles \emph{loose
  Banach bundles} (\cite{fd:representations1}*{Remark~C.1}) and they
are called \emph{(H)~Banach bundles} by Dupr\'e and Gillette
(\cite{dg:banach}*{p.~8}).\index{loose Banach bundle}\index{H-Banach
  bundle@(H)-Banach bundle} For
convenience, we have collated some of the relevant definitions and
results in Appendix~A.  The case of \usc{} \cs-bundles\index{upper
  semicontinuous-\cs-bundle@\usc-\cs-bundle} is treated in
more detail in \cite{wil:crossed}*{Appendix~C}.  As an illustration of
the appropriateness of \usc{} bundles in the theory, we mention the
fact that every \cox-algebra $A$ is the section algebra of an \usc{}
\cs-bundle over $X$ \cite{wil:crossed}*{Theorem~C.26}.  Furthermore,
in practice, most results for continuous bundles extend to \usc{}
bundles without significant change.

Extending the groupoid disintegration theorem to the setting of Fell
bundles is a formidable task.  The first obstacle one faces is the
necessity of developing a useful theory of generalized Radon measures
for linear functionals on the section algebras of \usc-Banach bundles.
This requires a not-altogether-straightforward extension to our
setting of a result of Dinculeanu
\cite{din:integration74}*{Theorem~28.32} for continuous Banach
bundles.  Some not-so-standard facts about complex Radon measures are
also required.
% and some of these properties are summarized in
% \cite{muhwil:nyjm08}*{Appendix~A.1}.

We must also cope with the reality that Ramsay's selection theorems
(\cite{ram:am71}*{Theorem~5.1} and \cite{ram:jfa82}*{Theorem~3.2}) ---
which show, for example, that an almost everywhere homomorphism is
equal almost everywhere to a homomorphism --- are not available for
Fell bundles.  Instead, we must finesse this with new techniques.
(These techniques are valuable in the scalar case as well; for
example, they play a role in the proof of Renault's Disintegration
Theorem given in Appendix~B of \cite{muhwil:nyjm08}.

% The ``philosophy'' of the proof in our setting is that a Fell bundle
% $p:\B\to G$ over $G$ is a realization of $G$ as a category in which
% the ``object'' $u\in\go$ is realized as the \cs-algebra $B(u)=
% p^{-1}(u)$ and the ``arrow'' $x\in G$ is realized as the
% $B\bigl(r(p(x))\bigr)\sme B\bigl(s(p(x)\bigr)$-imprimitivity
% bimodule $B(x):=p^{-1}(x)$. Then the idea is to ramp up the proof in
% the scalar case to achieve our result.  Some support for this point
% of view is that we can recover the scalar disintegration theorem
% from our approach --- see Theorem~\ref{thm-ren-4.2}.  With the
% disintegration result in hand, the path to a proof of an equivalence
% theorem is clear --- if a bit rocky.

The ``philosophy'' of the proof in our setting is that the total
space $\B$ of a Fell bundle $p:\B\to G$ over $G$ may be thought of as
a groupoid in the category whose objects are \cs-algebras and whose
morphisms are (isomorphism classes) of imprimitivity bimodules --- the
composition being induced by balanced tensor products. The bundle map
$p$ may then be viewed as a homomorphism or functor from $\B$ to $G$
Thus, the ``object'' mapped to $u\in\go$ is the \cs-algebra $B(u)=
p^{-1}(u)$ and the ``arrow'' $x\in G$ is the image of the
$B\bigl(r(p(x))\bigr)\sme B\bigl(s(p(x)\bigr)$-imprimitivity bimodule
$B(x):=p^{-1}(x)$. Then the idea is to ramp up the proof in the scalar
case to achieve our result.  Some support for this point of view is
that we can recover the scalar disintegration theorem from our
approach --- see Theorem~\ref{thm-ren-4.2}.  With the disintegration
result in hand, the path to a proof of an equivalence theorem is clear
--- if a bit rocky.

We have made an effort to use standard notation and conventions
throughout.  Because of the use of direct integrals and measure
theory, we want most of the objects with which we work to be
separable.  In particular, $G$ will always denote a second countable
locally compact Hausdorff groupoid with an Haar system
$\set{\lambda^{u}}_{u\in\go}$.  Homomorphisms between \cs-algebras are
always $*$-preserving, and representations of \cs-algebras are assumed
to be nondegenerate.  If $A$ is a \cs-algebra, then $\widetilde A$ is
the subalgebra of the multiplier algebra, $M(A)$, generated by $A$ and
$1_{A}$; that is, $\widetilde A$ is just $A$ if $A$ has an identity,
and it is $A$ with an identity adjoined otherwise.  The term
``pre-compact'' is used to describe a set contained in a compact set.

We start in Section~\ref{sec:definitions} with the definition of a
Fell bundle $p:\B\to G$ and its associated \cs-algebra $\cs(G,\B)$.
Although the definition might seem overly technical on a first
reading, we hope that the examples in Section~\ref{sec:examples} show
that the definition is in fact easy to apply in practice.  In
Section~\ref{sec:radon-measures-fell}, we develop the tools necessary
to associate measures to linear functionals, which we call
\emph{generalized Radon measures}, on the section algebra of an
\usc-Banach bundle.  In Section~\ref{sec:repr-fell-bundl} we formalize
various notions of representations of Fell bundles and formulate the
disintegration result.  We give the proof in
Section~\ref{sec:proof-disint-theor}.  In
section~\ref{sec:equiv-fell-bundl} we formulate the equivalence
theorem which is our ultimate goal.  The proof is given in
section~\ref{sec:proof-theorem} except for some important, but
technical details, about the existence of approximate identities of a
special form.  These details are dealt with in
section~\ref{sec:proof-prop-refpr}.

For convenience, we have included a brief appendix on \usc-Banach
bundles (Appendix~\ref{sec:usc-banach-bundles}).  A more detailed
treatment is available in the original papers
\citelist{\cite{hofkei:lnm79}\cite{hof:lnm77}\cite{hof:74}\cite{dg:banach}}.
An elementary summary for \usc{} \cs-bundles can be found in
\cite{wil:crossed}*{Appendix~C}.  There is also a short appendix
(Appendix~\ref{sec:an-example:-scalar}) showing how to derive the
scalar version of the disintegration theorem from our results.

\section{Preliminaries}
\label{sec:definitions}

There are a number of (equivalent) definitions of Fell bundles in the
literature.  For example, in addition to Yamagami's original
definition in \cite{yam:xx87}*{Definition~1.1}, there is Muhly's
version from \cite{muh:cm01}*{Definition~6}, and more recently
Deaconu, Kumjian and Ramazan have advanced one
\cite{dkr:ms08}*{Definition~2.1}.  Here we give another variant which
is formulated primarily with the goal of being easy to check in
examples rather than being the most succinct or most elegant.  We hope
that the formulation will be validated when we look at examples in
Section~\ref{sec:examples}.

\begin{definition}
  \label{def-fell-bund}
  Suppose that $p:\B\to G$ is a separable \usc-Banach bundle over a
  second countable locally compact Hausdorff groupoid $G$.  Let
  \begin{equation*}
    \B^{(2)}:=\set{(a,b)\in\B\times\B:\bigl(p(a),p(b)\bigr)\in G^{(2)}}.
  \end{equation*}
  We say that $p:\B\to G$ is a \emph{\index{Fell bundle}Fell bundle}
  if there is a continuous bilinear associative ``multiplication'' map
  $m:\B^{(2)}\to\B$, and a continuous involution $b\mapsto b^{*}$ such
  that
  \begin{enumerate}
  \item $p\bigl(m(a,b)\bigr)=p(a)p(b)$,
  \item $p(b^{*})=p(b)^{-1}$,
  \item $m(a,b)^{*}=m(b^{*},a^{*})$,
  \item for each $u\in\go$, the Banach space $B(u)$ is a \cs-algebra
    with respect to the $*$-algebra structure induced by the
    involution and multiplication $m$,
  \item for each $x\in G$, $B(x)$ is a $B\bigl(r(x)\bigr)\sme
    B\bigl(s(x)\bigr)$-\ib{} with the module actions determined by $m$
    and inner products
    \begin{equation*}
      \lip B(r(x))<a,b>=m(a,b^{*})\quad\text{and} \quad \rip
      B(s(x))<a,b> =m(a^{*},b).
    \end{equation*}
  \end{enumerate}

\end{definition}

We will normally suppress the map $m$, and simply write $ab$ in place
of $m(a,b)$.  As pointed out in \cite{muh:cm01}, the axioms imply that
$B(x^{-1})$ is (isomorphic to) the dual module to $B(x)$, and that
formulas like
\begin{gather*}
  \|b^{*}b\|_{B(s(p(b)))}=\|b\|^{2}_{B(p(b))} \quad\text{and}\quad
  b^{*}b\ge0 \quad\text{in $B(s(p(b)))$}
\end{gather*}
hold for all $b\in\B$. Note also that \partref1 implies that
$B(x)B(y)\subset B(xy)$. In fact, we have the following (cf.
\cite{muh:cm01}*{Definition~6(2)}):
\begin{lemma}\label{lem-extra}
  Multiplication induces an \ib{} isomorphism of
  $B(x)\tensor_{B(s(x))} B(y)$ with $B(xy)$.
\end{lemma}
\begin{proof}{\emergencystretch=100pt
  Recall that we are suppressing the map $m$.  Since $m$ is bilinear,
  we obtain a map $\bar m:B(x)\atensor B(y)\to B(xy)$, where
  $B(x)\atensor B(y)$ denotes the algebraic tensor product of $B(x)$
  and $B(y)$ balanced over $B(s(x))=B(r(y))$, and then 
  \begin{align*}
    \blip B(r(x))<\bar m(a\tensor b),\bar m(c\tensor d)>&=
    (ab)(cd)^{*} 
    = a(cdb^{*})^{*} 
    = \blip B(r(x))<a, c {\lip B(r(y))<d,b>}>.
  \end{align*}
  It follows that $\bar m$ maps $B(x)\tensor_{B(s(x))} B(y)$
  isometrically onto a closed $B\bigl(r(x)\bigr)\sme
  B\bigl(s(x)\bigr)$-sub-bimodule $Y$ of $B(xy)$.  Since $\blip
  B(r(x))<Y,Y>$ contains $\blip B(r(x)) < B(x)\atensor
  B(y),B(x)\atensor B(y)>$, the ideal of $B\bigl(r(x)\bigr)$
  corresponding to $Y$ in the Rieffel correspondence is all of
  $B\bigl(r(x)\bigr)$.  Hence $Y$ must be all of $B(xy)$.\par}
\end{proof}

\begin{conventions}
  In fact, to make formulas like the above easier to read,
  if\index{conventions} $b\in\B$, then we will usually write $s(b)$ in
  place of the more cumbersome $s(p(b))$ and similarly for $r(b)$.
\end{conventions}

Given a Fell bundle $p:\B\to G$, we want to make the section
algebra\index{Fell bundle!section algebra}
$\gcb$\index{gammacb@$\gcb$} into a topological $*$-algebra in the
\ilt.  The involution is not controversial.  We define\index{Fell
  bundle!involution}
\begin{equation}
  \label{eq:1}
  f^{*}(x):=f(x^{-1})^{*}\quad\text{for $f\in\gcb$.}
\end{equation}
The product is to be given by the convolution formula\index{Fell
  bundle!convolution} 
\begin{equation}
  \label{eq:2}
  f*g(x):=\int_{G}f(y)g(y^{-1}x)\,d\lambda^{r(x)}(y)\quad\text{for all
    $f,g\in\gcb$.}
\end{equation}
There is no issue seeing that $f*g(x)$ is a well-defined element of
$B(x)$.  Clearly $y\mapsto f(y)g(y^{-1}x)$ is in
$C_{c}(G^{r(x)},B(x))$.  Then \cite{wil:crossed}*{Lemma~1.91}) implies
that the integral converges to an element in the Banach space $B(x)$.
However, it is not so clear that $f*g$ is continuous.  For this we
need the following lemma.

\begin{lemma}
  \label{lem-conv}
  Let $G\starr G=\set{(x,y)\in G\times G:r(x)=r(y)}$ and define
  $q:G\starr G\to G$ by $q(x,y)=x$.  Let $q^{*}\B$ be the pull-back.
  If $F\in \sa_{c}(G\starr G;q^{*}\B)$, then
  \begin{equation*}
    f_{F}(x):=\int_{G}F(x,y)\,d\lambda^{r(x)}(y)
  \end{equation*}
  defines a section $f_{F}\in\gcb$.
\end{lemma}
\begin{proof}
  A partition of unity argument (see Lemma~\ref{lem-dense-sections})
  implies that sections of the form
  \begin{equation*}
    (x,y)\mapsto \psi(x,y)f(x)
  \end{equation*}
  for $\psi\in C_{c}(G\starr G)$ and $f\in \gcb$ span a dense subspace
  of $\sa_{c}(G\starr G;q^{*}\B)$.  Since we can approximate $\psi$ in
  the \ilt{} by sums of the form $(x,y)\mapsto
  \phi_{1}(x)\phi_{2}(y)$, and since
  $\phi_{1}(x)\phi_{2}(y)f(x)=\phi_{2}(y)(\phi_{1}\cdot f)(x)$, it
  follows that sections of the form
  \begin{equation*}
    (x,y)\mapsto \phi(y)f(x)
  \end{equation*}
  with $\phi\in\ccg$ and $f\in\gcb$ span a dense subspace $\AA$ of
  $\sa_{c}(G\starr G;q^{*}\B)$.  Since $f_{F_{i}}\to f_{F}$ uniformly
  if $F_{i}\to F$ in the \ilt{} in $\sa_{c}(G\starr G;q^{*}\B)$, it
  suffices to show that $f_{F}$ is continuous when $F\in\AA$.  But if
  $F(x,y)=\phi(y)f(x)$, then
  \begin{equation*}
    f_{F}(x) =\Bigl(\int_{G}\phi(y)\,d\lambda^{r(x)}(y)\Bigr) f(x)
    =\lambda(\phi)\bigl(r(x)\bigr)f(x) 
  \end{equation*}
  which is clearly in $\gcb$ (because $\lambda(\phi)$ is in
  $C_{c}(\go)$ since $\set{\lambda^{u}}$ is a Haar system).
\end{proof}

\begin{cor}
  \label{cor-conv}
  If $f$ and $g$ are in $\gcb$, then so is $f*g$.
\end{cor}
\begin{proof}
  Note that $(x,y)\mapsto f(y)g(y^{-1}x)$ defines a section in
  $\sa_{c}(G\starr G;q^{*}\B)$.
\end{proof}

Just as for groupoid algebras, the \emph{\index{I-norm}I-norm} on
$\gcb$\index{gammacb@$\gcb$!I-norm} is given
by\index{normI@$"\"|\cdot"\"|_{I}$} 
\begin{equation}
  \label{eq:32}
  \|f\|_{I}=\max\Bigl(\sup_{u\in\go}\int_{G}\|f(x)\|\,d\lambda^{u}(x),
  \sup_{u\in\go} \int_{G}\|f(x)\|\,d\lambda_{u}(x)\Bigr).
\end{equation}
A $*$-homomorphism $L:\gcb\to B(\H_{L})$ is called a \emph{\index{I-norm
  decreasing representation}I-norm
  decreasing representation}\index{representation!I-norm decreasing} if $\|L(f)\|\le\|f\|_{I}$ for all $f\in
\gcb$ and if $\operatorname{span}\set{L(f)\xi:\text{$f\in\gcb$ and
    $\xi\in\H_{L}$}}$ is dense in $\H_{L}$.  The \emph{\index{universal
  norm}universal
  \cs-norm} on $\gcb$\index{gammacb@$\gcb$!universal norm} is given by
\begin{equation*}
  \|f\|:=\sup\set{\|L(f)\|:\text{$L$ is a $I$-norm decreasing
      representation}}. 
\end{equation*}
The completion of $\gcb$ with respect to the universal norm is denoted
by $\cs(G,\B)$.\index{csgb@$\cs(G,\B)$}

\section{Examples}
\label{sec:examples}

We want to review some of the examples of Fell Bundles described by
Muhly in Example~7 of \cite{muh:cm01}*{\S3}. At the same time, we want
to add a bit of detail, and make a few alterations.

\begin{example}[Groupoid Crossed Products]
  \label{ex-cp}
  Let $\pi:\A\to\go$ be an \usc{} \cs-bundle over $\go$.  We assume
  that $(\A,G,\alpha)$ is a groupoid dynamical system.\index{groupoid
    dynamical system}\index{groupoid crossed product}  Unlike the
  treatment in \cite{muh:cm01}*{Example~7(1)} where the focus is on
  $s^{*}\A$, we want to work with the pull-back
  $\B:=r^{*}\A=\set{(a,x):\pi(a)=r(x)}$.  The first step is to define
  a multiplication on
  $\B^{(2)}:=\set{\bigl((a,x),(b,y)\bigr)\in\B\times\B:(x,y)\in
    G^{(2)}}$ as follows:
  \begin{equation}
    \label{eq:4}
    (a,x)(b,y):=\bigl(a\alpha_{x}(b),xy\bigr).
  \end{equation}
  (This formula looks ``even more natural'' if we write $x\cdot b$ in
  place of $\alpha_{x}(b)$.)  The involution is given by
  \begin{equation}
    \label{eq:5}
    (a,x)^{*}:=\bigl(\alpha_{x}^{-1}(a^{*}),x^{-1}\bigr). 
  \end{equation}

  To verify that $B(x)=\set{(a,x):a\in A\bigl(r(x)\bigr)}$ is a
  $A\bigl(r(x)\bigr)\sme A\bigl(s(x)\bigr)$-\ib, we proceed as
  follows.  Keep in mind that the left $A\bigl(r(x)\bigr)$- and right
  $A\bigl(s(x)\bigr)$-actions are determined by \eqref{eq:4}.  Since
  $\alpha_{u}=\id$ if $u\in \go$, we have
  \begin{align*}
    \label{eq:6}
    b\cdot
    (a,x)&:=(b,r(x))(a,x)=(ba,x)\quad\text{and}\\
    (a,x)\cdot c&= (a,x)(c,s(x))=\bigl(a\alpha_{x}(c),x\bigr).
  \end{align*}
  Again, by axiom, the inner products are supposed to be given by
  \begin{align*}
    \brip B(s(x)) <(a,x),(b,x)>&
    =(a,x)^{*}(b,x)=\bigl(\alpha_{x}^{-1}(a^{*}b),s(x)\bigr) \quad
    \text{and} \\
    \blip B(r(x))<(a,x),(b,x)>&= (a,x)(b,x)^{*} = \bigl(ab^{*}
    ,r(x)\bigr)
  \end{align*}
  Of course, these are the natural inner products and actions on
  $A\bigl(r(x)\bigr)$ making it into an \ib.\footnote{Recall that if
    $\theta:A\to B$ is an isomorphism, then there is a natural way to
    make $A$ into an $A\sme B$-\ib{} usually denoted $A_{\theta}$.}

  \begin{remark}
    \label{rem-cp}
    Notice that a section $f\in \sa_{}(G;\B)$
    \index{gammagb@$\sa_{}(G;\B)$}\index{section} is determined by a
    function $\check f:G\to \A$ such that $\check f(x)\in
    A\bigl(r(x)\bigr)$.  Then $f(x)=\bigl(\check f(x),x\bigr)$.  We
    will often not distinguish between $f$ and $\check f$.
  \end{remark}

\end{example}

\begin{example}[Twists]
  \label{ex-twist}
  The notion of a \index{twist}twist $E$ over $G$, or a $\T$-groupoid
  over $G$, is\index{T-groupoid@$\T$-groupoid} 
  due to Kumjian.  Recall that $E$ must be a principal circle bundle,
  say $j:E\to G$, over $G$ and that $E$ is also equipped with a
  groupoid structure such that
  \begin{equation*}
    \xymatrix{\go\ar[r]&\go\times\T\ar[r]^-{i}&E\ar[r]^-{j}&G\ar[r]&\go}
  \end{equation*}
  is a groupoid extension such that $t\cdot e = i\bigl(r(e),t)\bigr)e$
  and $(t\cdot e)(s\cdot f)=(ts)\cdot ef$ (see
  \cite{muhwil:plms395}*{p.~115}).  In this case, we let $\B$ be the
  complex line bundle over $G$ associated to $E$.\footnote{One obtains
    $\B$ from $E$ by taking the $\T$-valued transition functions for
    $E$, and viewing them as $\operatorname{GL}_{1}(\C)$-valued
    transition functions for a complex vector bundle.  Thus $\B$ is
    the quotient of $\coprod U_{i}\times \C$ where $(i,x,\lambda)\sim
    (j,x,\sigma_{ij}(x)\lambda)$.  Since the $\sigma_{ij}(x)$ act by
    multiplication, we get a vector bundle rather than a principal
    bundle with respect to the additive group action on the fibres.}
  The multiplication on $\B^{(2)}$ goes as follows. We can view $\B$
  as the quotient of $E\times\C$ by the action of $\T$ given by
  $z(e,\lambda):=(ze,\bar z\lambda)$.  Then the product is just
  $[e,\lambda][f,\beta]:= [ef,\lambda\beta]$.

  \begin{remark}
    \label{rem-sec-twist}
    Note that we can view $E\subset \B$ (in the model above, just send
    $e$ to $[e,1]$).  Furthermore, sections of $\B$ are identified in
    a natural way with continuous $\C$-valued functions $\check f$ on
    $E$ which transform as follows:
    \begin{equation}
      \label{eq:7}
      \check f(z e)=\bar z \check f(e).
    \end{equation}
    The section $f\in\sa_{\relax}(G;\B)$ associated to $\check f$,
    transforming as in \eqref{eq:7}, is given by
    \begin{equation}
      \label{eq:8}
      f\bigl(j(e)\bigr)=\check f(e)e:=[e,\check f(e)].
    \end{equation}
  \end{remark}
\end{example}

\begin{example}[Green-Renault]
  \label{ex-renault}
  As pointed out in \cite{muh:cm01}*{Example~7(3)}, Examples
  \ref{ex-cp}~and \ref{ex-twist} are subsumed by Renault's formalism
  from \citelist{\cite{ren:jot91}\cite{ren:jot87}}.  In this case, we
  have a groupoid extension
  \begin{equation*}
    \xymatrix{\go\ar[r]&S\ar[r]^{i}&\Sigma\ar[r]^{j}&G\ar[r]&\go}
  \end{equation*}
  of locally compact groupoids over $\go$ where $S$ is a group bundle
  of \emph{abelian} groups with Haar system.  We view $S$ as a closed
  subgroupoid of $\Sigma$.

  In the spirit of Green twisted dynamical systems, we assume that we
  have a groupoid dynamical system\index{groupoid dynamical system}
  $(\A,\Sigma,\alpha)$ (so that $\pi:\A\to\go=\Sigma^{(0)}$ is an
  \usc{} \cs-bundle).  We also need an element $\chi\in\prod_{s\in
    S}M\bigl(A\bigl(r(s)\bigr)\bigr)$ such that
  \begin{equation}
    \label{eq:9}
    (s,a)\mapsto \chi(s)a
  \end{equation}
  is continuous from $S*\A:=\set{(s,a):r(s)=\pi(a)}$ to $\A$, and such
  that
  \begin{gather}
    \label{eq:10}
    \alpha_{s}(a)=\chi(s)a\chi(s)^{*}\quad\text{for all $(s,a)\in
      S*\A$, and} \\
    \chi(\sigma s\sigma^{-1})=\bar\alpha_{\sigma}\bigl(\chi(s)\bigr)
    \quad\text{ for $(\sigma,s)\in \Sigma^{(2)}$}.\label{eq:11}
  \end{gather}
  A little computation shows that
  \begin{equation*}
    a\chi(s)^{*}=\bigl(\chi(s)a^{*}\bigr)^{*},
  \end{equation*}
  and it follows that
  \begin{equation}
    \label{eq:12}
    (s,a)\mapsto a\chi(s)^{*}
  \end{equation}
  is continuous.  Therefore we can define a $S$-action on
  $r^{*}\A=\set{(a,\sigma):\pi(a)=r(\sigma)}$:
  \begin{equation}
    \label{eq:13}
    (a,\sigma)\cdot s:=(a\chi(s)^{*},s\sigma).
  \end{equation}
  We define $\B$ to be the quotient $r^{*}\A/S$, and define $p:\B\to
  G$ by $p\bigl([a,\sigma]\bigr)=j(\sigma)$.
  \begin{lemma}
    \label{lem-one}
    With the set-up above, $p:\B\to G$ is an \usc\ Banach bundle over
    $G$.
  \end{lemma}
  \begin{proof}
    The proof is obtained by modifying the proof of
    \cite{kmrw:ajm98}*{Proposition~2.15}.  Specifically, we have to
    show that $p$ is open, and that axioms B1--B4 of
    Definition~\ref{def-usc-bundle} are satisfied.\footnote{In
      \cite{kmrw:ajm98}*{Proposition~2.15}, we are working with
      \cs-bundles.  Here we have to adjust to \usc-bundles which are
      merely Banach bundles.}
  \end{proof}

  To get a Fell bundle, we'll need a multiplication on
  \begin{equation*}
    \B^{(2)}:=\set{\bigl([a,\sigma],
      [b,\tau]\bigr):\bigl(j(\sigma),j(\tau)\bigr) \in G^{(2)}}. 
  \end{equation*}
  Since $\bigl(j(\sigma),j(\tau)\bigr) \in G^{(2)}$ exactly when
  $(\sigma,\tau)\in\Sigma^{(2)}$, we want to try
  \begin{equation}
    \label{eq:14}
    [a,\sigma][b,\tau]:=[a\alpha_{\sigma}(b),\sigma\tau].
  \end{equation}
  To see that \eqref{eq:14} is well defined, consider
  \begin{align*}
    [a\chi(s)^{*},s\sigma][b\chi(t)^{*},t\tau]&= [a\chi(s)^{*}
    \alpha_{s\sigma}\bigl(b\chi(t)^{*}\bigr),s\sigma t\tau] \\
    \intertext{which, using \eqref{eq:10}, is}
    &=[a\alpha_{\sigma}\bigl(b\chi(t)^{*}\bigr)\chi(s)^{*},s\sigma
    t\sigma^{-1}\sigma \tau] \\
    \intertext{which, using \eqref{eq:11}, is}
    &=[a\sigma_{\sigma}(b)\chi(\sigma
    t\sigma^{-1})^{*}\chi(s)^{*},(s\sigma t\sigma^{-1})\sigma\tau]\\
    &=[a\alpha_{\sigma}(b)\chi(s\sigma t\sigma^{-1})^{*}, (s\sigma
    t\sigma^{-1})\sigma\tau] \\
    &=[a\alpha_{\sigma}(b),\sigma\tau].
  \end{align*}
  Thus \eqref{eq:14} is well-defined and we can establish the
  following lemma without difficulty.
  \begin{lemma}
    \label{lem-two}
    With respect to the multiplication defined above, $p:\B\to G$ is a
    Fell bundle over $G$.
  \end{lemma}

  To get a section of $\B$, we need a continuous function $f:\Sigma\to
  \A$ such that
  \begin{gather}
    \label{eq:15}
    f(\sigma)\in A\bigl(r(\sigma)\bigr)\quad\text{and}\\
    \label{eq:16} f(s\sigma)=f(\sigma)\chi(s)^{*}\quad\text{for
      $(s,\sigma)\in\Sigma^{(2)}$.}
  \end{gather}
  The corresponding section is given by
  \begin{equation}
    \label{eq:17}
    \check f\bigl(j(\sigma)\bigr):=[f(\sigma),\sigma].
  \end{equation}
\end{example}

Now recall that the set $\sa_{c}(G;\B)$ of continuous compactly
supported sections is endowed with a $*$-algebra structure as follows:
\begin{gather}
  \label{eq:18}
  f*g(x):=\int_{G} f(y)g(y^{-1}x)\,d\lambda^{r(x)}(y),
  \quad\text{and}\\
  f^{*}(x):=f(x^{-1})^{*}.\label{eq:19}
\end{gather}

It is a worthwhile exercise to look a bit more closely at the
$*$-algebra $\sa_{c}(G;\B)$ in each of the basic examples above.

\begin{example}[The crossed product for Example~\ref{ex-cp}]
  \label{ex-cp-cp}
  Let $p:\B\to G$ be the Fell bundle associated to a dynamical
  system\index{groupoid dynamical system}\index{groupoid crossed product}
  $(\A,G,\alpha)$ as in Example~\ref{ex-cp}.  Let
  $f,g\in\sa_{c}(G;\B)$, and let $\check f$ and $\check g$ be the
  corresponding $\A$-valued functions on $G$ as in
  Remark~\ref{rem-cp}.  Then
  \begin{align*}
    f*g(x)&=\int_{G} (\check f(y),y)(g(y^{-1}x),y^{-1}x)
    \,d\lambda^{r(x)}(y) \\
    &=\int_{G} \bigl(\check f(y)\alpha_{y}\bigl(\check
    g(y^{-1}x)\bigr),x\bigr) \,\lambda^{r(x)} (y).
  \end{align*}
  Similarly,
  \begin{equation*}
    f^{*}(x)=f(x^{-1})^{*}=(\check
    f(x^{-1}),x^{-1})^{*} =\bigl(\alpha_{x}\bigl(f(x^{-1})\bigr)^{*},x\bigr)
  \end{equation*}
  Thus if we confound $f$ and $\check f$, as is usually the case, we
  obtain the ``usual'' convolution formula:
  \begin{equation}
    \label{eq:20}
    f*g(x):=\int_{G}f(y)\alpha_{y}\bigl(g(y^{-1}x)\bigr) \, d\lambda^{r(x)}(y),  \tag{\ref{eq:18}$'$}
  \end{equation}
  and the ``usual'' involution formula:
  \begin{equation}
    \label{eq:21}
    f^{*}(x)=\alpha_{x}\bigl(f(x^{-1})^{*}\bigr).
    \tag{\ref{eq:19}$'$} 
  \end{equation}
  In this case, after completing as in
  Section~\ref{sec:repr-fell-bundl}, we obtain the crossed product
  $\A\rtimes_{\alpha}G$ (or
  $A\rtimes_{\alpha}G$).\index{acrossg@$\A\rtimes_{\alpha}G$}   (For more on
  groupoid crossed products, see \cite{muhwil:nyjm08}.)
\end{example}

\begin{example}[The crossed product for Example~\ref{ex-twist}]
  \label{ex-twist-cp}
  In this case, we work with functions on $E$ transforming as in
  \eqref{eq:7}.  Then
  \begin{equation*}
    f*g(e)=\int_{G} [\check f(d)\check
    g(d^{-1}e),e]\,d\lambda^{r(e)}(j(d)) ,
  \end{equation*}
  and $f*g$ is represented by the function on $E$ given by
  \begin{equation*}
    e\mapsto \int_{G} \check f(d)\check
    g(d^{-1}e) \,d\lambda^{r(j(e))}(j(d)).
  \end{equation*}
  Thus the completion is the algebra $\cs(G;E)$\index{csge@$\cs(G;E)$} as in
  \cite{muhwil:plms395}.
\end{example}

\begin{example}
  [The crossed product for Example~\ref{ex-renault}] Here we
  replace\index{groupoid dynamical system}\index{groupoid crossed product}
  $f$ with $\check f$.  Then the functions transform as in
  \eqref{eq:16}, and the operations on functions on $\Sigma$ are given
  by
  \begin{gather*}
    f*g(\sigma)=\int_{G}f(\tau)\alpha_{\tau}(
    g(\tau^{-1}\sigma)\,d\lambda^{r(j(\sigma))}(\tau) \quad\text{and}\\
    f^{*}(\sigma)=\alpha_{\sigma}\bigl(f(\sigma^{-1})^{*}\bigr).
  \end{gather*}
  The completion is Renault's
  $\cs(G,\Sigma,\A,\lambda)$\index{csgsal@$\cs(G,\Sigma,\A,\lambda)$}
  from 
  \citelist{\cite{ren:jot91}\cite{ren:jot87}}.
\end{example}

\section{Generalized Radon Measures on Fell Bundles}
\label{sec:radon-measures-fell}

For the proof of the disintegration theorem for representations of
Fell bundles (Theorem~\ref{thm-fell-disintegration}), we will need a
version of Yamagami's \cite{yam:xx87}*{Lemma~2.2} suitable for
\usc-Banach bundles.  Note that \cite{yam:xx87}*{Lemma~2.2} is
intended to be a restatement of
\cite{din:integration74}*{Theorem~28.32} due to Dinculeanu.  In fact,
there is a bit of work to do just to coax out the result Yamagami
claims in the (continuous) Banach bundle case from
\cite{din:integration74}*{Theorem~28.32}.  Therefore it seems more
than reasonable to work out the details of the more general result
here.

First we need some terminology and notation. Let $p:\B\to G$ be an
\usc-Banach bundle over a \emph{second countable} Hausdorff groupoid
$G$ such that the corresponding Banach space $B:=\sa_{0}(G;\B)$ is
\emph{separable}.  We call a linear functional
\begin{equation*}
  \nu:\gcb\to\C
\end{equation*}
a \emph{\grm{}}\index{generalized Radon measure@\grm} provided that $\nu$ is continuous in the inductive
limit topology.  Of course, if $\B$ is the trivial (complex) line
bundle over $G$, then a \grm{} is just a complex Radon
measure.\index{complex Radon measure}  Some
useful comments on complex Radon measures can be found in
\cite{muhwil:nyjm08}*{Appendix~A.1}.

\begin{lemma}
  \label{lem-barnu}
  Let $\nu:\gcb\to\C$ be a \grm{} on $G$.  Then there is a Radon
  measure $\mu$ on $G$ such that for all $\phi\in C_{c}^{+}(G)$,
  \begin{equation}\label{eq:22}
    \mu(\phi):=\sup\set{|\nu(f)|:\|f\|\le\phi}.
  \end{equation}
\end{lemma}
\begin{proof}
  We will produce a function $\mu:\ccpg\to \R^{+}$ satisfying
  \eqref{eq:22} and such that for all $\alpha\ge0$ and
  $\phi_{i}\in\ccpg$ we have
  \begin{enumerate}
  \item $\mu(\alpha \phi_{1})=\alpha\mu(\phi_{1})$ and
  \item $\mu(\phi_{1}+\phi_{2})=\mu(\phi_{1})+\mu(\phi_{2})$.
  \end{enumerate}
  Then it is not hard to see that we can extend $\mu$ to all of $\ccg$
  in the expected way (cf.,
  \cite{din:integration74}*{Proposition~2.20} or
  \cite{hr:abstract}*{Theorem~B.38}).

  Naturally, we define $\mu$ on $\ccpg$ using \eqref{eq:22}.  Then
  part~\partref1 follows immediately from the definition of $\mu$.
  Note that if $K\subset G$ is compact, then the continuity of $\nu$
  implies that there is a constant $a_{K}\ge0$ such that
  \begin{equation*}
    |\nu(f)|\le a_{K}\|f\|_{\infty}\quad\text{for all $f\in \gcb$ such
      that $\supp f\subset K$.}
  \end{equation*}
  (If not, then for each $n$ we can find $f_{n}$ with $\supp
  f_{n}\subset K$, $\|f_{n}\|_{\infty}\le 1$ and $|\nu(f_{n})|\ge
  n^{2}$.  Then $\frac1n f_{n}$ tends to $0$ in the inductive limit
  topology, while $\nu(\frac1n f_{n})\not\to0$.)  It follows that
  $\mu(\phi)<\infty$ for all $\phi\in\ccpg$.

  Fix $\phi_{1},\phi_{2}\in\ccpg$ and $\epsilon>0$.  Choose
  $f_{i}\in\gcb$ such that $\|f_{i}\|\le\phi_{i}$ and
  \begin{equation*}
    \mu(\phi_{i})\le|\nu(f_{i})|+\frac\epsilon2.
  \end{equation*}
  Let $\tau_{i}$ be a unimodular scalar such that
  $|\nu(f_{i})|=\nu(\tau_{i}f)$.  Then
  \begin{equation*}
    \|\tau_{1}f_{1}+\tau_{2}f_{2}\|\le \|f_{1}\|+\|f_{2}\|\le
    \phi_{1}+\phi_{2}. 
  \end{equation*}
  Thus
  \begin{align*}
    \mu(\phi_{1})+\mu(\phi_{2})&\le |\nu(f_{1})|+|\nu(f_{2})|+\epsilon
    \\
    &= \nu(\tau_{1}f_{1}+\tau_{2}f_{2})+\epsilon\\
    &=\mu(\phi_{1}+\phi_{2})+\epsilon.
  \end{align*}
  Since $\epsilon$ was arbitrary, $\mu(\phi_{1})+\mu(\phi_{2})\le
  \mu(\phi_{1} +\phi_{2})$.

  Now suppose that $h\in \gcb$ is such that
  $\|h\|\le\phi_{1}+\phi_{2}$.  Define
  \begin{equation*}
    \sigma_{1}(x):=
    \begin{cases}
      \displaystyle{ \frac{\phi_{1}(x)}{\phi_{1}(x)+\phi_{2}(x)}}&
      \text{if $\phi_{1}(x) +\phi_{2}(x)>0$, and } \\
      0&\text{otherwise.}
    \end{cases}
  \end{equation*}
  Let $h_{i}:=\sigma_{i}\cdot h$.  We want to see that each
  $h_{i}\in\gcb$.  For this, we just need to see that $x\to h_{i}(x)$
  is continuous from $G$ to $\B$.  Fix $x_{0}\in G$.  If
  $h(x_{0})\not=0$, then $\phi_{1}(x_{0})+\phi_{2}(x_{0})>0$ and
  $\phi_{1}(x)+\phi_{2}(x) >0$ near $x_{0}$.  Consequently,
  $\sigma_{i}$ is continuous near $x_{0}$.  Therefore $h_{i}$ is
  continuous at $x_{0}$ (see \cite{wil:crossed}*{Proposition~C.17}).
  On the other hand, suppose that $h(x_{0})=0$ and $\epsilon>0$.
  Since $x\mapsto \|h(x)\|$ is upper semicontinuous, there is a
  neighborhood $V$ of $x_{0}$ such that $\|h(x)\|<\epsilon$ for all
  $x\in V$.  However, $|\sigma_{i}(x)|\le1$ for all $x$.  Hence
  \begin{equation*}
    \|h_{i}(x)\|=|\sigma_{i}(x)|\|h(x)\|\le
    \|h(x)\|<\epsilon\quad\text{provided $x\in V$.}
  \end{equation*}
  It follows that $h_{i}$ is continuous at $x_{0}$ (for example, see
  axiom~B4 of Definition~A.1).

  Furthermore, if $h_{i}(x)\not=0$, then $\phi_{1}(x)+\phi_{2}(x)>0$
  and
  \begin{equation*}
    \|h_{i}(x)\|= \phi_{i}(x)
    \frac1{\phi_{1}(x)+\phi_{2}(x)} \|h(x)\| \le \phi_{i}(x),
  \end{equation*}
  since $\|h(x)\|\le \phi_{i}(x)+\phi_{2}(x)$.  Thus if $\|h\|\le
  \phi_{1}+\phi_{2}$, then there are $h_{i}$ such that $h=h_{1}+h_{2}$
  and $\|h_{i}\|\le\phi_{i}$.  Therefore we can compute as follows:
  \begin{align*}
    \mu(\phi_{1}+\phi_{2})&=\sup\set{|\nu(h)|:\|h\|\le
      \phi_{1}+\phi_{2}} \\
    &= \sup\set{|\nu(h_{1}+h_{2})|:\|h_{i}\|\le\phi_{i}} \\
    &\le \sup\set{|\nu(h_{1})|+|\nu(h_{2})|: \|h_{i}\|\le\phi_{i}} \\
    &\le \mu(\phi_{1})+\mu(\phi_{2}).
  \end{align*}
  This establishes \partref2.  Since $G$ is second countable and since
  we established that $\mu$ was finite on $C_{c}(G)$, $\mu$ is a Radon
  measure (by, for example, \cite{rud:real}*{Theorem~2.18}).
\end{proof}

\begin{example}
  \label{ex-tot-var}
  Suppose that $\nu:\ccg\to\C$ is a \emph{complex} Radon measure on
  $G$.  Then the measure $\mu$ associated to $\nu$ via \eqref{eq:22}
  in Lemma~\ref{lem-barnu} is the total variation $|\nu|$ of $\nu$.
\end{example}
\begin{proof}
  It suffices to see that for all $\phi\in\ccpg$ we have
  \begin{equation*}
    |\nu|(\phi)=\sup\set{|\nu(f)|:\|f\|\le\phi}.
  \end{equation*}
  Since $\nu=\tau|\nu|$ for a unimodular function $\tau$ (see
  \cite{fol:real}*{Proposition~3.13} in the bounded case and
  \cite{muhwil:nyjm08}*{Appendix~A.1} in general) and since we can let
  $f=\bar\tau\phi$, we clearly get equality.
\end{proof}

In view of Example~\ref{ex-tot-var}, we call the measure
$\mu$\index{total variation measure}\index{generalized Radon
  measure@\grm!total variation of}
appearing in Lemma~\ref{lem-barnu} the \emph{total variation} of the
\grm{} $\nu$ and write $|\nu|$ in place of $\mu$.\index{barnu@$"|\nu"|$}  Then what we need
to prove is the following.
\begin{prop}
  \label{prop-din28.32}
  Suppose that $p:\B\to G$ is an \usc-Banach bundle over a second
  countable locally compact Hausdorff space $G$ such that the section
  algebra $B:=\sa_{0}(G;\B)$ is separable.\footnote{This hypothesis of
    separability seems to be crucial.  In the proof we will have to
    collect null sets for a countable dense subset reminiscent of
    direct integral type arguments.}  If $\nu:\gcb\to\C$ is a \grm{}
  on $\B$ with total variation $|\nu|$, then for all $x\in G$ there
  are linear functionals $\epsilon_{x}\in B(x)^{*}$ of norm at
  most\index{epsilonsubx@$\epsilon_{x}$} 
  one such that
  \begin{enumerate}
  \item for each $f\in\gcb$, $x\mapsto \epsilon_{x}\bigl(f(x)\bigr) $
    is in $\bb_{c}(G)$\footnote{We use $\bb_{c}(G)$ to denote the
      bounded Borel functions on $G$ which vanish off a compact set.}
    and
  \item such that
    \begin{equation*}
      \nu(f)=\int_{G} \epsilon_{x}\bigl(f(x)\bigr) \,
      d|\nu|(x) 
    \end{equation*}
    for all $f\in\gcb$.
  \end{enumerate}
\end{prop}

Before proceeding with the proof, we need to deal with the reality that
--- unlike the case in \cite{din:integration74} where continuous
Banach bundles are used --- $\|f\|$ need not be in $\ccg$ if
$f\in\gcb$.  But it is at least upper semicontinuous (by
axiom~B1 of Definition~\ref{def-usc-bundle}).  Therefore $x\mapsto
\|f(x)\|$ is in $\bb_{c}(G)$.  In particular, it is integrable with
respect to any Radon measure on $G$.

We need the following observations.

\begin{lemma}
  \label{lem-usc}
  Suppose that $f$ is a bounded nonnegative upper semicontinuous
  function with compact support on $G$.  Then if $\mu$ is a Radon
  measure on $G$,
  \begin{equation}
    \label{eq:23}
    \int_{G}f(x)\,d\mu=\inf\Bigl\{\, \int_{G} g(x)\,
    d\mu(x): \text{$g\in C_{c}^{+}(G)$ and $f\le g$}\Bigr\}.
  \end{equation}
\end{lemma}
\begin{proof}
  Fix $\phi\in\ccpg$ such that $f\le \phi$.  Then $\phi-f$ is a
  nonnegative \emph{lower} semicontinuous function on $G$.  By
  \cite{fol:real}*{Corollary~7.13}, given $\epsilon>0$, there is a
  $g\in\ccg$ such that $0\le g\le \phi-f$ and such that
  \begin{equation*}
    \int_{G} \bigl(\phi(x)-f(x)\bigr)\,d\mu(x) \le
    \int_{G}g(x) \,d\mu(x)+\epsilon.
  \end{equation*}
  But then
  \begin{equation*}
    \int_{G}\bigl(\phi(x)-g(x)\bigr) \le \int_{G}
    f(x)\,d\mu(x) +\epsilon.
  \end{equation*}
  Since $f\le\phi-g\le\phi$, we have $\phi-g\in\ccg$ and dominates
  $f$.  Since $\epsilon>0$ was arbitrary, the right-hand side of
  \eqref{eq:23} is at least the left-hand side.  The other inequality
  is clear, so the result is proved.
\end{proof}

\begin{cor}
  \label{cor-inequal}
  Suppose that $\nu:\gcb\to\C$ is a \grm{} on $G$.  Then for all
  $f\in\gcb$,
  \begin{equation*}
    |\nu(f)|\le|\nu|\bigl(\|f\|\bigr).
  \end{equation*}
  In particular, if $\mu$ is any Radon measure on $G$ such that
  $|\nu(f)|\le \mu\bigl(\|f\|\bigr)$ for all $f\in\gcb$, then
  $|\nu|(\phi)\le \mu(\phi)$ for all $\phi\in\ccpg$.
\end{cor}
\begin{proof}
  Suppose that $\phi\in\ccpg$ and $\|f\|\le\phi$.  Then by definition,
  $|\nu| (\phi)\ge|\nu(f)|$.  By the previous lemma,
  \begin{align*}
    |\nu|\bigl(\|f\|\bigr) &=
    \inf\set{|\nu|(\phi):\text{$\phi\in\ccpg$
        and $\|f\|\le\phi$}} \\
    &\ge|\nu(f)|.
  \end{align*}

  Now suppose that $\mu$ is as above.  Then for all $\phi\in\ccpg$,
  \begin{align*}
    |\nu|(\phi)&=\sup\set{|\nu(f)|:\|f\|\le\phi}\\
    &\le \sup \set{\mu(\|f\|):\|f\|\le\phi} \\
    &\le \mu(\phi)\qedhere
  \end{align*}
\end{proof}

\begin{proof}[Proof of Proposition~\ref{prop-din28.32}]
  Suppose $\nu:\gcb\to\C$ is a \grm{} as in the statement of the
  proposition.  For each $f\in\gcb$, define an (scalar-valued) complex
  Radon measure on $G$ by $\nu_{f}(\phi):=\nu(\phi\cdot f)$.  Then by
  Corollary~\ref{cor-inequal},
  \begin{equation}
    \label{eq:24}
    |\nu_{f}(\phi)|=|\nu(\phi\cdot f)|\le\nnu\bigl(|\phi|
    \|f\|\bigr). 
  \end{equation}
  Therefore
  \begin{equation}
    \label{eq:25}
    |\nu_{f}(\phi)|\le \|f\|_{\infty}\|\phi\|_{1},
  \end{equation}
  where $\|\phi\|_{1}=\nnu(|\phi|)$ is the norm of $\phi$ in
  $L^{1}(G,\nnu)$.

  It follows that $f\mapsto \nu_{f}$ is a bounded linear map $\Phi$ of
  $\gcb\subset B=\sa_{0}(G;\B)$ into $L^{1}(\nnu)^{*}$.  Of course, we
  can identify $L^{1}(\nnu)^{*}$ with $L^{\infty}(\nnu)$.  Then, since
  $B$ is separable and $\Phi$ is bounded, we can identify
  $S:=\Phi\bigl(\gcb\bigr)$ with a separable subspace of
  $L^{\infty}(\nnu)$.  By \cite{wil:crossed}*{Lemma~I.8}, there is a
  linear map $\rho:S\to \bb(G)$ such that the $\nnu$-almost everywhere
  equivalence class of $\rho(b)$ is $b$.\footnote{The function $\rho$
    is called a ``lift'' for $S\subset L^{\infty}(\nnu)$.  More
    comments and references about lifts can be found in the paragraph
    preceding \cite{wil:crossed}*{Lemma~I.8}.}

  Thus if we let $b_{f}:= \rho(\nu_{f})$, then $b_{f}$ is a bounded
  Borel function such that
  \begin{equation*}
    \nu_{f}(\phi) =\int_{G}\phi(x)b_{f}(x)\,d\nnu(x).
  \end{equation*}
  Furthermore, the linearity of $\rho$ implies that
  \begin{equation}
    \label{eq:26}
    b_{\alpha f + \beta g}=\alpha b_{f}+\beta b_{g}\quad\text{for all
      $f,g\in\gcb$ and $\alpha,\beta\in\C$.}
  \end{equation}
  It follows from Corollary~\ref{cor-inequal} (as well as
  Example~\ref{ex-tot-var}) and Equation~\eqref{eq:24} that for all
  $\phi\in\ccpg$,
  \begin{equation}
    \label{eq:27}
    |\nu_{f}|(\phi)\le \nnu\bigl(\phi \|f\|\bigr).
  \end{equation}
  Since $|\nu_{f}|=|b_{f}|\nnu$, \eqref{eq:27} amounts to
  \begin{equation*}
    \int_{G}\phi(x)|b_{f}(x)|\,d\nnu(x)\le\int_{G}\phi(x) \|f(x)\| \,d\nnu(x)\quad\text{for all $\phi\in\ccpg$.}
  \end{equation*}
  Therefore,
  \begin{equation*}
    |b_{f}(x)|\le \|f(x)\|\quad\text{for $\nnu$-almost all
      $x$.} 
  \end{equation*}
  Using \eqref{eq:26}, this means that
  \begin{equation*}
    |b_{f}(x)-b_{g}(x)|\le\|f(x)-g(x)\|\quad\text{for 
      $\nnu$-almost all $x$.} 
  \end{equation*}

  Since $B$ is separable, there is a sequence $\set{f_{n}}\subset
  \gcb$ such that $\set{f_{n}(x)}$ is dense in $B(x)$ for all $x$.
  Let $B_{0}$ be the rational span of the $f_{n}$.  Then
  \begin{equation*}
    B_{0}(x):=\set{f(x):f\in B_{0}}
  \end{equation*}
  is a vector space over $\mathbf{Q}$ which is dense in $B(x)$.  Since
  $B_{0}$ is countable, there is a $\nnu$-null set $N$ such that for
  all $x\notin N$ and all $f,g\in B_{0}$,
  \begin{align}
    \label{eq:28}
    |b_{f}(x)|&\le\|f(x)\|\quad\text{and}\\
    |b_{f}(x)-b_{g}(x)|&\le \|f(x)-g(x)\|.
    \label{eq:29}
  \end{align}
  We can simply alter each $b_{f}$ so that $b_{f}(x)=0$ if $x\in N$.
  Then, for all $f,g\in B_{0}$, equations \eqref{eq:26} (for $\alpha$
  and $\beta$ rational), \eqref{eq:28} and \eqref{eq:29} are valid for
  \emph{all} $x\in G$.

  In particular, \eqref{eq:29} implies that if $f,g\in B_{0}$ and
  $f(x)=g(x)$, then $b_{f}(x)=b_{g}(x)$.  Therefore we get a
  well-defined $\mathbf{Q}$-linear map
  \begin{equation*}
    \epsilon_{x}:B_{0}(x)\to\C
  \end{equation*}
  by letting $\epsilon_{x}(a)=b_{f}(x)$, where $f$ is any section in
  $B_{0}$ such that $f(x)=a$.  In view of \eqref{eq:28}, each
  $\epsilon_{x}$ has norm at most one.  It follows that $\epsilon_{x}$
  extends uniquely to an element in $B(x)^{*}$ of norm at most one
  which we continue to denote by $\epsilon_{x}$.

  If $f\in B_{0}$, then $b_{f}(x)=\epsilon_{x}\bigl(f(x)\bigr)$.
  Consequently, $x\mapsto \epsilon_{x}\bigl(f(x)\bigr)$ is in
  $\bb_{c}(G)$.  Furthermore, if $f\in B_{0}$ and $\phi\in \ccg$, then
  \begin{align*}
    \nu(\phi\cdot f)&= \nu_{f}(\phi) \\
    &= \int_{G}\phi(x) b_{f}(x)\,d\nnu(x) \\
    &= \int_{G} \phi(x) \epsilon_{x}\bigl(f(x)\bigr ) \,
    d\nnu(x) \\
    &= \int_{G} \epsilon_{x}\bigl(\phi\cdot f(x)\bigr)\, d\nnu(x) .
  \end{align*}
  But, if $f\in\gcb$, then there is a sequence $\set{g_{k}}\subset
  \gcb$ converging to $f$ in the inductive limit topology such that
  each $g_{k}$ is a finite sum of the form $\sum \phi_{j}\cdot f_{j}$
  with each $\phi_{j}\in \ccg$ and each $f_{j}\in B_{0}$.  By the
  above, we have
  \begin{equation*}
    \nu(g_{k})=\int_{G}\epsilon_{x}\bigl(g_{k}(x)\bigr)\,
    d\nnu(x), 
  \end{equation*}
  and since $g_{k}\to f$ in the inductive limit topology and since
  each $\epsilon_{x}$ has norm at most one,
  $\epsilon_{x}\bigl(g_{k}(x)\bigr) \to \epsilon_{x}\bigl(f(x)\bigr)$
  uniformly and the entire sequence vanishes off a compact set.  Hence
  $x\mapsto \epsilon_{x}\bigl(f(x)\bigr)$ is in $\bb_{c}(G)$ as
  claimed in part~\partref1 of the proposition, and
  \begin{align*}
    \nu(f)&=\lim_{k}\nu(g_{k}) \\
    &=\lim_{k} \int_{G} \epsilon_{x}\bigl(g_{k}(x)\bigr) \,
    d\nnu(x) \\
    &= \int_{G} \epsilon_{x}\bigl( f(x)\bigr) \, d\nnu(x).
  \end{align*}
  This completes the proof of the proposition.
\end{proof}

%%
%% Some grm stuff
%%
We close this section with some technicalities that will be needed
later.  In particular, we will need to deal with some not necessarily
continuous sections.\footnote{Note that we are trying to avoid dealing
  with measurability issues for sections of $\B$.  Consequently, we
  are only going to introduce only those potentially discontinuous
  sections that we absolutely have to.}

\begin{definition}
  \label{def-sigma1}
  We let $\sgcb$\index{sigmagcb@$\sgcb$} be the set of bounded
  sections $f$ of $p:\B\to G$ such that there is a uniformly bounded
  sequence $\set{f_{n}}\subset \gcb$ and a compact set $K\subset G$
  such that $\supp f_{n}\subset K$ for all $n$ and such that
  $f_{n}(x)\to f(x)$ for all $x\in G$.\footnote{Recall that we take
    the point of view that a section of $p:\B\to G$ is simply a
    function $f$ from $G$ to $\B$ such that $p(f(x))=x$.}
  Analogously, we let $\bboc(G)$\index{Bonec@$\bboc$} the family of
  Borel functions $\phi$ on $G$ such that there is a uniformly bounded
  sequence $\set{\phi_{n}}\subset \ccg$ and a compact set $K\subset G$
  such that $\supp \phi_{n}\subset K$ for all $n$ and such that
  $\phi_{n}(x)\to \phi(x)$ for all $x\in G$.
\end{definition}

\begin{example}
  \label{ex-pauls-b1}
  Suppose that $f\in\gcb$ and $\phi\in\bboc(G)$.  Then by considering
  $\set{\phi_{n}\cdot f}$ for appropriate $\phi_{n}$, we see that
  $\phi\cdot f\in\sgcb$.
\end{example}

\begin{lemma}
  \label{lem-tol-var-borel}
  Suppose that $\sigma$ is a \grm{} on $\B$ and that $\phi\in\bboc(G)$
  is such that there are $\set{\phi_{n}}\subset\ccpg$ with
  $\phi_{n}(x)\searrow \phi(x)$ for all $x\in G$.  Then
  \begin{equation*}
    \nsigma(\phi)=\sup\set{|\sigma(f)|: 
      \text{$f\in\sgcb$ and $\|f\|\le\sigma$.}}
  \end{equation*}
\end{lemma}

First, an even more specialized result.

\begin{lemma}
  \label{lem-split}
  Suppose that $\phi\in\bboc(G)$ is such that there is a sequence
  $\set{\phi_{n}}$ in $\ccpg$ such that $\phi_{n}\searrow \phi$.  If
  $f\in\sgcb$ is such that $\|f\|\le\phi_{n}$, then there are
  $f_{1},f_{2}\in\sgcb$ such that $f=f_{1}+f_{2}$, $\|f_{i}\|\le \phi$
  and $\|f_{2}\|\le\phi_{n}-\phi$.
\end{lemma}
\begin{proof}
  Define
  \begin{equation*}
    \sigma(x):=
    \begin{cases}
      \frac{\phi(x)}{\phi_{n}(x)}&\text{if $\phi_{n}(x)>0$ and}\\
      0&\text{otherwise}
    \end{cases}\quad \tau(x):=
    \begin{cases}
      \frac{\phi_{n}(x)-\phi(x)}{\phi_{n}(x)}&\text{if $\phi_{n}(x)>0$
        and} \\
      0&\text{otherwise}.
    \end{cases}
  \end{equation*}
  Then we clearly have $f=\sigma\cdot f+ \tau\cdot f$, $\|\sigma\cdot
  f\|\le\phi$ and $\|\tau\cdot f\|\le \phi_{n}-\phi$.  Therefore we
  just need to see that $\sigma\cdot f$ and $\tau\cdot f$ are in
  $\sgcb$.  But $\sigma(x)=\lim_{m\ge n}\sigma_{m}(x)$ where
  \begin{equation*}
    \sigma_{m}(x):=
    \begin{cases}
      \frac{\phi_{m}}{\phi_{n}}&\text{if $\phi_{n}(x)>0$ and}\\
      0&\text{otherwise.}
    \end{cases}
  \end{equation*}
  As in the proof of Lemma~\ref{lem-barnu}, we have $\sigma_{m}\cdot
  f\in\gcb$ and $\sigma_{m}\cdot f\to \sigma\cdot f$ pointwise.  It
  follows that $\sigma\cdot f\in\sgcb$.  A similar argument shows that
  $\tau\cdot f\in\sgcb$.  This is what we needed to show.
\end{proof}

\begin{proof} [Proof of Lemma~\ref{lem-tol-var-borel}]
  As in Proposition~\ref{prop-din28.32}, we have
  \begin{equation*}
    \sigma(f)=\int_{G}\epsilon_{x}\bigl(f(x)\bigr)\,d\nsigma(x)
  \end{equation*}
  for linear functionals $\epsilon_{x}$ in the unit ball of
  $B(x)^{*}$.  Thus if $\|f\|\le\phi$, then
  \begin{equation*}
    |\sigma(f)|\le \int_{G}\|f(x)\|\,d\nsigma(x) \le \nsigma(\phi).
  \end{equation*}
  By the dominated convergence theorem, $\nsigma(\phi_{n})\to
  \nsigma(\phi)$.  Therefore, there are $f_{n}\in\gcb$ such that
  $\|f_{n}\|\le \phi_{n}$ and such that $|\sigma(f_{n})| \to
  \nsigma(\phi)$.  By Lemma~\ref{lem-split}, we can decompose
  $f_{n}=f_{n}'+f_{n}''$ in $\sgcb$ such that $\|f_{n}'\|\le\phi$ and
  $\|f_{n}''\|\le \phi_{n}-\phi$.  By the above,
  \begin{equation*}
    |\sigma(f_{n}'')|\le\nsigma(\phi_{n}-\phi),
  \end{equation*}
  and $\sigma(f_{n}'')\to0$ (as $n\to\infty$).  It then follows that
  $|\sigma(f_{n}')|\to \nsigma(\phi)$.  The result follows.
\end{proof}

Part~\partref1 of Proposition~\ref{prop-din28.32} gives us ``just
enough'' measurability of the $\epsilon_{x}$ to get by.  We will need
to amplify this a bit for the proof of the disintegration theorem in the
next section.  What we need is provided below.

Let $m:G^{(2)}\to G$ be the multiplication map and let $m^{*}\B$ be
the pull back.

\begin{lemma}
  \label{lem-prod-borel}
  Suppose that $\nu:\gcb\to\C$ is a \grm{} given by
  \begin{equation}\label{eq:30}
    \nu(f)=\int_{G} \epsilon_{x}(f)\,d\nnu(x)
  \end{equation}
  as in Proposition~\ref{prop-din28.32}.  If
  $F\in\sa_{c}(G^{(2)};m^{*}\B)$, then
  \begin{equation*}
    (x,y)\mapsto \epsilon_{xy}\bigl(F(x,y)\bigr)
  \end{equation*}
  is a Borel function on $G^{(2)}$.
\end{lemma}
\begin{proof}
  Let $\AA$ be the subalgebra of $\sa_{c}(G^{(2)};m^{*}\B)$ spanned by
  sections of the form
  \begin{equation*}
    (x,y)\mapsto \phi(x,y)f(xy)
  \end{equation*}
  with $\phi\in C_{c}(G^{(2)})$ and $f\in\gcb$.  Clearly, $\AA$ is
  closed under multiplication by functions from $C_{c}(G^{(2)})$ and
  $\set{F(x,y):F\in\AA}$ is dense in $m^{*}\B_{(x,y)}=\B_{xy}$ for all
  $(x,y)$.  Then a partition of unity argument (see
  Lemma~\ref{lem-dense-sections}) implies that $\AA$ is dense in
  $\sa_{c}(G^{(2)};m^{*}\B)$ in the \ilt.
  Proposition~\ref{prop-din28.32} implies that for each $f\in\gcb$,
  $x\mapsto \epsilon_{x}\bigl(f(x)\bigr)$ is a bounded Borel function
  with compact support.  Since
  \begin{equation*}
    (x,y)\mapsto \epsilon_{xy}\bigl(f(xy)\bigr)
  \end{equation*}
  is the composition of two Borel functions, it too is Borel.  Thus
  \begin{equation*}
    (x,y)\mapsto \phi(x,y)\epsilon_{xy}\bigl(f(xy)\bigr)
  \end{equation*}
  is Borel for all $\phi\in C_{c}(G^{(2)}))$ and
  $f\in\gcb$. Consequently, $(x,y)\mapsto
  \epsilon_{xy}\bigl(F(x,y)\bigr)$ is Borel for all $F\in\AA$.  But if
  $F_{i}\to F$ in the \ilt{} on $\gcb$, then the functions
  $(x,y)\mapsto \epsilon_{xy}\bigl(F_{i}(x,y)\bigr)$ converge
  uniformly to $(x,y)\mapsto \epsilon_{xy} \bigl(F(x,y)\bigr)$.
  Therefore the result follows as $\AA$ is dense in
  $\sa_{c}(G^{(2)};m^{*}\B)$.
\end{proof}

\begin{lemma}
  \label{lemma-kappa-star}
  Suppose that $\nu$ is a \grm{} given by \eqref{eq:30} as in the
  statement of Lemma~\ref{lem-prod-borel}.  Let
  $G^{(3)}=\set{(x,y,z)\in G\times G\times G: \text{$s(x)=r(y)$ and
      $s(y)=r(z)$}}$, let $\kappa:G^{(3)}\to B$ be given by
  $\kappa(x,y,z):=y^{-1}x$ and let $\kappa^{*}\B$ be the pull-back.
  If $F\in \sa_{c}(G^{(3)};\kappa^{*}\B)$, then
  \begin{equation*}
    (x,y,z)\mapsto \epsilon_{y^{-1}x}\bigl(F(x,y,z)\bigr)
  \end{equation*}
  is Borel.
\end{lemma}
\begin{proof}
  Sections of the form
  \begin{equation*}
    (x,y,z)\mapsto \phi(x,y,z)f(y^{-1}x)
  \end{equation*}
  with $\phi\in C_{c}(G^{(3)})$ and $f\in\gcb$ are dense in
  $\sa_{c}(G^{(3)};\kappa^{*}\B)$.  Since
  \begin{equation*}
    (x,y,z)\mapsto \phi(x,y,z)\epsilon_{y^{-1}x}\bigl(f(y^{-1}x)\bigr)
  \end{equation*}
  is clearly Borel (it is even continuous in the $z$ variable), the
  result follows as in the proof of Lemma~\ref{lem-prod-borel}.
\end{proof}

\begin{example}
  \label{ex-b-borel}
  We assume that we have the same set-up as in the previous two
  lemmas.  If $f,g,h\in\gcb$, then $F(x,y,z):=g(y)^{*}f(z)h(z^{-1}x)$
  defines a section in $\sa_{c}(G^{(3)};\kappa^{*}\B)$.  Then
  \begin{equation*}
    (x,y,z)\mapsto
    \epsilon_{y^{-1}x}\bigl(g(y^{*})f(z)h(z^{-1}x)\bigr) 
  \end{equation*}
  is Borel.
\end{example}

\section{Representations of Fell Bundles}
\label{sec:repr-fell-bundl}

In order to prove an equivalence theorem generalizing that for
groupoids (c.f., \citelist{\cite{ren:jot87}\cite{mrw:jot87}}), we will
need to work with the sort of ``weak representations'' introduced by
Renault in \cite{ren:jot87}.
\begin{definition}
  \label{def-pre-rep}
  Let $\H_{0}$ be a dense subspace of a Hilbert space $\H$, and denote
  by $\Lin(\H_{0})$\index{Lin@$\Lin(\H_{0})$} the collection of
  \emph{all} linear operators on the vector space $\H_{0}$.  A
  \emph{\prerep}{}\index{prerepresentation@\prerep}\index{Fell
    bundle!\prerep} of $\B$ on $\H_{0}\subset \H$ is a homomorphism
  $L:\gcb\to\Lin(\H_{0})$ such that for all $\xi,\eta\in\H_{0}$,
  \begin{enumerate}
  \item $f\mapsto \bip(L(f)\xi|\eta)$ is continuous in the \ilt{} on
    $\gcb$,
  \item $\bip(L(f)\xi|\eta)=\bip(\xi|{L(f^{*})\eta})$ and such that
  \item $\hoo:=\operatorname{span}\set{L(f)\zeta:\text{$f\in\gcb$ and
        $\zeta\in\H_{0}$}}$ is dense in $\H$.
  \end{enumerate}
  Of course, two \prerep s $(L,\H_{0},\H)$ and $(L',\H_{0}',\H')$ are
  equivalent if there is a unitary $U:\H\to\H'$ intertwining $L$ and
  $L'$ on $\H_{0}$ and $\H_{0}'$, respectively.
\end{definition}

The next result implies that each \prerep{} has associated to
it a unique measure class on $\go$.  This will be the measure class
that appears in the disintegration.\footnote{For the basics on Borel Hilbert
bundles and direct integrals, see \cite{wil:crossed}*{Appendix~F}.}

\begin{prop}
  \label{prop-measure-class}
  Suppose that $L:\gcb\to\Lin(\H_{0})$ is a \prerep{} of $\B$ on
  $\H_{0}\subset \H$.  Then there is a representation $M:C_{0}(\go)\to
  B(\H)$ such that for all $h\in C_{0}(\go)$, $f\in\gcb$ and
  $\xi\in\H_{0}$ we have
  \begin{equation}
    \label{eq:31}
    M(h)L(f)\xi=L\bigl((h\circ r)\cdot f\bigr)\xi.
  \end{equation}
  In particular, after replacing $L$ by an equivalent representation,
  we may assume that $\H=L^{2}(\go*\VV,\mu)$ for a Borel Hilbert
  bundle $\go*\VV$ and a finite Radon measure $\mu$ on $\go$ such that
  \begin{equation*}
    M(h)\xi(u)=h(u)\xi(u)\quad\text{for all $h\in C_{0}(\go)$ and
      $\xi\in L^{2}(\go*\VV,\mu)$.}
  \end{equation*}

\end{prop}
\begin{proof}
  We can easily make sense of $(h\circ r)\cdot f$ for $h\in
  C_{0}(\go)^{\sim}$.  Furthermore, we can compute that
  \begin{equation*}
    \bip(L((h\circ r)\cdot f)\xi|{L(g)\eta})=\bip(L(f)\xi|{L((\bar
      h\circ r)\cdot g)\eta}).
  \end{equation*}
  Then, if $k\in C_{0}(\go)^{\sim}$ is such that
  \begin{equation*}
    \|h\|_{\infty}^{2}1-|h|^{2}=|k|^{2},
  \end{equation*}
  we can compute that
  \begin{align*}
    \|h\|_{\infty}^{2}\Bigl\|\sum_{i=1}^{n} L(f_{i})\xi_{i}\Bigr\|^{2}
    -{}& \Bigl\|\sum _{i=1}^{n} L\bigl((h\circ r)\cdot
    f_{i}\bigr)\xi_{i}\Bigr\|^{2}  \\
    &=\sum_{ij} \bip(L\bigl(\bigl((\|h\|_{\infty}^{2}1-|h|^{2})\circ
    r\bigr)\cdot f_{i}\bigr)\xi_{i} | {L(f_{j}})\xi_{j}
    )\\
    &= \Bigl\|\sum_{i=1}^{n} L\bigl((k\circ r)\cdot f_{i}\bigr)
    \xi_{i}
    \Bigr\|^{2} \\
    &\ge 0.
  \end{align*}
  Since $\hoo$ is dense in $\H$, it follows that there is a
  well-defined bounded operator $M(h)$ on all of $\H$ satisfying
  \eqref{eq:31}.  It is not hard to see that $M$ is a
  $*$-homomorphism.  To see that $M$ is a representation, by
  convention, we must also see that $M$ is nondegenerate.  If $f\in
  \gcb$, then $r\bigl(\supp f\bigr)$ is compact.  Hence there is a
  $h\in C_{0}(\go)$ such that $M(h)f=f$.  From this, it is
  straightforward to see that $M$ is nondegenerate and therefore a
  representation.

  Since $M$ is a representation of $C_{0}(\go)$, it is equivalent to a
  multiplication representation on $L^{2}(\go*\VV,\mu)$ for an
  appropriate Borel Hilbert bundle $\go*\VV$ and finite Radon measure
  $\mu$ by, for example, \cite{wil:crossed}*{Example~F.25}.  The
  second assertion follows, and this completes the proof.
\end{proof}

\begin{definition}
  \label{def-endh}
  Suppose that $\go*\HH$ is a Borel Hilbert bundle over $\go$.  We
  let\index{Endgo@$\End(\go*\HH)$} 
  \begin{equation*}
    \End(\go*\HH)=\set{(u,T,v):T\in B\bigl(\H(v),\H(u)\bigr)}.
  \end{equation*}
   We give $\End(\go*\HH)$ the smallest Borel structure
  such that
  \begin{equation*}
    \psi_{f,g}(u,T,v):=\bip(Tf(u)|g(v))
  \end{equation*}
  is Borel for all Borel sections $f,g\in B(\go,\HH)$ (see
  \cite{wil:crossed}*{Definitions F.1~and F.6}).
\end{definition}

\begin{remark}
  \label{rem-end-std}
  The Borel structure on $\End(\go*\HH)$ is standard.  To see this,
  first note that $B_{k}(\H,\K):=\set{T\in B(\H,\K):\|T\|\le k}$ can
  be viewed as a closed subset of $B_{k}(\H\oplus\K)$ in the weak
  operator topology.  Since the latter is a compact Polish space (see
  \cite{wil:crossed}*{Lemma~D.37}), so is $B_{k}(\H,\K)$.  On the
  other hand, $\go*\HH$ is isomorphic to $\coprod_{n}
  X_{n}\times\H_{n}$ for Hilbert spaces $\H_{n}$ and a Borel partition
  $\set{X_{n}}$ of $X$ \cite{wil:crossed}*{Corollary~F.12}.  Then it
  is not hard to convince yourself that $\End(\go*\HH)$ is Borel
  isomorphic to $\bigcup_{k}\coprod_{m,n} X_{n}\times
  B_{k}(\H_{m},\H_{n})\times X_{m}$.
\end{remark}

\begin{definition}
  \label{def-star-functor}
  We say that a map $\hat\pi:\B\to \End(\go*\HH)$ is a
  \emph{$*$-functor}\index{star-functor@$*$-functor} if $\hat
  \pi(b)=\bigl(r(b),\pi(b),s(b)\bigr)$ for 
  an operator $\pi(b):\H(s(b))\to\H(r(b))$\footnote{Here we have
    adopted the notation that $r(b)$ in place of the cumbersome
    $r\bigl(p(b)\bigr)$ --- and similarly for $s(b)$.} such that
  \begin{enumerate}
  \item $\pi(\lambda a+b)=\lambda\pi(a)+\pi(b)$ if $p(a)=p(b)$,
  \item $\pi(ab)=\pi(a)\pi(b)$ if $(a,b)\in\B^{(2)}$ and
  \item $\pi(b^{*})=\pi(b)^{*}$.
  \end{enumerate}
  We say that a $*$-functor is \emph{Borel} if $x\mapsto
  \hat\pi\bigl(f(x)\bigr)$ is Borel from $G$ to $\End(\go*\HH)$ for
  all $f\in \gcb$.
\end{definition}

\begin{remark}
  \label{rem-bounded}
  Notice that if $\pi$ is
  a $*$-functor, then for each $u\in\go$, $\pi\restr{B(u)}$
  is a $*$-homomorphism and therefore bounded (since $B(u)$ is a
  \cs-algebra).  But then
  \begin{align*}
    \|\pi(b)\|^{2}=\|\pi(b)^{*}\pi(b)\| 
    =\|\pi(b^{*}b)\|
    \le\|b^{*}b\| 
    =\|b\|^{2}.
  \end{align*}
  Hence $*$-functors are ``naturally'' norm decreasing.
\end{remark}

\begin{definition}
  \label{def-repn}
  Suppose that $p:\B\to G$ is a Fell bundle.  Then a
  $*$-homomorphism\index{representation}\index{Fell
    bundle!representation} 
  $L$ of $\gcb$ into $B(\H)$ is a called a \emph{representation} if
  \begin{enumerate}
  \item it is continuous from $\gcb$ equipped with the \ilt{} into
    $B(\H)$ with the weak operator topology, and
  \item it is nondegenerate in that
    \begin{equation*}
      \operatorname{span}\set{L(f)\xi:\text{$f\in\gcb$ and $\xi\in\H$}}
    \end{equation*}
    is dense in $\H$.
  \end{enumerate}
\end{definition}

\begin{example}
  \label{ex-I-norm}
  If $L:\gcb\to B(\H)$ is a nondegenerate $*$-homomorphism which is
  bounded with respect to the $I$-norm --- that is, if $\|L(f)\|\le
  \|f\|_{I}$ for all $f\in\gcb$ --- then it is easy to see that $L$ a
  representation as defined in Definition~\ref{def-repn}.
\end{example}

Recall that a measure $\mu$ on $\go$ is called
\emph{quasi-invariant}\index{quasi-invariant measure}
if the Radon measure $\nu=\mu\circ\lambda$\index{nu@$\nu$} on $G$
defined by 
\begin{equation}
  \label{eq:33}
  \nu(f):=\int_{\go}\int_{G}f(x)\,d\lambda^{u}(x)\,d\mu(u)
\end{equation}
is equivalent to the measure $\nu^{-1}$ defined by
$\nu^{-1}(f)=\nu(\tilde f)$, where $\tilde f(x):=f(x^{-1})$.
(Alternatively, $\nu^{-1}$ is the push-forward of $\nu$ by the
inversion map on $G$: $\nu^{-1}(E)=\nu(E^{-1})$ for Borel sets
$E\subset G$.)  Note that
\begin{equation}
  \label{eq:34}
  \begin{split}
    \nu^{-1}(f) &=\int_{\go}\int_{G}
    f(x^{-1})\,d\lambda^{u}(x)\,d\mu(u)\\
    &= \int_{\go}\int_{G} f(x)\,d\lambda_{u}(x)\,d\mu(u).
  \end{split}
\end{equation}

If $\mu$ is quasi-invariant, so that $\nu$ is equivalent to
$\nu^{-1}$, then we can let
$\Delta=\frac{d\nu}{d\nu^{-1}}$\index{Delta@$\Delta$} be a 
Radon-Nikodym derivative for $\nu$ with respect to $\nu^{-1}$.  It
will be important for the calculations below to note that can choose
$\Delta:G\to\R^{+}$ to be a \emph{bona fide} homomorphism (see
\cite{muh:cbms}*{Theorem~3.15} or \cite{hah:tams78}*{Corollary~3.14}).
As is standard, we will write $\nu_{0}$\index{nu0@$\nu_{0}$} for the
symmetrized measure 
$\Delta^{-\half}\nu$.\footnote{We say the $\nu_{0}$ is ``symmetrized''
  because it is invariant under the inverse map:
  \begin{align*}
    \label{eq:35}
    \int_{G}f(x)\,d\nu_{0}(x)&=\int_{G}f(x)\Delta(x)^{-\half}\,d\nu(x)
    =\int_{G} f(x)\Delta(x)^{-\half}\Delta(x)\,d\nu^{-1}(x) \\
    &=\int_{G} f(x^{-1})\Delta(x^{-1})^{-\half}\Delta(x^{-1}) \,d\nu(x) 
    =\int_{G} f(x^{-1}) \Delta(x)^{-\half}\nu(x) \\
    &=\int_{G}f(x^{-1})\,d\nu_{0}(x).
  \end{align*}}

\begin{definition}
  \label{def-strict-rep}
  Suppose that $p:\B\to G$ is a Fell bundle over $G$.  Then
  a\index{integrated form}
  \emph{strict representation}\index{strict representation}\index{Fell
  bundle!strict representation} of $\B$ is a triple $(\mu,\go*\HH,\pi)$
  consisting of a quasi-invariant measure $\mu$ on $\go$, a Borel
  Hilbert bundle $\go*\HH$ and a $*$-functor $\pi:\B\to\End(\go*\HH)$.
\end{definition}

\begin{prop}
  \label{prop-strict-rep}
  Suppose that $p:\B\to G$ is a Fell bundle and that
  $(\mu,\go*\HH,\pi)$ is a strict representation of $\B$.  Then there
  is an associated $I$-norm bounded $*$-homomorphism $L$, called the
  \emph{integrated form of $\pi$}, on $L^{2}(\go*\HH,\mu)$ given by
  \begin{equation}
    \label{eq:36}
    \bip(L(f)\xi|\eta)= \int_{G} \bip(\pi\bigl(f(x)\bigr)
    \xi\bigl(s(x)\bigr) | {\eta\bigl(r(x)\bigr)})
    \Delta(x)^{-\half}\,d\nu(x) ,
  \end{equation}
  where $\Delta$ is the Radon-Nikodym derivative of $\nu^{-1}$ with
  respect to $\nu$ and $\nu=\mu\circ \lambda$.
\end{prop}
\begin{remark}
  \label{rem-alt-disc}
  Using vector-valued integrals, we can also write
  \begin{equation}
    \label{eq:37}
    L(f)\xi(u)=\int_{G}
    \pi\bigl(f(x)\bigr)\xi
    \bigl(s(x)\bigr)\Delta(x)^{-\half} \,d\lambda^{u}(x)
  \end{equation}
\end{remark}

\begin{remark}
  \label{rem-degeneracy}
  Note that there is no reason to suspect that $L$ is nondegenerate,
  and hence a representation, without some additional hypotheses on
  $\pi$.
\end{remark}

\begin{proof}
  The proof that $L$ is bounded is standard and uses the
  quasi-invariance of $\mu$.  The ``trick'' (due to Renault) is to
  apply the Cauchy-Schwartz inequality in $L^{2}(\nu)$ (and to recall
  that $\pi$ must be norm decreasing by Remark~\ref{rem-bounded}):
  \begin{align*}
    \bigl|\bip(L(f)\xi|\eta)\bigr| &\le \int_{G}
    \|f(x)\|\|\xi(s(x))\|\|\eta(r(x))\|\Delta(x)^{-\half}(x)\,d\nu(x) \\
    &\le \Bigl(
    \int_{G}\|f(x)\|\|\xi(s(x))\|^{2}\Delta(x)^{-1}\,d\nu(x)\Bigr)^{\half} \\
    &\hskip 1in
    \Bigl( \int_{G}\|f(x)\|\|\eta(r(x))\|^{2} \,d\nu(x)\Bigr)^{\half} \\
    &\le \bigl(\|f\|_{I}\|\xi\|_{2}^{2}\bigr)^{\half}\bigl(\|f\|_{I}
    \|\eta\|_{2}^{2} \bigr)^{\half} \\
    &\le \|f\|_{I}\|\xi\|_{2}\|\eta\|_{2}.
  \end{align*}
  To show that $L$ is multiplicative, we use \eqref{eq:37} and compute
  as follows:
  \begin{align*}
    L(f*g)\xi(u)&=\int_{G}\pi\bigl(f*g(x)\bigr)\xi\bigl(s(x)\bigr)
    \Delta(x)^{-\half} \, d\lambda^{u}(x) \\
    \intertext{which, since $\pi\restr{B(x)}$ is a bounded linear map
      of $B(x)$ into $B(\H(s(x)),\H(r(x))$, is} &= \int_{G}\int_{G}
    \pi\bigl(f(y)g(y^{-1}x)\bigr) \xi\bigl(s(x)\bigr)
    \Delta(x)^{-\half} \, d\lambda^{u}(y)\,d\lambda^{u}(x) \\
    \intertext{which, after using Fubini and sending $x\mapsto yx$,
      is} &= \int_{G}\int_{G} \pi\bigl(f(y)g(x)\bigr)
    \xi\bigl(s(x)\bigr)
    \Delta(yx)^{-\half} \,d\lambda^{s(y)}(x) \,d\lambda^{u} (y) \\
    &=\int_{G}\pi\bigl(f(y)\bigr) \int_{G} \pi\bigl(g(x)\bigr)
    \xi\bigl(s(x)\bigr) \Delta(x)^{-\half} \,d\lambda^{s(y)}(x)
    \Delta(y)^{-\half}
    \,d\lambda^{u}(y) \\
    &= \int_{G}\pi\bigl(f(y)\bigr) L(g)\xi\bigl(s(y)\bigr)
    \Delta(y)^{-\half} \,
    d\lambda^{u}(y) \\
    &=L(f)L(g)\xi(u).
  \end{align*}
  To see that $L$ is $*$-preserving, we will need to use the
  quasi-invariance of $\mu$ in the form of the invariance of
  $\Delta^{-\half}\,d\nu$ under the inversion map.  We compute that
  \begin{align*}
    \bip(L(f^{*}\xi|\eta) &= \int_{G} \bip(\pi\bigl(f^{*}(x)\bigr)
    \xi\bigl(s(x)\bigr) |{\eta\bigl(r(x)\bigr)}) \Delta(x)^{-\half}
    \,d\nu(x) \\
    &= \int_{G} \bip(\pi\bigl(f(x^{-1})^{*}\bigr)\xi\bigl(s(x)\bigr)|
    {\eta\bigl(r(x)\bigr)}) \Delta(x)^{-\half} \,d\nu(x) \\
    \intertext{which, after sending $x\to x^{-1}$, is} &=\int_{G} \ip
    ( \xi\bigl(r(x)\bigr) | {\pi\bigl(f(x)\bigr)
      \eta\bigl(s(x)\bigr)}) \Delta(x)^{-\half} \,d\nu(x).\qedhere
  \end{align*}
\end{proof}

The next step is to prove a very strong
converse of Proposition~\ref{prop-strict-rep} modeled after Renault's
\cite{ren:jot87}*{Proposition~4.2}.  The extra generality is needed to
prove the equivalence theorem --- which is our eventual goal.  The
proof given in the next section follows Yamagami's suggestion that
Proposition~\ref{prop-din28.32} ought to ``replace'' the
Radon-Nikodym theorem in Renault's proof in the presence of suitable
approximate identities (see Proposition~\ref{prop-approx-id}).
The argument here follows Muhly's version of Renault's argument 
(see \cite{muh:cbms}*{Theorem~3.32} or
\cite{muhwil:nyjm08}*{Theorem~7.8})
with a 
couple of ``vector upgrades''.

\begin{thm}[Disintegration Theorem]
  \label{thm-fell-disintegration}\index{disintegration theorem}
  Suppose that $L:\gcb\to \Lin(\H_{0})$ is a \prerep{} of $\B$ on
  $\H_{0}\subset \H$.  Then $L$ is bounded in the sense that
  $\|L(f)\|\le\|f\|_{I}$ for all $f\in\gcb$.  Therefore $L$ extends to
  a \emph{bona fide} representation of $\gcb$ on $\H$ which is
  equivalent to the integrated form of a strict representation
  ($\mu,\go*\HH,\pi)$ of $\B$ where $\mu$ is the measure defined in
  Proposition~\ref{prop-measure-class}.\footnote{Notice that we
    are asserting that $\mu$ is necessarily quasi-invariant.} In
  particular, $L$ is bounded with respect to the universal \cs-norm on
  $\gcb$. 
\end{thm}

We will take up the proof of the Disintegration Theorem in the next
section.  However, once we have Theorem~\ref{thm-fell-disintegration}
in hand, we can ``adjust'' our definition of the universal norm as
follows.\index{universal norm}\index{gammacb@$\gcb$!universal norm}

\begin{remark}
  \label{rem-norm-bdd}
  Since any \prerep,  and \emph{a priori} any representation, $L$, of
  $\B$ is equivalent to the integrated form of a strict
  representation, it follows that $L$ is $I$-norm bounded by
  Proposition~\ref{prop-strict-rep}.  Conversely, $I$-norm
  bounded representations are clearly representations.  Therefore, we
  could have defined the universal norm on $\gcb$ via
\begin{equation*}
  \|f\|:=\sup\set{\|L(f)\|:\text{$L$ is a representation of $\B$.}}
\end{equation*}
\end{remark}

\section{Proof of the Disintegration Theorem}
\label{sec:proof-disint-theor}

Naturally, we will break the proof up into a number of steps.  The
first is, as suggested by Yamagami, to produce a two sided approximate
identity for $\gcb$ in the inductive limit topology.  This is a highly
nontrivial result.  However, it is an immediate
consequence of the rather special approximate identities that we need
for our proof of the equivalence theorem.  Thus the next result is an
immediate consequence of 
Proposition~\ref{prop-ai} (see the comments immediately following that
proposition).

\begin{prop}
  \label{prop-approx-id}
  There is a self-adjoint approximate identity for $\gcb$ in the
  inductive limit topology.\index{gammacb@$\gcb$!approximate identity}
\end{prop}

\begin{cor}
  \label{cor-dense}
  Suppose that $L$ is a \prerep{} of $\B$ on $\H_{0}\subset\H$.  Let
  $\hoo$ be the necessarily dense subspace
  $\operatorname{span}\set{L(f)\xi:\text{$f\in\gcb$ and
      $\xi\in\H_{0}$.}}$ If $\hoo'$ is a dense subspace of $\hoo$,
  then
  \begin{equation*}
    \operatorname{span}\set{L(f)\xi:\text{$f\in \gcb$ and $\xi\in\hoo'$}}
  \end{equation*}
  is dense in $\H$.
\end{cor}

\begin{proof}
  Let $\set{e_{i}}$ be a self-adjoint approximate identity for $\gcb$
  in the \ilt.  Then if $L(f)\xi\in\hoo$, we see that
  \begin{multline*}
    \|L(e_{i})L(f)\xi- L(f)\xi\|^{2} =\\
    \bip(L(f^{*}*e_{i}*e_{i}*f)\xi|\xi)
    -2\operatorname{Re}\bip(L(f^{*}*e_{i}*f)\xi|\xi) +
    \bip(L(f^{*}*f)\xi|\xi),
  \end{multline*}
  which tends to zero by part~\partref1 of
  Definition~\ref{def-pre-rep}.  It follows that $\hoo' \subset
  \overline{\operatorname{span}}\set{L(f)\xi:\text{$\xi\in\hoo'$ and
      $f\in\gcb$}}$.  Since $\hoo'$ is dense, the result follows.
\end{proof}

\begin{lemma}
  \label{lem-alg-tensor}
  Suppose that $L$ is a \prerep{} of $\gcb$ on $\H_{0}\subset\H$.
  Then there is a positive sesquilinear form $\rip<\cdot,\cdot>$ on
  $\gcbatho$ such that
  \begin{equation}
    \label{eq:38}
    \rip<f\tensor \xi,g\tensor\eta>=\bip(L(g^{*}*f)\xi|\eta).
  \end{equation}
  Furthermore, the Hilbert space completion $\mathcal{K}$ of
  $\gcbatho$ is isomorphic to $\H$.  In fact, if $[f\tensor\xi]$ is
  the class of $f\tensor \xi$ in $\mathcal{K}$, then
  $[f\tensor\xi]\mapsto L(f)\xi$ is well-defined and induces an
  isomorphism of $\mathcal{K}$ with $\H$ which maps the quotient
  $\gcbatho/\N$, where $\N$ is the subspace
  $\N=\set{\sum_{i}f_{i}\tensor \xi_{i}:\sum_{i}L(f_{i})\xi_{i}=0}$ of
  vectors in $\gcbatho$ of length zero, onto $\hoo$ (as defined in
  part~\partref3 of Definition~\ref{def-pre-rep}).
\end{lemma}
\begin{proof}
  Using the universal properties of the \emph{algebraic} tensor
  product, as in the proof of \cite{rw:morita}*{Proposition~2.64} for
  example, it is not hard to see that there is a unique sesquilinear
  form on $\gcbatho$ satisfying \eqref{eq:38}.\footnote{For fixed $g$
    and $\eta$, the left-hand side of \eqref{eq:38} is bilinear in $f$
    and $\xi$.  Therefore, by the universal properties of the
    algebraic tensor product, \eqref{eq:38} defines linear map
    $m(g,\eta):\gcbatho\to\C$.  Then $(g,\eta)\mapsto
    \overline{m(g,\eta)}$ is a bilinear map into the space
    $\operatorname{CL}(\gcbatho)$ of conjugate linear functionals on
    $\gcbatho$.  Then we get a linear map $N:\gcbatho\to
    \operatorname{CL}(\gcbatho)$.  We can then define
    $\rip<\alpha,\beta>:=\overline{N(\beta)(\alpha)}$.  Clearly
    $\alpha\mapsto \rip<\alpha,\beta>$ is linear and it is not hard to
    check that $\overline{\rip<\alpha,\beta>}=\rip<\beta,\alpha>$.}
  Thus to see that $\rip<\cdot,\cdot>$ is a pre-inner product, we just
  have to see that it is positive.  But
  \begin{equation}\label{eq:39}
    \begin{split}
      \brip<\sum_{i}f_{i}\tensor\xi_{i},\sum_{i} f_{i}\tensor\xi_{i}>
      &= \sum_{ij}\bip(L(f_{j}^{*}*f_{i})\xi_{i}|\xi_{j}) \\
      &= \sum_{ij} \bip(L(f_{i})\xi_{i}|{L(f_{j})\xi_{i}}) \\
      &=\bigl\|\sum_{i}L(f_{i})\xi_{i}\bigr\|^{2}.
    \end{split}
  \end{equation}
  As in \cite{rw:morita}*{Lemma~2.16}, $\rip<\cdot,\cdot>$ defines an
  inner-product on $\gcbatho/\N$, and $[f_{i}\tensor\xi]\mapsto
  L(f_{i})\xi$ is well-defined in view of \eqref{eq:39}.  Since this
  map has range $\hoo$ and since $\hoo$ is dense in $\H$ by
  definition, the map extends to an isomorphism of $\mathcal K$ onto
  $\H$ as claimed.
\end{proof}

\begin{conventions}
  Using Lemma~\ref{lem-alg-tensor}, we will\index{conventions}
  identify $\H$ with $\K$, and $\hoo$ with $\gcbatho/\N$.  Thus we
  will interpret $[f\tensor\xi]$ as a vector in
  $\hoo\subset\H_{0}\subset \H$.  Then we have
  \begin{align}
    \label{eq:40}
    L(g)[f\tensor\xi]&=[g*f\tensor \xi]\quad\text{and}\\
    M(h)[f\tensor\xi]&=[(h\circ r)\cdot f\tensor\xi],\label{eq:41}
  \end{align}
  where $M$ is the representation of $C_{0}(\go)$ defined in
  Proposition~\ref{prop-measure-class}, $g\in\gcb$ and $h\in
  C_{0}(\go)$.
\end{conventions}

\begin{remark}
  \label{rem-extend-m}
  In view of Proposition~\ref{prop-measure-class}, $M$ extends to a
  $*$-homomorphism of $\bb_{c}(G)$ into $B(\H)$ such that $M(h)=0$ if
  $h(u)=0$ for $\mu$-almost all $u$ (where $\mu$ is the measure
  defined in that proposition).  However, at this point, we can not
  assert that \eqref{eq:41} holds for any $h\notin C_{0}(\go)$.
\end{remark}

A critical step in producing a strict representation is producing a
quasi-invariant measure class.  While we have the measure $\mu$
courtesy of Proposition~\ref{prop-measure-class}, showing that $\mu$
is quasi-invariant requires that we extend equations \eqref{eq:40} and
\eqref{eq:41} to a larger class of functions.  This can't be done
without also enlarging the domain of definition of $L$.  This is
problematic as we don't as yet know that each $L(f)$ is bounded in any
sense, nor have we assumed that
$\H_{0}$ is complete.  Motivated by Muhly's proof in \cite{muh:cbms},
we have introduced $\sgcb$ and $\bboc(G)$ in
Definition~\ref{def-sigma1} in order to deal with only those additional
functions that we absolutely need.

``Not to put to fine a point on it,'' $\sgcb$ is not a well-behaved
class of functions on $G$. For example, there is no reason to suspect
that it is closed under the sort of uniformly bounded pointwise
convergence used in its definition.  Nevertheless, we have the
following useful observation.

\begin{lemma}
  \label{lem-grm-ext}
  Suppose that $\sigma$ is a \grm{} on $\B$ given by
  \begin{equation}
    \label{eq:42}
    \sigma(f)=\int_{G} \epsilon_{x}\bigl(f(x)\bigr) \,d|\sigma|(x)
  \end{equation}
  as in Proposition~\ref{prop-din28.32}.  If $f\in\sgcb$, then
  $x\mapsto \epsilon_{x}\bigl(f(x)\bigr)$ is in $\bb_{c}(G)$ and we
  can extend $\sigma$ to a linear functional on $\sgcb$.  In
  particular, if $\set{f_{n}}\subset\sgcb$ is a uniformly bounded
  sequence whose supports are contained in a fixed compact set and
  which converges pointwise to $f\in\sgcb$, then $\sigma(f_{n})\to \sigma(f)$.
\end{lemma}
\begin{proof}
  Let $\set{f_{n}}\subset\gcb$ be as in the second part of the lemma.
  Let $\tau_{n}(x):=\epsilon_{x}\bigl(f_{n}(x)\bigr)$.
  Proposition~\ref{prop-din28.32} implies that
  $\tau_{n}\in\bb_{c}(G)$ and clearly $\tau_{n}(x)\to
  \tau(x):=\epsilon_{x}\bigl(f(x)\bigr)$ for each $x\in G$.  Moreover
  the sequence $\set{\tau_{n}}$ is uniformly bounded and vanishes off
  some compact set.  Therefore $\tau\in\bb_{c}(G)$ and we can extend
  $\sigma$ using \eqref{eq:42}.  Furthermore, using the dominated
  convergence theorem we have
  \begin{align*}
    \sigma(f)&:=\int_{G}\epsilon_{x}\bigl(f(x)\bigr)\,d|\sigma|(x) \\
    &= \lim_{n}\int_{G}\epsilon_{x}\bigl(f_{n}(x)\bigr) \,d|\sigma|(x) \\
    &= \lim_{n}\sigma(f_{n}).
  \end{align*}
  Using this, it is not hard to see that the extension of $\sigma$ is
  a linear functional.

  If $\set{f_{n}}\subset \sgcb$ converges to $f\in\sgcb$ as in the
  second part of the lemma, then we can define $\tau_{n}$ and $\tau$
  as above.  Thus just as above, the dominated convergence theorem
  implies the final assertion in the lemma.
\end{proof}

\begin{lemma}
  \label{lem-sgcb-conv}
  Suppose that $f,g\in\sgcb$.  Then
  \begin{equation*}
    f*g(x):=\int_{G} f(y)g(y^{-1}x)\,d\lambda^{r(x)}(y)
  \end{equation*}
  is a well-defined element of $B(x)$, and $f*g$ defines a section in
  $\sgcb$.
\end{lemma}
\begin{proof}
  Let $\set{f_{n}}$ and $\set{g_{n}}$ be uniformly bounded sequences
  in $\gcb$ with supports all in the same compact set such that
  $f_{n}\to f$ and $g_{n}\to g$ pointwise.  Then for each $x$,
  $y\mapsto f_{n}(y)g_{n}(y^{-1}x)$ converges pointwise to $y\mapsto
  f(y) g(y^{-1}x)$.  Therefore the latter is a bounded Borel function
  from $G^{r(x)}$ to $B(x)$ vanishing off a compact set .  Thus
  $f*g(x)$ is a well-defined element of $B(x)$ (for example, by
  \cite{wil:crossed}*{Lemma~1.91}).

  Furthermore, since
  \begin{equation*}
    \|f_{n}*g_{n}(x)\|\le\|f_{n}\|_{\infty}\|g_{n}\|_{\infty} \sup_{u}
    \lambda^{u} \bigl((\supp f_{n})(\supp g_{n})\bigr),
  \end{equation*}
  $\set{f_{n}*g_{n}}$ is uniformly bounded sequence in $\gcb$ (by
  Corollary~\ref{cor-conv}) whose supports are all contained in a
  fixed compact set and which converges pointwise to $f*g$.  Hence
  $f*g\in\sgcb$ as claimed.
\end{proof}

For each $\xi$ and $\eta$ in $\H_{0}$, part~\partref2 of
Definition~\ref{def-pre-rep} implies that
\begin{equation}
  \label{eq:43}
  L_{\xi,\eta}(f):=\bip(L(f)\xi|\eta)
\end{equation}
defines a \grm{} on $\B$.  We will use Lemma~\ref{lem-grm-ext} to
   extend $L_{\xi,\eta}$ to $\sgcb$.

\begin{lemma}
  \label{lem-paul2}
  Suppose that $L$ is a \prerep{} of $\B$ on $\H_{0}\subset\H$.  Then
  there is a positive sesquilinear form on $\sgcbatho$, extending that
  on $\gcbatho$, such that
  \begin{equation*}
    \rip<f\tensor \xi,g\tensor
    \eta>=L_{\xi,\eta}(g^{*}*f)\quad\text{for all $f,g\in\sgcbatho$ and
      $\xi,\eta\in \H_{0}$.}
  \end{equation*}
  In particular, if
  \begin{equation*}
    \N_{b}:=\set{\sum_{i}f_{i}\tensor\xi\in\sgcbatho
      :\brip<\sum_{i}f_{i}\tensor \xi, \sum_{i}f_{i}\tensor\xi_{i}>=0}
  \end{equation*}
  is the subspace of vectors of zero length, then the quotient
  $\sgcbatho/\N_{b}$ can be identified with a subspace of $\H$
  containing $\hoo:=\gcbatho/\N$.
\end{lemma}
\begin{remark}
  \label{rem-notation}
  As before, we will write $[f\tensor\xi]$ for the class of
  $f\tensor\xi$ in the quotient $\sgcbatho/\N_{b}\subset \H$.
\end{remark}

\begin{proof}
  Just as in Lemma~\ref{lem-alg-tensor}, there is a well-defined
  sesquilinear form on $\sgcbatho$ extending that on $\gcbatho$.
  (Note that the right-hand side of \eqref{eq:38} can be rewritten as
  $L_{\xi,\eta}(g^{*}*f)$.) In particular, we have
  \begin{equation*}
    \brip<\sum_{i} f_{i}\tensor\xi_{i},\sum_{j} g_{j}\tensor\eta_{j}> =
    \sum_{ij} L_{\xi_{i},\eta_{j}}(g_{j}^{*}*f_{i}).
  \end{equation*}
  We need to see that the form is positive.  Let $\alpha:=\sum_{i}
  f_{i}\tensor \xi_{i}$, and let $\set{f_{i,n}}$ be a uniformly
  bounded sequence in $\gcb$ converging pointwise to $f_{i}$ with all
  the supports contained in a fixed compact set.  Then for each $i$
  and $j$, $f_{j,n}^{*}*f_{i,n}\to f_{j}^{*}*f_{i}$ in the appropriate
  sense.  In particular, Lemma~\ref{lem-grm-ext} implies that
  \begin{align*}
    \rip<\alpha,\alpha>&= \sum_{ij}
    L_{\xi_{i},\xi_{j}}(f_{j}^{*}*f_{i}) \\
    &= \lim_{n} \sum_{ij} L_{\xi_{i},\xi_{j}}(f_{j,n}^{*}*f_{i,n}) \\
    &= \lim_{n} \rip<\alpha_{n},\alpha_{n}>,
  \end{align*}
  where $\alpha_{n}:=\sum_{i}f_{i,n}\tensor \xi_{i}$.
  Since$\rip<\cdot,\cdot>$ is positive on $\gcbatho$ by
  Lemma~\ref{lem-alg-tensor}, we have $\rip<\alpha_{n},
  \alpha_{n}>\ge0$, and we've shown that $\rip<\cdot,\cdot>$ is still
  positive on $\sgcbatho$.

  Clearly the map sending the class $f\tensor\xi+\N$ to
  $f\tensor\xi+\N_{b}$ is isometric and therefore extends to an
  isometric embedding of $\H$ into the Hilbert space completion
  $\H_{b}$ of $\sgcbatho$ with respect to $\rip<\cdot,\cdot>$.
  However if $g\tensor\xi \in\sgcbatho$ and if $\set{g_{n}}$ is a
  sequence in $\gcb$ such that $g_{n}\to g$ in the usual way, then
  \begin{equation*}
    \|(g_{n}\tensor\xi+\N_{b})-(g\tensor\xi+\N_{b})\|^{2} =
    L_{\xi,\xi}(g_{n}^{*}*g_{n} -g_{n}^{*}*g -g*g_{n}^{*}+g^{*}*g),
  \end{equation*}
  and this tends to zero by Lemma~\ref{lem-grm-ext}.  Thus the image
  of $\H$ in $\H_{b}$ is all of $\H_{b}$.  Consequently, we can
  identify the completion of $\sgcbatho$ with $\H$ and
  $\sgcbatho/\N_{b}$ with a subspace of $\H$ containing $\hoo$.
\end{proof}

The ``extra'' vectors provided by $\sgcbatho/\N_{b}$ are just enough
to allow us to use a bit of general nonsense about unbounded operators
to extend the domain of each $L(f)$.  More precisely, for $f\in \gcb$,
we can view $L(f)$ as an operator in $\H$ with domain $\D(L(f))=\hoo$.
Then using part~\partref2 of Definition~\ref{def-pre-rep}, we see
that
\begin{equation*}
  L(f^{*})\subset L(f)^{*}.
\end{equation*}
This implies that $L(f)^{*}$ is a densely defined operator.  Hence
$L(f)$ is closable \cite{con:course}*{Proposition~X.1.6}.
Consequently, the closure of the graph of $L(f)$ in $\H\times\H$ is
the graph of the closure $\overline{L(f)}$ of $L(f)$
\cite{con:course}*{Proposition~X.1.4}.

Suppose that $g\in\sgcb$.  Let $\set{g_{n}}$ be a uniformly bounded
sequence in $\gcb$ all supported in a fixed compact set such that
$g_{n}\to g$ pointwise.  Then
\begin{equation}\label{eq:44}
  \|[g_{n}\tensor\xi]-[g\tensor\xi]\|^{2} =
  L_{\xi,\xi}(g_{n}^{*}*g_{n}-g^{*}*g_{n} -g_{n}^{*}*g +g*g). 
\end{equation}
However $\set{g_{n}^{*}*g_{n}-g^{*}*g_{n} -g_{n}^{*}*g +g*g}$ is
uniformly bounded and converges pointwise to zero.  Since the supports
are all contained in a fixed compact set, the left-hand side of
\eqref{eq:44} tends to zero by Lemma~\ref{lem-grm-ext}.  Similarly,
\begin{equation*}
  \|[f*g_{n}\tensor\xi]-[f*g\tensor\xi]\|^{2}\to 0.
\end{equation*}
If follows that
\begin{equation*}
  \bigl([g_{n}\tensor\xi,L(f)[g_{n}\tensor \xi]\bigr) \to
  \bigl([g\tensor\xi],[f*g\tensor\xi]\bigr) 
\end{equation*}
in $(\sgcbatho/\N_{b})\times (\sgcbatho/\N_{b})\subset\H\times\H$.
Therefore $[g\tensor\xi]\in\D\bigl(\overline{L(f)}\bigr)$ and
$\overline{L(f)}[g\tensor\xi]= [f*g\tensor\xi]$.  We have proved the
following.

\begin{lemma}
  \label{lem-ext-lb}
  For each $f\in \gcb$, $L(f)$ is a closable operator in $\H$ with
  domain $\D(L(f))=\hoo=\gcbatho/\N$.  Furthermore $\sgcbatho/\N_{b}$
  belongs to $\D\bigl(\overline{L(f)}\bigr)$, and
  \begin{equation*}
    \overline{L(f)}[g\tensor\xi]=[f*g\tensor\xi]\quad\text{for all
      $f\in \gcb$, $g\in\sgcb$ and $\xi\in\H_{0}$.}
  \end{equation*}

\end{lemma}

\begin{lemma}
  \label{lem-lb}
  For each $f\in\sgcb$, there is a well-defined operator $L_{b}(f)\in
  \Lin(\sgcbatho)/N_{b})$ such that
  \begin{equation}
    \label{eq:45}
    L_{b}(f)[g\tensor\xi]=[f*g\tensor\xi].
  \end{equation}
  If $f\in\gcb$, then $L_{b}(f)\subset\overline{L(f)}$.
\end{lemma}
\begin{proof}
  To see that \eqref{eq:45} determines a well-defined operator, we
  need to see that
  \begin{equation}
    \label{eq:46}
    \sum_{i}[g_{i}\tensor\xi_{i}]=0\quad\text{implies}\quad \sum_{i}
    [f*g_{i}\tensor \xi_{i}]=0.
  \end{equation}
  However,
  \begin{equation}\label{eq:47}
    \bigl\|\sum_{i} [f*g_{i}\tensor\xi_{i}]\bigr\|^{2}=\sum_{ij}
    L_{\xi_{i},\xi_{j}} (g_{j}^{*}*f^{*}*f*g_{i}).
  \end{equation}
  Since $f\in\sgcb$, we can approximate the right-hand side of
  \eqref{eq:47} by sums of the form
  \begin{equation}
    \label{eq:48}
    \sum_{ij} L_{\xi_{i},\xi_{j}} (g_{j}^{*}*h^{*}*h*g_{i}),
  \end{equation}
  where $h\in\gcb$.  But \eqref{eq:48} equals
  \begin{equation*}
    \bigl\| \overline{L(h)} \sum_{i}[g_{i}\tensor\xi_{i}]\bigr\|^{2}
  \end{equation*}
  which is zero if the left-hand side of \eqref{eq:46} is zero.  Hence
  the right-hand side of \eqref{eq:46} is also zero and $L_{b}(f)$ is
  well-defined.

  If $f\in\gcb$, then $L_{b}(f)\subset \overline{L(f)}$ by
  Lemma~\ref{lem-ext-lb}.
\end{proof}

Now we prove the analogue of
Muhly's technical lemma (\cite{muh:cbms}*{Lemma~3.33} or
\cite{muhwil:nyjm08}*{Lemma~B.12}) 
which will allow us to compute with Borel functions as we would
expect.

\begin{lemma}
  \label{lem-paul-tech}
  Suppose that $f\in\sgcb$ and that $k$ is a bounded Borel function on
  $\go$ which is the pointwise limit of a uniformly bounded sequence
  from $C_{0}(\go)$.  Then for all $g,h\in\gcb$ and
  $\xi,\eta\in\H_{0}$, we have the following.
  % \begin{enumerate}
  % \item
  %   \begin{align}
  %     \bip(L_{b}(f)[g\tensor\xi]|[h\tensor\eta]) &=
  %     \bip([f*g\tensor\xi]|[ h\tensor\eta]) \\
  %     &= L_{\xi,\eta}(g^{*}*f*h) \\
  %     &= L_{[g\tensor\xi],[h\tensor\eta]}(f)
  %   \end{align}
  % \item
  %   \begin{align*}
  %     \big(M(k)[g\tensor\xi]|[h\tensor\eta]) &=
  %     L_{\xi,\eta}(h^{*}*((k\circ r)\cdot g)) \\
  %     &= \bip([(k\circ r)\cdot g\tensor \xi]|[h\tensor\eta]) \\
  %     &= \bip(M(k)L(g)\xi|L(h)\eta).
  %   \end{align*}
  % \item
  %   \begin{equation*}
  %     \bip(M(k)L_{b}(f)[g\tensor\xi]|[h\tensor\eta]) =
  %     \bip(L_{b}((k\circ r)\cdot f)[g\tensor\xi]|[h\tensor\eta]).
  %   \end{equation*}
  % \end{enumerate}
  \begin{align}
    \bip(L_{b}(f)[g\tensor\xi]|[h\tensor\eta]) &=
    \bip([f*g\tensor\xi]|[ h\tensor\eta]) \tag{a}\\
    &= L_{\xi,\eta}(g^{*}*f*h) \notag\\
    &= L_{[g\tensor\xi],[h\tensor\eta]}(f)\notag \\[1ex]
    \big(M(k)[g\tensor\xi]|[h\tensor\eta]) &=
    L_{\xi,\eta}(h^{*}*((k\circ r)\cdot g)) \tag{b}\\
    &= \bip([(k\circ r)\cdot g\tensor \xi]|[h\tensor\eta])\notag \\
    &= \bip(M(k)L(g)\xi|L(h)\eta)\notag\\[1ex]
    \bip(M(k)L_{b}(f)[g\tensor\xi]|[h\tensor\eta]) &=
    \bip(L_{b}((k\circ r)\cdot f)[g\tensor\xi]|[h\tensor\eta]).\tag{c}
  \end{align}
\end{lemma}
\begin{proof}
  We start with \partref1.  The first equality is just the definition
  of $L_{b}(f)$.  The second follows from the definition of the inner
  product on $\sgcbatho/\N_{b}$.  If $f$ is in $\gcb$, then the third
  equation holds just by untangling the definition of the \grm s
  $L_{\xi,\eta}$ and using property \partref2 of
  Definition~\ref{def-pre-rep}. Therefore the third equality holds
  for $f\in\sgcb$ by applying the continuity assertion in
  Lemma~\ref{lem-grm-ext}.

  Part~\partref2 is proved similarly.  The first equation holds if
  $k\in C_{0}(\go)$ by definition of $M(k)$ and $L_{\xi,\eta}$.  If
  $\set{k_{n}}\subset C_{0}(\go)$ is a bounded sequence converging
  pointwise to $k$, then $M(k_{k})\to M(k)$ in the weak operator
  topology by the dominated convergence theorem.  On the other hand
  $g^{*}*(k_{n}\circ r)\cdot h\to g^{*}*(k\circ r)h$ in the required
  way.  Thus $L_{\xi,\eta}(g^{*}*(k_{n}\circ r)\cdot h) \to
  L_{\xi,\eta}(g^{*}*(k\circ r)h)$ by Lemma~\ref{lem-grm-ext}.  Thus
  the first equality is valid.  The second equality if clear if $k\in
  C_{0}(\go)$ and passes to the limit as above.  The third equality is
  simply our identification of $[g\tensor\xi]$ with $L(g)\xi$ as in
  Lemma~\ref{lem-alg-tensor}.

  For part~\partref3, first note that if $f_{n}\to f$ and $k_{n}\to k$
  are uniformly bounded sequences converging pointwise with supports
  in fixed compact sets independent of $n$, then $(k\circ r)\cdot
  f=\lim_{n} (k_{n}\circ r)\cdot f_{n}$.  It follows that $(k\circ
  r)\cdot f\in\sgcb$.  Also, $[f\tensor\xi]=\lim[f_{n}\tensor \xi]$,
  and since $M(k)$ is bounded, part~\partref2 implies that
  \begin{align*}
    M(k)[f\tensor\xi]&=\lim_{n}M(k)[f_{n}\tensor\xi] \\
    &= \lim_{n}[(k\circ r)\cdot f_{n}\tensor\xi]\\
    &= [(k\circ r)\cdot f\tensor \xi].
  \end{align*}
  Since it is not hard to verify that $M(k)^{*}[f\tensor\xi]= (\bar
  k\circ r)\cdot f\tensor \xi]$, we can compute that
  \begin{align*}
    \bip(M(k)L_{b}(f)[g\tensor\xi]|[h\tensor\eta]) &=
    \bip([f*g\tensor\xi] | (\bar k\circ r)\cdot h\tensor\eta]) \\
    &= \bip([k\circ r)\cdot(f*g)\tensor\xi]|[h\tensor\eta]) \\
    &=\bip([((k\circ r)\cdot f)*g\tensor \xi]|[h\tensor\eta]) \\
    &= \bip(L_{b}((k\circ r)\cdot
    f)[g\tensor\xi]|[h\tensor\eta]).\qedhere
  \end{align*}
\end{proof}

\begin{prop}
  \label{prop-quasi-inv}
  Let $\mu$ be the Radon measure on $\go$ associated to the \prerep{}
  $L$ by Proposition~\ref{prop-measure-class}.  Then $\mu$ is
  quasi-invariant.
\end{prop}
\begin{proof}
  We need to show that measures $\nu$ and $\nu^{-1}$ (defined in
  \eqref{eq:33} and \eqref{eq:34}, respectively) are equivalent.
  Therefore, we have to show that if $A$ is pre-compact in $G$, then
  $\nu(A)=0$ if and only if $\nu(A^{-1})=0$.  Since $(A^{-1})^{-1}=A$,
  it's enough to show that $A$ $\nu$-null implies that $A^{-1}$ is
  too.  Since $\nu$ is regular, we may as well assume that $A$ is a
  $G_{\delta}$-set so that $\phi:=\charfcn A$ is in $\bboc(G)$.  Let
  $\tilde\phi(x)=\phi(x^{-1})$.  We need to show that $\tilde
  \phi(x)=0$ for $\nu$-almost all $x$.  Since $A$ is a $G_{\delta}$,
  we can find a sequence $\set{\phi_{n}}\subset \ccpg$ such that
  $\phi_{n}\searrow \phi$.

  If $\psi$ is any function in $\ccg$, then
  $\psi\phi\bar\psi=|\psi|^{2}\phi\in\bboc(G)$ and vanishes
  $\nu$-almost everywhere.  By the monotone convergence theorem,
  \begin{equation*}
    \lambda(\psi\phi\bar\psi)(u) 
    :=\int_{G}|\psi(x)|^{2}\phi(x)\,d\lambda^{u}(x) 
  \end{equation*}
  defines a function in $\bboc(\go)$ which is equal to $0$ for
  $\mu$-almost all $u$.  In particular,
  $M(\lambda(\psi\phi\bar\psi))=0$.

  Therefore
  \begin{equation}
    \label{eq:49}
    0=\bip(M\bigl(\lambda(\psi\phi\bar\psi)\bigr)L(g)\xi|L(g)\xi) =
    L_{\xi,\xi} (g^{*}*\bigl(\lambda(\psi\phi\bar\psi)\circ r)\cdot
    g\bigr) 
  \end{equation}
  for all $g\in\gcb$ and $\xi\in\H_{0}$.  On the other hand, if
  \eqref{eq:49} holds for all $g\in\gcb$, $\xi\in\H_{0}$ and
  $\psi\in\ccg$, then % since
  % $M\bigl(\lambda(\psi\phi\bar\psi)\bigr)\ge0$,
  we must have $M\bigl(\lambda(\psi\phi\bar\psi)\bigr)=0$ for all
  $\psi\in\ccg$.  Since $\phi(x)\ge0$ everywhere, this forces
  $|\psi(x)|^{2}\phi(x)=0$ for $\nu$-almost all $x$.  Since $\psi$ is
  arbitrary, we conclude that $\phi(x)=0$ for $\nu$-almost all $x$.
  Therefore it will suffice to show that
  \begin{multline}
    \label{eq:50}
    L_{\xi,\xi}\bigl(g^{*}*\bigl(\lambda(\psi\tilde\phi\bar\psi)\circ
    r\bigr)\cdot g\bigr)=0 \\ \quad\text{for all $g\in\gcb$,
      $\xi\in\H_{0}$ and $\psi\in\ccg$,}
  \end{multline}
  where we have replaced $\phi$ with $\tilde \phi$ in the right-hand
  side of \eqref{eq:49}.  First, we compute that with $\phi$ in
  \eqref{eq:49} we have
  \begin{align}
    g^{*}*\bigl(\lambda(\psi\phi&\bar\psi)\circ r\bigr)\cdot g(z) =
    \int_{G}g(x^{-1})^{*}\bigl(\lambda(\psi\phi\bar\psi)\circ r\bigr)
    \cdot g(x^{-1}z) \,d\lambda^{r(z)}(x)\notag \\
    &= \int_{G } g(x^{-1})^{*}
    \lambda(\psi\phi\bar\psi)\bigl(s(x)\bigr)
    g(x^{-1}z) \,d\lambda^{r(z)}(x)\notag\\
    &=\int_{G}\int_{G} g(x^{-1})^{*} \overline{\psi(y)}\phi(y)\psi(y)
    g(x^{-1}z) \,d\lambda^{s(x)}(y)\,d\lambda^{r(z)}(x) \notag\\
    \intertext{which, after sending $y\mapsto x^{-1}y$ and using
      left-invariance of the Haar system, is} &=\int_{G}\int_{G}
    g(x^{-1})^{*} \overline{\psi(x^{-1}y)} \phi(x^{-1}y) \psi(x^{-1}y)
    g(x^{-1}z)
    \,d\lambda^{r(z)}(y)\,d\lambda^{r(z)}(x) \notag\\
    \intertext{which, after defining $F(x,y):=\psi(x^{-1}y)g(x^{-1})$
      and $\phi\cdot F(x,y):=\phi(x^{-1}y)F(x,y)$ for $(x,y)\in
      G\starr G$, is} &=\int_{G}\int_{G}F(x,y)^{*}\phi\cdot
    F(z^{-1}x,z^{-1}y) \,d\lambda^{r(z)}(y) \,
    d\lambda^{r(z)}(x).\label{eq:51}
  \end{align}

  We will have to look at vector-valued integrals of the form
  \eqref{eq:51} in some detail.  Define $\iota:G\starr G\to G$ by
  $\iota(x,y)=x^{-1}$, and let $\iota^{*}\B =
  \set{(x,y,b):p(b)=x^{-1}}$ be the pull-back bundle.  If
  $\psi\in\ccg$ and $g\in\gcb$, then
  \begin{equation*}
    (x,y)\mapsto \psi(y)g(x^{-1})
  \end{equation*}
  is a section in $\sa_{c}(G\starr G;\iota^{*}\B)$, and it is not hard
  to see that such sections span a subspace dense in the
  \ilt.\footnote{Our by now standard partition of unity argument
    (Lemma~\ref{lem-dense-sections}) shows that sections of
    the form
    \begin{equation*}
      (x,y)\mapsto \theta(x,y)g(x^{-1})
    \end{equation*}
    for $\theta\in C_{c}(G\starr G)$ and $g\in\gcb$ span a dense
    subspace.  We can approximate $\theta$ by sums of the form
    $\psi_{1}(x)\psi_{2}(y)$.  But $\psi_{1}(x)\psi_{2}(y)g(x^{-1}) =
    \psi_{2}(y)\tilde\psi_{1}\cdot g(x^{-1})$.}

  \begin{lemma}
    \label{lem-sp-F-one}
    Suppose that $F_{1},F_{2}\in\sa_{c}(G\starr G;\iota^{*}\B)$.  Then
    \begin{equation*}
      z\mapsto
      \int_{G}\int_{G}F_{1}(x,y)^{*}F_{2}(z^{-1}x,z^{-1}y)\,d\lambda^{r(z)}(y)
      \, d\lambda^{r(z)}(x)
    \end{equation*}
    defines a section in $\gcb$ which we denote by
    $F_{1}^{*}*_{\lambda*\lambda}F_{2}$.
  \end{lemma}
  \begin{proof}
    Let $K=K^{-1}$ be a compact set in $G$ such that $\supp
    F_{i}\subset K\times K$.  Then $\supp
    F_{1}*_{\lambda*\lambda}F_{2}\subset K^{2}$, and
    \begin{equation*}
      \|F_{1}^{*}*_{\lambda*\lambda}F_{2}\|_{\infty}
      \le\|F_{1}\|_{\infty}\|F_{2}\|_{\infty} \lambda^{u} (K)^{2}. 
    \end{equation*}
    Thus it suffices to show the result for $F_{i}$ which span a dense
    subspace in the \ilt.  Therefore we may as well assume that
    $F_{i}(x,y)=\psi_{i}(y)g_{i}(x^{-1})$ as above.  Then
    \begin{align*}
      F_{1}^{*}*_{\lambda*\lambda}F_{2}(z)&=\int_{G}\int_{G}
      \overline{\psi_{1}(y)} \psi_{2}(z^{-1}y)
      g_{1}(x^{-1})^{*}g_{2}(x^{-1}z)
      \,d\lambda^{r(z)}(y)\,d\lambda^{r(z)}(x) \\
      &= \bar\psi_{1}*\tilde\psi_{2}(z)g_{1}^{*}*g_{2}(z),
    \end{align*}
    and the result follows by Corollary~\ref{cor-conv}.
  \end{proof}

  \begin{lemma}
    \label{lem-dense}
    If $\psi\in\ccg$ and $g\in\gcb$, then
    \begin{equation*}
      (x,y)\mapsto \psi(x^{-1}y)g(x^{-1})
    \end{equation*}
    is a section in $\sa_{c}(G\starr G;\iota^{*}\B)$ and sections of
    this form space a dense subspace in the \ilt.
  \end{lemma}
  \begin{proof}
    It suffices to see that we can approximate sections of the form
    \begin{equation}
      \label{eq:52}
      (x,y)\mapsto \theta(x,y)f(x^{-1})
    \end{equation}
    with $\theta\in C_{c}(G\starr G)$ and $f\in\gcb$.  A
    Stone-Weierstrass argument shows that we can approximate $\theta$
    with sums of functions of the form $(x,y)\mapsto
    \psi(x^{-1}y)\psi'(x^{-1})$.  Then we can approximate
    \eqref{eq:52} by sums of sections of the required form:
    $(x,y)\mapsto \psi(x^{-1}y)g(x^{-1})$ where $g(x):=\psi'(x)f(x)$.
    This completes the proof.
  \end{proof}

  Let $\AA\subset\sa_{c}(G\starr G;\iota^{*}\B)$ be the dense subspace
  of sections of the form considered in Lemma~\ref{lem-dense}.  We
  continue to write $\phi$ for the characteristic function of our
  fixed pre-compact, $\nu$-null set.  Then we know from \eqref{eq:49}
  that
  \begin{equation}
    \label{eq:53}
    L_{\xi,\xi}\bigl(F^{*}*_{\lambda*\lambda}(\phi\cdot
    F)\bigr)=0\quad\text{for all $F\in\AA$.}
  \end{equation}
  It is not hard to check that, if $\phi'\in\bboc(G)$, then
  $F^{*}*_{\lambda*\lambda}(\phi'\cdot F)\in \sgcb$ and that if
  $F_{n}\to F$ in the \ilt{} in $\sa_{c}(G\starr G;\iota^{*}\B)$, then
  $\set{F_{n}^{*}*_{\lambda*\lambda}(\phi'\cdot F)}$ is uniformly
  bounded and converges pointwise to $F*_{\lambda*\lambda}(\phi'\cdot
  F)$.  In particular the continuity of the $L_{\xi,\xi}$ (see
  Lemma~\ref{lem-grm-ext}) implies that \eqref{eq:53} holds for all
  $F\in\sa_{c}(G\starr G;\iota^{*}\B)$.  But if we define $\tilde
  F(x,y):=F(y,x)$, then we see from the definition that
  \begin{equation*}
    \tilde F^{*}*_{\lambda*\lambda}(\phi\cdot \tilde F)=
    F^{*}*_{\lambda*\lambda}(\tilde \phi\cdot F)),
  \end{equation*}
  where we recall that $\tilde \phi(x):=\phi(x^{-1})$.  Thus
  \begin{equation*}
    L_{\xi,\xi}\bigl(F^{*}*_{\lambda*\lambda}(\tilde \phi\cdot F))\bigr)
    =0 \quad\text{for all $F\in\sa_{c}(G\starr G;\iota^{*}\B)$}.
  \end{equation*}
  Since the above holds in particular for $F\in\AA$, this implies
  \eqref{eq:50} and completes the proof.
\end{proof}

We can now turn our attention to creating the Borel Hilbert bundle.
We still need some ``Borelogy'' in the form of 
Lemma~\ref{lem-paul-tech}.

\begin{lemma}
  \label{lem-abs-cont}
  Let $a$ and $b$ be vectors in $\hoo$ (identified with
  $\gcbatho/\N$).  Let $L_{a,b}$\index{Lab@$L_{a,b}$} be the \grm{}
  given by
  \begin{equation}\label{eq:54}
    L_{a,b}(f):=\bip(L(f)a|b).
  \end{equation}
  Then $|L_{a,b}|\ll\nu$, where $\nu$ is the measure on $G$ given by
  \eqref{eq:33}, and $|L_{a,b}|$ is the total variation of $L_{a,b}$
  as defined in the paragraph preceding\index{Lab@$"|L_{a,b}"|$}
  Proposition~\ref{prop-din28.32}.
\end{lemma}
\begin{remark}
  \label{rem-hoo}
  Although the \grm{} defined in \eqref{eq:54} makes perfectly good
  sense for any $a,b\in\H_{0}$, notice that we are only claiming the
  result of the proposition when $a,b\in\hoo$ (because that is all we
  are able to  prove).
\end{remark}

\begin{proof}
  It is enough to show that if $M$ is a pre-compact $\nu$-null set,
  then $|L_{a,b}|(M)=0$.  Since $\nu$ is a Radon measure, and
  therefore regular, we may as well assume that $M$ is a
  $G_{\delta}$-set.  Thus if $\phi:=\charfcn M$, then there are
  $\phi_{n}\in\ccpg$ such that $\phi_{n}(x)\searrow \phi(x)$ as in
  Lemma~\ref{lem-tol-var-borel} (which has been cooked up for just
  this purpose).  Then, using Lemma~\ref{lem-tol-var-borel}, it will
  suffice to see that $L_{a,b}(f)=0$ whenever $f\in\sgcb$ and
  $\|f\|\le\phi$.

  On the other hand,
  \begin{equation*}
    0=\int_{\go}\int_{G} \phi(x)\,d\lambda^{u}(x)\,d\mu(u),
  \end{equation*}
  so there is a $\mu$-null set $N\subset\go$ such that
  $\lambda^{u}(M\cap G^{u})=0$ if $u\notin N$.  As above, we can
  assume that $N$ is a $G_{\delta}$ set.  Thus if $f\in\sgcb$ is such
  that $\|f\|\le\phi$ and if $g\in\gcb$, then
  \begin{equation*}
    f*g(x)=\int_{G}f(y)g(y^{-1}x)\,d\lambda^{r(x)}(y)=0
  \end{equation*}
  whenever $r(x)\notin N$.  Since $\supp \lambda^{r(x)}=G^{r(x)}$, it
  follows that for all $x\in G$ (without exception),
  \begin{equation}\label{eq:55}
    f*g(x)= \charfcn N\bigl(r(x)\bigr) f*g(x) = \bigl((\charfcn N\circ
    r)\cdot f\bigr)*g(x).  
  \end{equation}

  Since $a,b\in\hoo$, it suffices to consider $a=[g\tensor\xi]$ and
  $b=[h\tensor \eta]$ (with $g,h\in\gcb$ and $\xi,\eta\in\H_{0})$.
  Note that $f$ and $\charfcn N$ satisfy the hypotheses of
  Lemma~\ref{lem-paul-tech}.  Therefore, by part~\partref1 of that
  lemma,
  \begin{align*}
    L_{[g\tensor\xi,h\tensor\eta}(f)&=\bip([f*g\tensor
    \xi]|[h\tensor\eta]) \\
    \intertext{which, by \eqref{eq:55}, is}
    &= \bip([((\charfcn N\circ r)\cdot f)*g\tensor\xi]|[h\tensor\eta]) \\
    \intertext{which, by part~\partref1 of Lemma~\ref{lem-paul-tech},
      is} &= \bip(L_{b}((\charfcn N\circ r)\cdot
    f)[g\tensor\xi]|[h\tensor\eta])
    \\
    \intertext{which, by part~\partref3 of Lemma~\ref{lem-paul-tech},
      is}
    &= \bip(M(\charfcn N)L_{b}(f)[g\tensor\xi]|[h\tensor\eta]) \\
    \intertext{which, since $M(\charfcn N)=0$, is} &=0,
  \end{align*}
  as desired.  This completes the proof.
\end{proof}

Since the measures $\nu$ and $\nu_{0}$ are equivalent, for each
$\xi,\eta\in\hoo$, we can, in view of Lemma~\ref{lem-abs-cont}, let
$\rho_{\xi,\eta}$\index{rho@$\rho_{\xi,\eta}$} be the Radon-Nikodym
derivative of $|L_{\xi,\eta}|$ 
with respect to $\nu_{0}$.  Then for each $\xi,\eta\in\hoo$, we have
\begin{align*}
  \bip(&L(f)\xi|\eta)=
  L_{\xi,\eta}(f) \\
  &= \int_{G} \epsilon(\xi,\eta)_{x}\bigl(f(x)\bigr)
  \,d|L_{\xi,\eta}|(x)
  \\
  &= \int_{G} \epsilon(\xi,\eta)_{x}\bigl(f(x)\bigr)
  \rho_{\xi,\eta}(x) \Delta(x)^{-\half} \,d\nu(x)
  \\
  &=\int_{\go}\int_{G}\epsilon(\xi,\eta)_{x}\bigl(f(x)\bigr)
  \rho_{\xi,\eta}(x) \Delta(x)^{-\half}\,d\lambda^{u}(x)\,d\mu(u) ,
\end{align*}
where, of course, $\epsilon(\xi,\eta)_{x}$ denotes a linear functional
in the unit ball of $B(x)^{*}$ which depends on the choice of $\xi$
and $\eta$ in $\hoo$.

\begin{remark}
  \label{rem-uniqueness}
  It is surely the case that there are interesting uniqueness
  conditions satisfied by the
  $\rho_{\xi,\eta}$ and the $\epsilon(\xi,\eta)_{x}$
  that would make our subsequent calculations a bit tidier.  That is,
  we would expect that, in some ``almost
  everywhere'' sense, $\epsilon(\xi,\eta)$ is linear in $\xi$ and
  conjugate linear in $\eta$.  A similar statement should hold for
  $\rho_{\xi,\eta}$.  If true, this would make defining an
  inner product in Lemma~\ref{lem-tensor-fix} more straightforward. 
We will finesse these issues below by restricting to a countable set
  of $\xi$'s and $\eta$'s.
\end{remark}

Our next computation serves to motivate the construction in
Lemma~\ref{lem-tensor-fix}.  We need $\xi,\eta\in\hoo$ to be able to
apply Lemma~\ref{lem-abs-cont}.  To simplify the notation, we write
$\epsilon$ in place of $\epsilon(\xi,\eta)$ and $\rho$ in place of
$\rho_{\xi,\eta}$.  Then
\begin{align*}
  \bip(L(f)&\xi|L(g)\eta)= \bip(L(g^{*}*f)\xi|\eta) =L_{\xi,\eta}(g^{*}*f) \\
  &= \int_{\go}\int_{G} \epsilon_{x}\bigl(g^{*}*f(x)\bigr) \rho(x)
  \Delta(x)^{-\half} \,d\lambda^{u}(x)\,d\mu(u) \\
  &= \int_{\go}\int_{G}\int_{G}
  \epsilon_{x}\bigl(g(y^{-1})^{*}f(y^{-1}x)\bigr)
  \rho(x)\Delta(x)^{-\half}
  \,d\lambda^{u}(u)\,d\lambda^{u}(x)\,d\mu(u)
   \\
  \intertext{which, by Fubini\footnote{Empty} and sending $x\mapsto
    yx$, is} &= \int_{\go}\int_{G}\int_{G}
  \epsilon_{yx}\bigl(g(y^{-1})^{*}f(x)\bigr) \rho(yx)
  \Delta(yx)^{-\half} \, d\lambda^{s(y)}(x) \,d\lambda^{u}(y)\,d\mu(u)
 \\
  \intertext{which, after sending $y\mapsto y^{-1}$, and using the
    symmetry of $\nu_{0}$, is} &= \int_{\go}\int_{G}\int_{G}
  \epsilon_{y^{-1}x}\bigl(g(y)^{*}f(x)\bigr) \rho(y^{-1}x)
  \Delta(y)^{-\half} \Delta(x)^{-\half} \\
  &\hskip2.5in \,d\lambda^{u}(x)\,d\lambda^{u}(y)\, d\mu(u).
\end{align*}

\footnotetext{The integrand is Borel by Lemma~\ref{lem-prod-borel}, so
  there is no problem applying Fubini's Theorem here.}  It will be
convenient to establish some additional notation.  We fix once and for
all a countable orthonormal basis
$\set{\spxi_{i}}$\index{zeta@$\spxi_{i}$} for $\hoo$.  (Actually, any
countable linearly independent set whose span is dense in $\hoo$ will
do.)  We let\index{Hoop@$\hoop$}
\begin{equation*}
  \hoop:=\operatorname{span}\set{\spxi_{i}}.
\end{equation*}
To make the subsequent formulas a bit easier to read, we will write
$\rho_{ij}$\index{rho@$\rho_{ij}$} in place of the Radon-Nikodym
derivative $\rho_{\spxi_{i},\spxi_{j}}$ and
$\eij_{x}$\index{epsilon@$\eij_{x}$} in place of the linear functional
$\epsilon(\spxi_{i},\spxi_{j})_{x}$.  The linear independence of the
$\spxi_{i}$ guarantees that each $\alpha\in\gcb\atensor\hoo'$ can be
written \emph{uniquely} as
\begin{equation*}
  \alpha=\sum_{i} f_{i}\tensor\spxi_{i}
\end{equation*}
where all by finitely many $f_{i}$ are zero.

\begin{lemma}
  \label{lem-tensor-fix}
  For each $u\in\go$, there is a sesquilinear form $\ipu<\cdot,\cdot>$
  on $\gcbathoop$ such that
  \begin{multline}
    \label{eq:56}
    \ipu<f\tensor\spxi_{i},g\tensor \spxi_{j}>\\ = \int_{G}\int_{G}
    \eij_{y^{-1}x} \bigl(g(y)^{*}f(x)\bigr) \rho_{ij}(y^{-1}x)
    \Delta(yx)^{-\half} \,d\lambda^{u}(x)\,d\lambda^{u}(y).
  \end{multline}
  Furthermore, there is a $\mu$-conull set $F\subset \go$ such that
  $\ipu<\cdot,\cdot>$ is a pre-inner product for all $u\in F$.
\end{lemma}
\begin{remark}
  \label{rem-tensor-prob}
  As we noted in Remark~\ref{rem-uniqueness}, we fixed the $\spxi_{i}$
  because it isn't clear that the right-hand side of \eqref{eq:56} is
  linear in $\spxi_{i}$ or conjugate linear in $\spxi_{j}$.
\end{remark}

\begin{proof}
  The integrand in \eqref{eq:56} has compact support and
is Borel by
  Lemma~\ref{lem-prod-borel}.  Therefore the integrals there and
  below are well-defined.  Given $\alpha=\sum_{i}f_{i}\tensor
  \spxi_{i}$ and $\beta=\sum_{j} g_{j}\tensor \spxi_{j}$, we get a
  well-defined form via the definition
  \begin{equation*}
    \ipu<\alpha,\beta> = 
    \sum_{ij} \int_{G}\int_{G} \eij_{y^{-1}x}\bigl(g_{j}(y)^{*}f_{i}(x)\bigr)
    \rho_{ij}(y^{-1}x) \Delta(yx)^{-\half}
    \,d\lambda^{u}(x)\,d\lambda^{u}(y).  
  \end{equation*}
  This clearly satisfies \eqref{eq:56}, and the linearity of the
  $\eij_{z}$ imply that $(\alpha,\beta)\mapsto \ipu<\alpha,\beta>$ is
  linear in $\alpha$ and conjugate linear in $\beta$.  It only remains
  to provide a conull Borel set $F$ such that $\ipu<\cdot,\cdot>$ is
  positive for all $u\in F$.

  However, \eqref{eq:56} was inspired by the calculation preceding the
  lemma.  Hence if $\alpha:=\sum_{i}f_{i}\tensor\spxi_{i}$, then
  \begin{equation}\label{eq:57}
    \begin{split}
      \Bigl\|\sum L(f_{i})\spxi_{i}\Bigr\|^{2}&=
      \sum_{ij} \bip(L(f_{i})\spxi_{i}|{L(f_{j})\spxi_{j}}) \\
      &= \sum_{ij} \bip(L(f_{j}^{*}*f_{i})\spxi_{i}|\spxi_{j}) \\
      &= \sum_{ij} \int_{\go}\int_{G} \int_{G} \eij_{y^{-1}x}\bigl(
      f_{j}(y)^{*}f_{i}(x) \bigr) \rho_{ij}(y^{-1}x) \Delta(xy)^{-\half}\\
      &\hskip2in
      \,d\lambda^{u}(x) \, d\lambda^{u}(y) \,d\mu(u) \\
      &= \sum_{ij} \int_{\go} \ipu<f_{i}\tensor \spxi_{i},f_{j}\tensor
      \spxi_{j}>\,d\mu(u) \\
      &= \int_{\go} \ipu<\alpha,\alpha>\,d\mu(u).
    \end{split}
  \end{equation}
  Thus, for $\mu$-almost all $u$, we have $\ipu<\alpha,\alpha>\ge0$.
  The difficulty is that the exceptional null set depends on $\alpha$.
  However, since we have assumed $B:=\sa_{0}(G;\B)$ is separable,
  there is a sequence $\set{f_{i}}\subset\gcb$ which is dense in
  $\gcb$ in the \ilt. %\altdb\footnote{Do we need a reference for this?}
  Let $\AA$ be the rational vector space spanned by the countable set
  $\set{f_{i}\tensor\spxi_{j}}_{i,j}$.  Since $\AA$ is countable,
  there is a $\mu$-conull set $F$ such that $\ipu<\cdot,\cdot>$ is a
  positive $\mathbf{Q}$-sesquilinear form on $\AA$.  However, if
  $g_{i}\to g$ and $h_{i}\to h$ in the \ilt{} in $\gcb$, then, since
  $\lambda^{u}\times\lambda^{u}$ is a Radon measure on $G\times G$, we
  have $\ipu<g_{i}\tensor\spxi_{j}, h_{i}\tensor \spxi_{k} > \to
  \ipu<g\tensor\spxi_{j},h\tensor\spxi_{k}>$.  It follows that for all
  $u\in F$, $\ipu<\cdot,\cdot>$ is a positive sesquilinear form (over
  $\C$) on the complex vector space generated by
  \begin{equation*}
    \set{f\tensor \spxi_{i}:f\in\gcb}.
  \end{equation*}
  However, as that is all of $\gcbathoop$, the proof is complete.
\end{proof}

Note that for \emph{any} $u\in\go$, the value of
$\ipu<f\tensor\spxi_{i},g \tensor\spxi_{j}>$ depends only on $f\restr
{G^{u}}$ and $g\restr{G^{u}}$.
Furthermore, using a suitable vector-valued Tietze Extension Theorem
(see Proposition~\ref{prop-tietze}), we can view $\ipu<\cdot,\cdot>$
as a sesquilinear form on $\gucb$.  (Clearly, each $f\in\gcb$
determines a section of $\gucb$.  We need the extension theorem to
know that every section in $\gucb$ arises in this fashion.)

Using our Tietze Extension Theorem as above, given $f\in\gcb$ and
$b\in\B$ there is a section, denoted by $\hatpi(b)f$ such
that\footnote{Recall that $r(b)$ is a shorthand for $r(p(b))$, and
  similarly for $s(b)$.}
\begin{equation*}
  (\hatpi(b)f)(x)=\Delta(z)^{\half}bf(z^{-1}x)\quad\text{for all 
    $x\in G^{r(b)}$.}
\end{equation*}
Of course, $\hatpi(b)f$ is only well-defined on $G^{r(b)}$.  Then, if
$b\in\B$, we can compute that
\begin{align*}
  \bipb r(b)<\hatpi(b)f\tensor\spxi_{i},g\tensor\spxi_{j}&> =
  \int_{G}\int_{G} \eij_{y^{-1}x}\bigl(g(y)^{*}bf(z^{-1}x)\bigr)
  \rho_{ij} (x^{-1}y) \Delta(z^{-1}xy)^{-\half}\\
  &\hskip2.25in \, d\lambda^{r(b)}(y) \,
  d\lambda^{r(b)} (x) \\
  \intertext{which, after sending $x\mapsto zx$, is} &=
  \int_{G}\int_{G} \eij_{y^{-1}zx}\bigl(g(y)^{*}bf(x)\bigr)
  \rho_{ij}(x^{-1}z^{-1}y) \Delta(xy)^{-\half} \\
  &\hskip2.25in
  \,d\lambda^{r(b)}(y)\,d\lambda^{s(b)}(x) \\
  \intertext{which, after sending $y\mapsto zy$, is} &=
  \int_{G}\int_{G} \eij_{y^{-1}x}\bigl(g(zy)^{*}bf(x)\bigr)
  \rho_{ij}(x^{-1}y) \Delta(z)^{-\half}\Delta(xy)^{-\half}\\
  &\hskip2.25in \,
  d\lambda^{s(b)}(y) \,d\lambda^{s(b)}(x) \\
  &= \int_{G}\int_{G}
  \eij_{y^{-1}x}\bigl(\bigl(\hatpi(b^{*})g\bigr)(y)f(x)\bigr)
  \Delta(xy)^{-\half} \,
  d\lambda^{s(b)}(y) \,d\lambda^{s(b)}(x) \\
  &= \bipb s(b)<f\tensor\spxi_{i},\hatpi(b^{*})g\tensor\eta>.
\end{align*}
In particular, if $a\in B(u)$, then $(\hatpi(a)f)(x)=af(x)$ and
\begin{equation}
  \label{eq:58}
  \bipu<\hatpi(a)f\tensor\spxi_{i},g\tensor\spxi_{j}> = \bipu<f\tensor
  \spxi_{i}, \hatpi(a^{*})g\tensor\spxi_{j}>.
\end{equation}

Recall that $G$ acts continuously on the left of $\go$: $x\cdot
s(x)=r(x)$.  In particular, if $C$ is compact in $G$ and if $K$ is
compact in $\go$, then
\begin{equation*}
  C\cdot K=\set{x\cdot u:(x,u)\in G^{(2)}\cap(C\times K)}
\end{equation*}
is compact.  If $U\subset \go$, then we say that $U$ is
\emph{saturated}\index{saturated} if $U$ is $G$-invariant.  More
simply, $U$ is saturated if $s(x)\in U$ implies $r(x)$ is in $U$.  If
$V\subset\go$, then its \emph{saturation} is the set $[V]=G\cdot
V$\index{V@$[V]$} which is the smallest saturated set containing $V$.

The next result is a key technical step in our proof and takes the
place of the \index{Ramsay selection theorems}Ramsay selection theorems
(\citelist{\cite{ram:jfa82}*{Theorem~3.2}\cite{ram:am71}*{Theorem~5.1}})
used in Muhly's and Renault's proofs.

\begin{lemma}
  \label{lem-saturated}
  We can choose the $\mu$-conull Borel set $F\subset\go$ in
  Lemma~\ref{lem-tensor-fix} to be saturated for the $G$-action on
  $\go$.
\end{lemma}
\begin{proof}
  Let $F$ be the Borel set from Lemma~\ref{lem-tensor-fix}.  We want
  to see that $\ipv<\cdot,\cdot>$ is positive for all $v$ in the
  saturation of $F$.  To this end, suppose that $u\in F$ and that
  $z\in G$ is such that $s(z)=u$ and $r(z)=v$.  If $b\in B(z)$, then
  \begin{equation*}
    x\mapsto \Delta(z)^{\half} b f(z^{-1}x)
  \end{equation*}
  is a section in $\gucb$, and an application of
  Lemma~\ref{lem-dense-sections} shows that such sections span a
  dense subspace of $\gucb$ in the \ilt.  Moreover, as we observed at
  the end of the proof of Lemma~\ref{lem-tensor-fix},
  \begin{equation*}
    \ipv<f_{i}\tensor\spxi_{j},g_{i}\tensor\spxi_{k}>\to
    \ipv<f\tensor\spxi_{j} , g\tensor\spxi_{k}>
  \end{equation*}
  provided $f_{i}\to f$ and $g_{i}\to g$ in the \ilt{} in
  $\sa_{c}(G^{v};\B\restr{G^{v}})$.  Therefore, to show that
  $\ipv<\cdot,\cdot>$ is positive, it will suffice to check on vectors
  of the form $\alpha:=\sum_{i} \hatpi(b_{i})(f_{i})\tensor
  \spxi_{i}$.  Then using the calculation preceding \eqref{eq:58}, we
  have
  \begin{equation} \label{eq:59} \ipv<\alpha,\alpha>=\sum_{ij}
    \bipu<\hatpi(b_{j}^{*}b_{i})f_{i}\tensor\spxi_{i} ,
    f_{j}\tensor\spxi_{j}>.
  \end{equation}
  However, since $B(z)$ is, in particular, a right Hilbert
  $B(u)$-module with inner product $\rip
  B(u)<b_{j},b_{i}>=b_{j}^{*}b_{i}$ (by part~\partref5 of
  Definition~\ref{def-fell-bund}), the matrix $(b_{j}^{*}b_{i})$ is
  positive in $M_{n}\bigl(B(u)\bigr)$ by
  \cite{rw:morita}*{Lemma~2.65}.  Therefore there are $d_{rs}\in B(u)$
  such that $b_{j}^{*}b_{i}=\sum_{k} d_{kj}^{*}d_{ki}$.  Then, using
  \eqref{eq:59}, the right-hand side of \eqref{eq:59} is
  \begin{align*}
    \sum_{ijk} \bipu<\hatpi(d_{kj}^{*}d_{ki})f_{i}\tensor\spxi_{i} ,
    f_{j}\tensor \spxi_{j}> &= \sum_{ijk} \bipu <
    \hatpi(d_{ki})f_{i}\tensor\spxi_{i},
    \hatpi(d_{kj}) f_{j}\tensor\spxi_{j}> \\
    &= \sum_k \bipu< \sum_{i} \hatpi(d_{ki})f_{i}\tensor\spxi_{i},
    \sum_{i} \hatpi(d_{ki})f_{i}\tensor\spxi_{i}>
  \end{align*}
  which is positive since $u\in F$.

  It only remains to verify that the saturation of $F$ is Borel.
  Since $\mu$ is a Radon measure --- and therefore regular --- we can
  shrink $F$ a bit, if necessary, and assume it is $\sigma$-compact.
  Say $F=\bigcup K_{n}$.  On the other hand, $G$ is second countable
  and therefore $\sigma$-compact.  If $G=\bigcup C_{m}$, then
  $[F]=\bigcup C_{m}\cdot K_{n}$.  Since each $C_{m}\cdot K_{n}$ is
  compact, $[F]$ is $\sigma$-compact and therefore Borel.  This
  completes the proof.
\end{proof}

From here on, we will assume that $F$ is saturated.  In view of
Lemma~\ref{lem-tensor-fix}, for each $u\in F$ we can define $\H(u)$
to be the Hilbert space completion of $\gcbathoop$ with respect to
$\ipu<\cdot,\cdot>$.  We will denote the image of $f\tensor\spxi_{i}$
in $\H(u)$ by $f\tensor_{u}\spxi_{i}$.  Since the complement of $F$ is
$\mu$-null and also saturated, what we do off $F$ has little
consequence.  In particular, $G$ is the disjoint union of $G\restr F$
and the $\nu$-null set $G\restr{\go\setminus F}$.\footnote{The
  saturation of $F$ is critical to what follows.  If $F$ is not
  saturated, then in general $G$ is not the union of $G\restr F$ and $
  G\restr{\go\setminus F}$.  But as $F$ \emph{is} saturated, note that
  a homomorphism $\phi:G\restr F\to H$ can be trivially extended to a
  homomorphism on all of $G$ by letting $\phi$ be suitably trivial on
  $G\restr{\go\setminus F}$.  This is certainly not the case unless $F$ is
  saturated.}  Nevertheless, for niceties sake, we let
$\mathcal{V}$ be a Hilbert space with an orthonormal basis
$\set{e_{ij}}$ doubly indexed by the same index sets as for
$\set{f_{i}}$ and $\set{\spxi_{j}}$, and set $\H(u):=\mathcal{V}$ for
all $u\in\go\setminus F$.  We then let
\begin{equation*}
  \go*\HH=\set{(u,h):\text{$u\in F$ and $h\in\H(u)$}},
\end{equation*}
and define $\Phi_{ij}:F\to F*\HH$ by
\begin{equation*}
  \Phi_{ij}(u):=
  \begin{cases}
    f_{i}\tensor_{u}
    \spxi_{j}&\text{if $u\in F$ and} \\
    e_{ij}&\text{if $u\notin F$.}
  \end{cases}
\end{equation*}
(Technically, $\Phi_{ij}(u)=(u, f_{i}\tensor_{u}\spxi_{j})$ --- at
least for $u\in F$ --- but we have agreed to obscure this subtlety.)
Then \cite{wil:crossed}*{Proposition~F.8} implies that we can make
$\go*\HH$ into a Borel Hilbert bundle over $\go$ in such a way that
the $\set{\Phi_{ij}}$ form a fundamental sequence (see
\cite{wil:crossed}*{Definition~F.1}).  Note that if
$f\tensor\spxi_{i}\in \gcbathoop$ and if
$\Phi(u):=f\tensor_{u}\spxi_{i}$, then
\begin{equation*}
  u\mapsto \bipu<\Phi(u),\Phi_{ij}(u)>
\end{equation*}
is Borel on $F$.\footnote{We can define $\Phi(u)$ to be zero off $F$.
  We are going to continue to pay as little attention as possible to
  the null complement of $F$ in the sequel.}  It follows that $\Phi$
is a Borel section of $\go*\HH$ and defines a class in
$L^{2}(\go*\HH,\mu)$.

We can extend $\hatpi$ so that \eqref{eq:58} holds for all $a\in
B(u)^{\sim}$.  If $p(b)\in G^{v}_{u}$, then $\|b\|^{2}1_{B(u)}-b^{*}b$
is a positive element in the \cs-algebra $B(u)$.  Therefore, there is
a $k\in B(u)^{\sim}$ such that
\begin{equation*}
  \|b\|^{2}1_{B(u)}-b^{*}b=k^{*}k.
\end{equation*}
Then, using \eqref{eq:58} and the computation that preceded it, we
have
\begin{align*}
  \|b\|^{2}\bipu<\sum_{i} f_{i}\tensor\spxi_{i}, \sum_{j}
  f_{j}\tensor\spxi_{j}>&{} - \bipv<\sum_{i}
  \hatpi(b)f_{i}\tensor\spxi_{i} , \sum_{j} \hatpi(b)
  f_{j}\tensor\spxi_{j}>
  \\
  &=\sum_{ij} \bipu< \hatpi\bigl(\|b\|^{2}1-b*b\bigr)f_{i}\tensor
  \spxi_{i}, f_{j}\tensor \spxi_{j}> \\
  &=\sum_{ij}
  \bipu<\hatpi(k)f_{i}\tensor\spxi_{i},\hatpi(k)f_{j}\tensor
  \spxi_{j}> \\
  &= \bipu<\sum_{i} \hatpi(k)f_{i}\tensor\spxi_{i}, \sum_{j}
  \hatpi(k)f_{j}\tensor \spxi_{j}> \\
  &\ge 0.
\end{align*}
In other words, if we define $\pi(b)$ by
\begin{equation*}
  \pi(b)\Bigl(\sum_{i}f_{i}\tensor\spxi_{i}\Bigr)=
  \sum_{i}\hatpi(b)f_{i}\tensor \spxi_{i},
\end{equation*}
then
\begin{equation}\label{eq:60}
  \bipb r(b)<\pi(b)(\alpha),\pi(b)(\alpha)> \le
  \|b\|^{2}\bipb s(b)<\alpha,\alpha> \quad\text{for $b\in
    p^{-1}(G\restr F)$.} 
\end{equation}
Therefore we get a bounded operator, also denoted by $\pi(b)$, from
$\H(u)$ to $\H(v)$ with $\|\pi(b)\|\le\|b\|$.  Since \eqref{eq:60}
implies that $\pi(b)$ takes vectors of length zero to vectors of
length zero, we have
\begin{equation*}
  \pi(b)(f\tensor_{s(b)}\spxi_{i})=
  \hatpi(b)f\tensor_{r(b)}\spxi_{i}\quad\text{for all $b\in p^{-1}(G\restr F)$}.
\end{equation*}

If $b\notin p^{-1}(G\restr F)$, then $\H(s(b))=\H(r(b))=\mathcal V$,
and we can let $\pi(b)$ be the identity operator.

\begin{lemma}
  \label{lem-pi}
  The map $\hat\pi$ from $\B$ to $ \End(\go*\HH)$ defined by
  $\hat\pi(b):= \bigl(r(b), \pi(b), s(b)\bigr)$ is a Borel
  $*$-functor.  Hence $(\mu,\go*\HH,\hat\pi)$ is a strict
  representation of $\B$ on $L^{2}(\go*\HH,\mu)$.
\end{lemma}
\begin{proof}
  If $f\in\gcb$ and if $z\in G\restr F$, then
  \begin{multline*}
    \bip(\pi\bigl(f(z)\bigr)\Phi_{ij}\bigl(s(z)\bigr)
    |{\Phi_{kl}\bigl(r(z)\bigr) }) = \\
    \int_{G}\int_{G}
    \epsilon^{jl}_{y^{-1}x}\bigl(f_{k}^{*}(y)f(z)f_{i}(z^{-1}x) \bigr)
    \rho_{jl}(y^{-1}x)\Delta(z^{-1}xy)^{-\half} \, d\lambda^{r(z)}(y)\,
    d\lambda^{r(z)}(x).
  \end{multline*}
  Thus $z\mapsto \bip(\pi\bigl(f(z)\bigr)\Phi_{ij}\bigl(s(z)\bigr)
  |{\Phi_{kl}\bigl(r(z)\bigr) })$ is Borel on $F$ by
  Lemma~\ref{lemma-kappa-star} (in the form of
  Example~\ref{ex-b-borel}) and Fubini's Theorem.  Since it is
  clearly Borel on the complement of $F$, $\hat\pi$ satisfies the
  Borel condition in Definition~\ref{def-star-functor}.\footnote{It
    suffices to check on a fundamental sequence.}  It is straightforward to verify the algebraic
  properties (that is, properties
  \partref1, \partref2 and \partref3 of
  Definition~\ref{def-star-functor}).  For example, assuming that
  $x\in G^{r(a)}$, we have on the one hand,
  \begin{equation*}
    \bigl(\hatpi(ab)f\bigr)(x)=\Delta\bigl(p(ab)\bigr)^{\half} ab
    f\bigl( p(ab)^{-1}x\bigr),
  \end{equation*}
  while
  \begin{align*}
    \bigl(\hatpi(a)\hatpi(b)f\bigr)(x)&=
    \Delta\bigl(p(a)\bigr)^{\half}a\bigl(\hatpi(b)f\bigr)\bigl(p(a)^{-1}x\bigr) \\
    &= \Delta\bigl(p(a)p(b)\bigr)^{\half} ab
    f\bigl(p(b)^{-1}p(a)^{-1}x\bigr).
  \end{align*}
  Since $p(ab)=p(a)p(b)$, it follows that $\hat\pi$ is multiplicative
  on $p^{-1}(G\restr F)$.  Of course, it is clearly multiplicative on
  the complement (which is $p^{-1}(G\restr{\go\setminus F})$ since $F$
  is saturated).  The other properties follow similarly.
\end{proof}

\begin{lemma}
  \label{lem-V}
  Each $f\tensor\spxi_{i}\in\gcbathoop$ determines a Borel section
  $\Phi(u):=f\tensor_{u}\spxi_{i}$ whose
  class in $L^{2}(\go*\HH,\mu)$ depends only on the class of
  $[f\tensor\spxi_{i}]\in \gcbathoop/\N\subset \gcbatho/\N=\hoo$.
  Furthermore, there is a unitary isomorphism $V$ of $\H$ onto
  $L^{2}(\go*\HH,\mu)$ such that $V(L(f)\spxi_{i})=[\Phi]$.
\end{lemma}
\begin{proof}
  We have already seen that $\Phi$ is in $L^{2}(F*\HH,\mu)\cong
  L^{2}(\go*\HH,\mu)$.  More generally, the computation \eqref{eq:57}
  in the proof of Lemma~\ref{lem-tensor-fix} shows that if
  $\alpha=\sum_{i} f_{i}\tensor\spxi_{i}$ and $\Psi(u):=\sum_{i}
  f_{i}\tensor_{u}\spxi_{i}$, then
  \begin{equation*}
    \|\Psi\|_{2}^{2}=\Bigl\|\sum_{i} L(f_{i})\spxi_{i}\Bigr\|^{2}
  \end{equation*}
  Thus there is a well defined isometric map $V$ as in the statement
  of lemma mapping $\operatorname{span}\set{L(f)\spxi_{i}:f\in\gcb}$
  onto a dense subspace of $L^{2}(F*\HH,\mu)$.  Since $\hoop$ is dense
  in $\hoo$, and therefore in $\H$, the result follows by
  Corollary~\ref{cor-dense}.
\end{proof}

The proof of Theorem~\ref{thm-fell-disintegration} now follows almost
immediately from the next proposition.

\begin{prop}
  \label{prop-main}
  The unitary $V$ defined in Lemma~\ref{lem-V} intertwines $L$ with a
  representation $L'$ which in the integrated form of the strict
  representation $(\mu,\go*\HH,\hat\pi)$ from Lemma~\ref{lem-pi}.
\end{prop}
\begin{proof}
  We have $L'(f_{1})=VL(f_{1})V^{*}$.  On the one hand,
  \begin{multline*}
    \bip(L(f_{1})[f\tensor\spxi_{i}]|[g\tensor \spxi_{j}])_{\H}
    = \\ \bip(VL(f_{1})[f\tensor\spxi_{i}]|V[g\tensor\spxi_{j}]) =\\
    \bip(L'(f_{1})V[f\tensor\spxi_{1}]|V[g\tensor\spxi_{j}]).
  \end{multline*}
  But the left-hand side is
  \begin{align*}
    \bip(L(f_{1}&*f)\spxi_{i}|L(g)\spxi_{j}) =
    L_{\spxi_{i},\spxi_{j}}(g^{*}*f_{1}*f) \\
    &= \int_{\go}\int_{G}\int_{G}
    \eij_{y^{-1}x}\bigl(g(y)^{*}f_{1}*f(x)\bigr)
    \rho_{ij}(y^{-1}x)\Delta(yx)^{-\half}
    \,d\lambda^{u}(x)\,d\lambda^{u}(y)\,d\mu(u) \\
    &= \int_{\go}\int_{G}\int_{G}\int_{G}
    \eij_{y^{-1}x}\bigl(g(y)^{*}f_{1}(z)f(z^{-1}x)\bigr)
    \rho_{ij}(y^{-1}x) \Delta(yx)^{-\half} \\
    &\hskip2in
    \,d\lambda^{u}(z) \,d\lambda^{u}(x)\,d\lambda^{u}(y)\,d\mu(u) \\
    &= \int_{\go}\int_{G}\int_{G}\int_{G} \eij_{y^{-1}x}\bigl(g(y)^{*}
    \hatpi\bigl(f_{1}(z)\bigr)(f)(x)\bigr)
    \rho_{ij}(y^{-1}x) \Delta(yx)^{-\half} \Delta(z)^{-\half}\\
    &\hskip2in
    \,d\lambda^{u}(z) \,d\lambda^{u}(x)\,d\lambda^{u}(y)\,d\mu(u) \\
    &= \int_{F}\int_{G} \bipu<\hatpi(f_{1}(z))f\tensor\spxi_{i},
    g\tensor\spxi_{j}> \Delta(z)^{-\half} \,d\lambda^{u}(z)\,d\mu(u) \\
    &= \int_{F}\int_{G} \bipu<\pi\bigl(f_{1}(z)\bigr)
    (f\tensor_{s(z)}\spxi_{i}), (g\tensor_{u} \spxi_{j}>
    \Delta(z)^{-\half}
    \,d\lambda^{u}(z)\,d\mu(u) \\
    &= \int_{G}\bipb
    r(z)<\pi\bigl(f_{1}(z)\bigr)V[f\tensor\spxi_{i}]\bigl(s(z)\bigr) ,
    V[g\tensor\spxi_{j}]\bigl(r(z)\bigr)> \Delta(z)^{-\half}
    \,d\nu(z).
  \end{align*}
  Thus $L'$ is the integrated form as claimed.
\end{proof}

\section{Equivalence of Fell Bundles}
\label{sec:equiv-fell-bundl}

We want to formalize the notion of the equivalence of two Fell
bundles.  As with the definition of a Fell bundle presented in
Definition~\ref{def-fell-bund}, our formulation will be a
modification of the existing definitions (cf.,
\cite{yam:xx87}*{Definition~1.5} and \cite{muh:cm01}*{Definition~10}).
We are aiming to give a
version which will be readily ``check-able'' even at the expense of
length or elegance.  First, however, we need
to establish some notation and terminology.

% Now we want to recall the notion of an equivalence between two Fell
% bundles.  What we present here is a modification of Muhly's Definition
% in \cite{muh:cm01}*{Definition~10} and Yamagami's
% \cite{yam:xx87}*{Definition~1.5}.  As with our definition of a Fell
% bundle from Definition~\ref{def-fell-bund}, we are aiming to give a
% version which will be readily ``check-able'' even at the expense of
% length or elegance.  First, however, we need
% to establish some notation and terminology.

Suppose that $p:\B\to G$ is a Fell bundle and that $q:\E\to T$ is a
\usc-Banach bundle over a left $G$-space $T$.  Then we say that $\B$
acts on (the left) of $\E$ if there is a continuous map $(b,e)\mapsto
b\cdot e$ from\index{Fell bundle!action of}
\begin{equation*}
  \B*\E:=\set{(b,e)\in\B\times\E:s\bigl(p(b)\bigr)=r\bigl(q(e)\bigr)}
\end{equation*}
to $\E$ such that
\begin{enumerate}
\item $q(b\cdot e)=p(b)\cdot q(e)$,
\item $a\cdot (b\cdot e)=(ab)\cdot e$ for appropriate $a,b\in \B$ and
  $e\in\E$, and
\item $\|b\cdot e\|=\|b\|\|e\|$.
\end{enumerate}
Of course, there is an analogous notion of a right action of a Fell
bundle.

If $T$ is a $(G,H)$-equivalence, it is important to keep in mind that
$[x\cdot t,t]\mapsto x$ is an isomorphism of the orbit space
$(T*_{s}T)/H$ onto $G$ and that $[t,t\cdot h]\mapsto h$ is an
isomorphism of $G\backslash (T*_{r}T)$ onto $H$ (see
\cite{mrw:jot87}*{\S2}).  We will write
$[t,s]_{G}$\index{tsg@$[t,s]_{G}$} for the image of the orbit of
$(t,s)$ in $G$.  Similarly, we will write
$[t,s]_{H}$\index{tsh@$[t,s]_{H}$} for the image of 
the corresponding orbit in $H$.

If $q:\E\to T$ is a Banach bundle over a $(G,H)$-equivalence, then we
will write $\E*_{s}\E$ for $\set{(e,f)\in\E\times\E:s\bigl(q(e)\bigr)
  = s\bigl(q(f)\bigr)}$, and similarly for $\E*_{r}\E$.

\begin{definition}
  \label{def-equivalence}
  Suppose that $T$ is a $(G,H)$-equivalence,\index{equivalence of Fell
    bundles} that $p_{G}:\B\to G$ and $p_{H}:\CC\to H$ are Fell
  bundles\index{Fell bundle!equivalence of}, and that $q:\E\to T$ is a
  Banach bundle.  We say that $q:\E\to T$ is a $\B\sme\CC$-equivalence
  if the following conditions hold.
  \begin{enumerate}
  \item There is a left $\B$-action and a right $\CC$-action on $\E$
    such that $b\cdot (e\cdot c)=(b\cdot e)\cdot c$ for all $b\in\B$,
    $e\in\E$ and $c\in\CC$.
  \item There are sesquilinear maps $(e,f)\mapsto \lip\B<e,f>$ from
    $\E*_{s}\E$ to $\B$ and $(e,f)\mapsto \rip\CC<e,f>$ from
    $\E*_{r}\E$ to $\CC$ such that
    \begin{enumerate}
    \item
      $p_{G}\bigl(\lip\B<e,f>\bigr)=[q(e),q(f)]_{G}$\quad\text{and}\quad
      $p_{H}\bigl(\rip\CC<e,f>\bigr)=[q(e),q(f)]_{H}$,
    \item $\lip\B<e,f>^{*}=\lip\B<f,e>$\quad\text{and}\quad
      $\rip\CC<e,f>^{*}=\rip\CC<f,e>$,
    \item $\lip\B<b\cdot e,f>=b\lip\B<e,f>$\quad\text{and}\quad
      $\rip\CC<e,f\cdot c>=\rip\CC<e,f>c$,
    \item $\lip\B<e,f>\cdot g=e\cdot\rip\CC<f,g>$.
    \end{enumerate}
  \item With the actions coming from \partref1 and the inner products
    coming from \partref2, each $E(t)$ is a $B\bigl(r(t)\bigr)\sme
    B\bigl(s(t)\bigr)$-\ib.
  \end{enumerate}
\end{definition}

\begin{lemma}
  The map $(b,e)\mapsto b\cdot e$ induces an \ib{} isomorphism of
  $B\bigl(p(b)\bigr) \tensor_{B(s(p(b)))} E(t)$ onto $E(p(b)\cdot t)$.
\end{lemma}
\begin{proof}
  The proof is similar to that for Lemma~\ref{lem-extra}.
\end{proof}

Our next observation is straightforward, but it will be helpful to
keep it in mind in the sequel.
\begin{example}
  \label{ex-self-equi}
  Suppose that $p:\B\to G$ is a Fell bundle over $G$.  Since $G$ is
  naturally a $(G,G)$-equivalence, we see immediately that $\B$ acts
  on the right and the left of itself with $\B*\B=\B^{(2)}$.  If we
  define
  \begin{equation*}
    \lip\B<a,b>:=ab^{*}\quad\text{and}\quad \rip\B<a,b>:=a^{*}b,
  \end{equation*}
then it is a simple matter to check that $p:\B\to G$ is a
$\B\sme\B$-equivalence. 
\end{example}

We can now state our main result.

\begin{thm}[The Equivalence Theorem for Fell Bundles]
  \label{thm-yam-equi}
  Suppose that $G$ and $H$ are second countable groupoids with Haar
  systems $\set{\hg^{u}}_{u\in\go}$ and $\set{\hh^{v}}_{v\in\ho}$,
  respectively. Suppose also that $p_{\B}:\B\to G$ and $p_{H}:\CC\to
  H$ are Fell bundles and that $q:\E\to T$ is a
  $\B\sme\CC$-equivalence.  Then $\X_{0}:=\sa_{c}(T;\E)$ becomes a
  $\fcs(G,\B)\sme\fcs(H,\CC)$-pre-\ib{} with respect to the actions
  and inner products given by
  \begin{align}
    \label{eq:68}
    f\cdot \xi(t)&:= \int_{G}f(x)\xi(x^{-1}\cdot t)\,d\hg^{r(t)} (x)
    \\
\xi\cdot g(t)&:= \int_{H}\xi(t\cdot h)g(h^{-1})\,d\hh^{s(t)}(h)\label{eq:69} \\
\tlips<\xi,\eta>(x)&:=\int_{H}
\lip\B<\xi(xth),\eta(th)>\,d\hh^{s(t)}(h)\label{eq:70} \\
\trips<\xi,\eta>(h)&= \int_{G} \rip\CC<\xi(x^{-1}t),\eta(x^{-1}th)> \,
d\hg^{r(t)} (x).\label{eq:71}
  \end{align}
Consequently, $\fcs(G,\B)$ and $\fcs(H,\CC)$ are Morita equivalent.
\end{thm}

\begin{remark}
  Since $T$ is a $(G,H)$-equivalence, $r_{T}(t)=r_{T}(s)$ implies that
  $s=t\cdot h$ for some $h\in H$.  Thus, we are free to choose any $t$
  in \eqref{eq:70} satisfying $r_{T}(t)=s_{G}(x)$.  Similarly, in
  \eqref{eq:71}, any $t$ satisfying $s_{T}(t)=r_{H}(h)$ will do.
\end{remark}

The proof of Theorem~\ref{thm-yam-equi} is a bit involved --- for
example, it is not even obvious that
\eqref{eq:68}--\eqref{eq:71} define continuous sections of the
appropriate bundles.  In any event, we require some preliminary
comments and set-up before launching into the proof in the next
section.  However, the next example is fundamental and should help to
motivate some of what follows.

\begin{example}
  \label{ex-self-dual-equi}
  If $p:\B\to G$ is a Fell bundle, then we can also view it as a
  $\B\sme\B$-equivalence as described in Example~\ref{ex-self-equi}.
  Then the pre-\ib{} structure imposed on $\sa_{c}(G;\B)$ by
  Theorem~\ref{thm-yam-equi} is the expected one.  The left and right
  actions are 
  given by convolution as are the inner products:
  \begin{equation*}
    \tlips<f,g>=f*g^{*}\quad\text{and}\quad\trips<f,b>=f^{*}*g.
  \end{equation*}
\end{example}

Another example that will be useful in the sequel is the following
which illustrates the symmetry inherent in the definition of equivalence.
Sadly, the notation obscures what is actually quite straightforward.

\begin{example}
  \label{ex-symmetry}
  If $T$ is a $(G,H)$-equivalence, then we can make the same space
  into a $(H,G)$-equivalence $\op{T}$ as follows.  Let
  $\iota:T\to\op{T}$ be the identity map and define a left $H$-action
  and a left $G$-action by
  \begin{align*}
    r\bigl(\iota(t)\bigr)&:=s(t)&s\bigl(\iota(t)\bigr)&:=r(t) \\
h\cdot \iota(t)&:=\iota(t\cdot h^{-1}) & \iota(t)\cdot x&:=
\iota(x^{-1}c\dot t).
  \end{align*}
Similarly, if $p_{\B}:\B\to G$ and $p_{\CC}:\CC\to H$ are Fell bundles
and if $q:\E\to T$ is a $\B\sme\CC$-equivalence, then we can form a
$\CC\sme\B$-equivalence $\bq:\bE\to\op{T}$.  Let
$\bE$ be the conjugate vector space to $\E$.  Thus if $\dmap:\E\to\bE$
is the identity map, then scalar multiplication is given by $z\cdot
\dmap(e) := \dmap(\bar z\cdot e)$ for all $z\in\C$.  Let
$\bq\bigl(\dmap(e)\bigr) := \iota(q(e))$.  We then define
$\CC$- and $\B$-actions by
\begin{equation*}
  c\cdot\dmap(e):= \dmap(e\cdot c^{*})\quad\text{and}\quad
  \dmap(e)\cdot b:=\dmap(b^{*}\cdot e),
\end{equation*}
and sesquilinear forms by
\begin{equation*}
  \blip\CC<\dmap(e),\dmap(f)>:=\rip\CC<e,f>\quad\text\quad
  \brip\B<\dmap(e),\dmap(f)>:= \lip\B<e,f>.
\end{equation*}
Now it is simply a matter of verifying the axioms to see that
$\bq:\bE\to\op{T}$ is a $\CC\sme\B$-equivalence such that
$\bE(t)={E(t)^{\sim} }$, where ${E(t)^{\sim}}$ is the dual \ib{}
to $E(t)$ (see \cite{rw:morita}*{pp.~49--50}).

Of course, we obtain actions and inner products on
$\sa_{c}(\op{T};\bE)$ using \eqref{eq:68}--\eqref{eq:71}.  When using
these side-by-side with those on $\sa_{c}(T;\E)$ it will usually be
clear from context which we are applying.  Nevertheless, to make
matters a bit easier to decode, the inner products on
$\sa_{c}(\op{T};\bE)$ will be denoted by $\tlipsd<\cdot,\cdot>$ and
$\tripsd<\cdot,\cdot>$. 
\end{example}

The following technical lemma will be useful in exploiting the
symmetry in Theorem~\ref{thm-yam-equi}.
\begin{lemma}
  \label{lem-symmetry}
  Define $\Phi:\sa_{c}(T;\E)\to\sa_{c}(\op{T};\bE)$ by
  $\Phi(\xi)\bigl(\iota(t)\bigr) := \dmap\bigl(\xi(t)\bigr)$.  Then we have
  \begin{align*}
    \Phi(\xi\cdot g)&=g^{*}\cdot \Phi(\xi)&\Phi(f\cdot
    \xi)&=\Phi(\xi)\cdot f^{*} \\
\tlipsd<\Phi(\xi),\Phi(\eta)>&:=\trips<\xi,\eta>&
\tripsd<\Phi(\xi),\Phi(\eta)>&= \tlips<\xi,\eta>.
  \end{align*}
\end{lemma}
\begin{proof}
  The lemma follows from routine computations.  For example,
  \begin{align*}
    g^{*}\cdot \Phi(\xi)\bigl(\iota(t)\bigr) &= \int_{H} g(h^{-1})^{*}
    \Phi(\xi) \bigl(h^{-1}\cdot\iota(t)\bigr) \,d\hh^{s(t)}(h) \\
&= \int_{H} \dmap\bigl(\xi(t\cdot h)g(h^{-1})\bigr)  \,d\hh^{s(t)}(h)
\\
&= \Phi(\xi\cdot g)\bigl(\iota(t)\bigr).
  \end{align*}
While
\begin{align*}
  \tlipsd<\Phi(\xi),\Phi(\eta)>(h)&= \int_{G} \blip \CC <\Phi(\xi)
  \bigl(h\cdot \iota(t')\cdot x\bigr) ,
  {\Phi(\eta)\bigl(\iota(t')\cdot x\bigr)}> \, d\hg^{r(t')}(x) \\
&= \int_{G} \blip\CC<\Phi\bigl(\iota(x^{-1}\cdot t'\cdot h^{-1})\bigr)
, {\Phi(\eta)\bigl(\iota\bigl(x^{-1}\cdot t'\bigr)\bigr)}> \,
d\hg^{r(t')}(x) \\
\intertext{which, after letting $t=t'\cdot h^{-1}$, is}
&= \int_{G} \blip\CC<\dmap\bigl(\xi(x^{-1}\cdot t)\bigr),
{\dmap\bigl(\eta(x^{-1} \cdot t\cdot h)\bigr)}> \,d\hg^{r(t)}(x) \\
&= \int_{G}\brip\CC<\xi(x^{-1}\cdot t), {\eta(x^{-1}\cdot t\cdot h)}>
\,d \hg^{r(t)}(x) \\
&= \trips<\xi,\eta>(h).
\end{align*}
The other formulas are established similarly.
\end{proof}

\begin{remark}
  \label{rem-sym-pos}
  Lemma~\ref{lem-symmetry} will be very useful in the proof of
  Theorem~\ref{thm-yam-equi}.  For example, once we establish that
  $\trips<\cdot,\cdot>$ is positive, it follows immediately that
  $\tripsd<\cdot,\cdot>$ is also positive.  Then, since
  \begin{equation*}
    \tlips<\xi,\xi>=\tripsd<\Phi(\xi),\Phi(\xi)>,
  \end{equation*}
we can say that the positivity of $\tlips<\cdot,\cdot>$ ``follows by
symmetry''. 
\end{remark}

As is now standard, the key result needed for the proof of
Theorem~\ref{thm-yam-equi} is that we have approximate identities for
$\sa_{c}(G;\B)$ in the inductive limit topology of a special form.  In
the case of Fell bundles, even the existence of a one-sided
approximate identity for $\sa_{c}(G;\B)$ of any form is not so
obvious.\footnote{The ``standard'' technique of using functions with
  small support and unit integral fail as integrals of form
  \begin{equation*}
    \int_{G}f(x)\,d\hg^{u}(x)
  \end{equation*}
are meaningless for $f\in\gcb$.
}  Nevertheless, the result we want is the following.
\begin{prop}
  \label{prop-ai}\index{gammacb@$\gcb$!approximate identity}
  Suppose that $q:\E\to T$ is a $\B\sme\CC$-equivalence.  Then there
  is a net $\set{e_{\lambda}}$ in $\gcb$ consisting of elements of the
  form
  \begin{equation*}
    e_{\lambda}=\sum_{i=1}^{n_{\lambda}}\;\tlips<\xi_{i}^{\lambda},\xi_{i}^{\lambda}>, 
  \end{equation*}
with each $\xi_{i}^{\lambda}$ in $\sa_{c}(T;\E)$, which is an
approximate identity in the inductive limit topology for the left
action of $\gcb$ on itself and on $\sa_{c}(T;\E)$.
\end{prop}

Since the proof of Proposition~\ref{prop-ai} is rather technical, we
will postpone it to Section~\ref{sec:proof-prop-refpr}.  However, if
$p:\B\to G$ is a Fell bundle, then, since $\B$ is naturally a $\B\sme
\B$-equivalence and since each $e_{\lambda}$ in
Proposition~\ref{prop-ai} is self-adjoint by construction, we obtain
as an immediate corollary that $\gcb$ itself has a two-sided
approximate identity in the inductive limit topology. (As promised,
this proves Proposition~\ref{prop-approx-id}.)

\section{Proof of the Main Theorem}
\label{sec:proof-theorem}

It is high time to see that \eqref{eq:68}--\eqref{eq:71} define
continuous sections as claimed.
To see that \eqref{eq:68} defines an element of $\sa_{c}(T;\E)$, we
note that our assumptions on $\E$ imply that $(x,t)\mapsto
f(x)\xi(x^{-1}\cdot t)$ defines a section in $\sa_{c}(G*T;\tau^{*}\E)$
where $\tau:G*T\to T$ is the projection map.  Therefore, we will get
what we need from the following lemma:
\begin{lemma}
  \label{lem-tech1}
  If $f\in\sa_{c}(G*T;\tau^{*}\E)$ and if
  \begin{equation*}
    \theta(f)(t):=\int_{G} f(x,t)\,d\hg^{r(t)}(x),
  \end{equation*}
then $\theta(f)\in\sa_{c}(T;\E)$.
\end{lemma}
\begin{proof}
  Suppose that $\supp f\subset K_{G}\times K_{T}$ with each factor
  compact.  Then $\supp\theta(f)\subset K_{T}$ and
  \begin{equation*}
    \|\theta(f)\|_{\infty}\le M \|f\|_{\infty}
  \end{equation*}
where $M$ is an upper bound for $\hg^{u}(K_{G})$ (for $u\in\go$).  It
follows that the collection of $f$ for which
$\theta(f)\in\sa_{c}(T;\E)$ is closed in the inductive limit
topology.  Since sections of the form $(x,t)\mapsto
\phi(x)\psi(t)e(t)$ for $\phi\in C_{c}(G)$, $\psi\in C_{c}(T)$ and
$e\in \sa_{c}(T;\E)$ span a dense subspace of
$\sa_{c}(G*T;\tau^{*}\E)$ by Lemma~\ref{lem-dense-sections}, it
suffices to consider $f$ of the form $f(x,t)=\phi(x)\psi(t)e(t)$.  But
then $\theta(f)=\psi(t)\hg(\phi)\bigl(r(t)\bigr)e(t)$, which defines a
continuous section on $T$ since $\set{\hg^{u}}_{u\in\go}$ is a Haar
system. 
\end{proof}

Of course, the argument for \eqref{eq:69} is similar as are the
arguments for \eqref{eq:70} and \eqref{eq:71}.  For example, to
establish \eqref{eq:70}, let $\sigma:X\stars X\to G$ be given by
$\sigma(t,t'):=[t,t']_{G}$.  Then we'll need the following analogue of
Lemma~\ref{lem-tech1}.
\begin{lemma}
  \label{lem-tech2}
  Suppose that $F\in\sa_{c}(X\stars X;\sigma^{*}\B)$.  Then
  \begin{equation*}
    \lambda(F)\bigl([t,t']_{G}\bigr):=\int_{H}F(t\cdot h,t'\cdot
    h)\,d\lambda_{H}^{s(t')}(h) 
  \end{equation*}
defines a section $\lambda(F)\in\gcb$.
\end{lemma}
\begin{proof}[Sketch of the Proof]
  Using Lemma~\ref{lem-dense-sections}, we see that sections of the
  form $F(t,t')=\phi(t,t')b\bigl([t,t']_{G}\bigr)$, with $\phi\in
  C_{c}(X\stars X)$ and $b\in\gcb$ are dense in $\sa_{c}(X\stars
  X;\sigma^{*}\B)$ in the inductive limit topology.  But if $F$ has
  this form, then $\lambda(F)=\lambda(\phi)\cdot b$ where
  \begin{equation*}
    \lambda(\phi)\bigl([t,t']_{G}\bigr):=\int_{H}\phi(t\cdot h,t'\cdot h)\,d\lambda_{H}^{s(t')}(h).
  \end{equation*}
Since $\lambda(\phi)\in C_{c}(X\stars X)$ by
\cite{mrw:jot87}*{Lemma~2.9(b)}, $\lambda(F)\in\gcb$.  We we proceed
as in Lemma~\ref{lem-tech1}.
\end{proof}

We apply this to \eqref{eq:70} as follows.  Notice that
$F(t,t'):=\blip\B<\xi(t), \eta(t')>$ defines a section in
$\sa_{c}(X\stars X;\sigma^{*}\B)$, and then $\tlips<\xi,\eta>(x) =
\lambda(F)\bigl([xt,t]_{G}\bigr)$.  Of course, \eqref{eq:71} is dealt
with in a similar fashion.  

To complete the proof, we are going to appeal to
\cite{rw:morita}*{Proposition~3.12}.  Thus we must to do the following.
\begin{enumerate}
  \item [IB1:] Show that $\sa_{c}(T;\E)$ is a both a left $\gcb$- and a
    right $\sa_{c}(H;\CC)$-pre-inner product module (as in
    \cite{rw:morita}*{Lemma~2.16}).
  \item [IB2:] Show that the inner products span dense ideals.
  \item [IB3:]
Show that the actions are bounded in that
\begin{equation*}
  \trips<f\cdot
  \xi,f\cdot\xi>\le\|f\|^{2}_{\fcs(G,\B)}\trips<\xi,\xi>\quad\text{and}
  \quad 
\tlips<\xi\cdot g,\xi\cdot g>\le \|g\|^{2}_{\fcs(H,\CC)}\tlips<\xi,\xi>.
\end{equation*}
\item [IB4:] Show that
  $\tlips<\xi,\eta>\cdot\zeta=\xi\cdot\trips<\eta,\zeta>$. 
\end{enumerate}

Verifying IB4 is a nasty little computation: recall that
\begin{equation*}
  \tlips<\xi,\eta>(x)=\int_{H}\lip\B<\xi(x\cdot t'\cdot
  h),\eta(t'\cdot h)>\, d \hh^{s(t')},
\end{equation*}
where we are free to choose any $t'$ such that $r_{T}(t')=r_{G}(x)$.
Thus if $r(x)=r(t)$, then we can let $t'=x^{-1}\cdot t$.  Therefore
\begin{align*}
  \tlips<\xi,\eta>\cdot\zeta(t)&=\int_{G}\tlip<\xi,\eta>x)\cdot
  \zeta(x^{-1}\cdot t)\,d\hg^{r(t)}(x) \\
&= \int_{G}\Bigl(\int_{H} \;\lip\B<\xi(t\cdot h),{\eta(x^{-1}\cdot
  t\cdot h)}>\,d\hh^{s(t)}(h)\Bigr)\cdot \zeta(x^{-1}\cdot t)\,d\hg^{r(t)}(x).
\end{align*}
But if $e\in E(x^{-1}\cdot t)$, then $b\mapsto b\cdot e$ is a bounded
linear map from $B(x)$ into $E(t)$.  Thus the properties of
vector-valued integrals (cf. \cite{wil:crossed}*{Lemma~1.91}) imply that
\begin{align*}
  \tlips<\xi,\eta>\cdot \zeta(t)&= \int_{G}\int_{H} \;\lip\B<\xi(t\cdot
  h), {\eta(x^{-1}\cdot t\cdot h)}>\cdot\zeta(x^{-1}\cdot t)
  \,d\hh^{s(t)}(h) \,d\hg^{r(t)}(x) \\
\intertext{which, since $E(t)$ is an \ib, is}
&=\int_{G}\int_{H} \xi(t\cdot h)\cdot \rip\CC<\eta(x^{-1}\cdot t\cdot
h), { \zeta(x^{-1}\cdot t)}>\,d\hh^{s(t)}(h)\,d\hg^{r(t)}(x) \\
\intertext{which, after using Fubini and the properties of
  vector-valued integrals, is}
&=\int_{H}\xi(t\cdot h) \cdot \Bigl( \int_{G} \rip\CC<\eta(x^{-1}\cdot
t\cdot h),{\zeta(x^{-1}\cdot
  t)}>\,d\hg^{r(t)}(x)\Bigr)\,d\hh^{s(t)}(h) \\
\intertext{which, after replacing $t\cdot h$ by $t'$ in the inner integral, is}
& =\int_{H}\xi(t\cdot h) \cdot \Bigl( \int_{G} \rip\CC<\eta(x^{-1}\cdot
t'),{\zeta(x^{-1}\cdot
  t'\cdot h^{-1})}>\,d\hg^{r(t')}(x)\Bigr)\,d\hh^{s(t)}(h) \\
&= \int_{H}\xi(t\cdot h)\cdot \trips<\eta,\zeta>(h^{-1})\, d\hh^{s(t)}
(h) \\
&= \xi\cdot\trips<\eta,\zeta>(t).
\end{align*}
This proves IB4.

To verify IB1, it suffices, by symmetry, to consider only
$\tlips<\cdot,\cdot>$.  The algebraic properties are routine.  For
example, the axioms of an equivalence guarantee that $f\mapsto
\lip\B<f, e>$ is a bounded linear map of $E(x\cdot t) $ into
$B(x)$.  Thus the usual properties of vector-valued integration (cf. 
\cite{wil:crossed}*{Lemma~1.91}) imply that
\begin{equation*}
  \tlips<f\cdot \xi(t),\eta(t)>=\int_{G}\lip\B<f(x)\cdot
  \xi(x^{-1}\cdot t),{\eta(t)}>\,d\hg^{r(t)}(x).
\end{equation*}
Thus we can compute as follows:
\begin{align*}
  \tlips<f\cdot \xi,\eta>(x)&= 
\int_{H}\lip\B<f\cdot \xi(x\cdot t\cdot h), {\eta(t\cdot h)}>\,d\hh^{s(t)}(h) \\
&=\int_{H}\int_{G} \lip\B<f(y)\cdot \xi(y^{-1}x\cdot t\cdot h),
{\eta(t\cdot h)}>\,d\hg^{r(x)}(y)\,d\hh^{s(t)}(h) \\
\intertext{which, in view of part~\partref2(iii) of
  Definition~\ref{def-equivalence}, is}
&= \int_{H}\int_{G} f(y)\lip\B<\xi(y^{-1}x\cdot t\cdot h),{\eta(t\cdot
h)} > \,d\hg^{r(x)}(y)\,d\hh^{s(t)}(h) \\
&= \int_{G} f(y)\cdot \tlips<\xi,\eta>(y^{-1}x)\,d\hg^{r(x)}(y) \\
&= f*\tlips<\xi,\eta>(x).
\end{align*}
To complete the verification of IB1, we only need to see that the
pre-inner products are positive.  We will show this (as well as IB2)
using the approximate identity developed in
Proposition~\ref{prop-ai}.  

It is not hard to see that the inner products respect the inductive
limit topology, and convergence in the inductive limit topology in
$\gcb$ certainly implies convergence in the \cs-norm.  Thus we have
\begin{align*}
  \trips<\xi,\xi>&=\lim_{\lambda} \trips<e_{\lambda}*\xi,\xi> \\
&= \lim_{\lambda}\sum_{i=1}^{n_{\lambda}}
\trips<{\tlips<\xi_{i}^{\lambda},\xi_{i}^{\lambda}>}\cdot\xi,\xi> \\
&= \lim_{\lambda}\sum_{i} \trips<\xi_{i}^{\lambda}\cdot
{\trips<\xi_{i}^{\lambda},\xi>}, \xi> \\
&= \lim_{\lambda}\sum_{i} {\trips<\xi_{i}^{\lambda},\xi>^{*} }
\trips<\xi_{i}^{\lambda},\xi>,
\end{align*}
which shows positivity of $\trips<\cdot,\cdot>$.  The positivity of
$\tlips<\cdot,\cdot>$ follows by symmetry.

Since $\xi=\lim_{\lambda}e_{\lambda}*\xi$, similar considerations show
the span of the inner products are dense as required in IB2.

This leaves only IB3 to establish, and by symmetry, it is enough to
show that
\begin{equation}
  \label{eq:83}
  \trips<f\cdot \xi,f\cdot\xi>\le \|f\|^{2}\trips<\xi,\xi>.
\end{equation}
If $\rho$ is a state on $\fcs(G,\B)$, then the pre-inner product
\begin{equation*}
  \ip(\cdot | \cdot)_{\rho}:=\rho\bigl(\trips<\cdot,\cdot>\bigr)
\end{equation*}
makes $\sa_{c}(T;\E)$ into a pre-Hilbert space.  Let $\H_{0}$ be the
dense image of $\sa_{c}(T;\E)$ in the Hilbert space completion.  It is
not hard to check that the
left action of $\gcb$ on $\sa_{c}(T;\E)$ defines a pre-representation $L$
of $\B$ on $\H_{0}$.  Now we employ the full power of the 
Disintegration Theorem (Theorem~\ref{thm-fell-disintegration})
to conclude that $L$ extends to a bona fide representation of $\gcb$
(which is therefore norm-reducing for the universal norm).  Therefore
\begin{equation*}
  \rho\bigl(\trips<f\cdot \xi,f\cdot \xi>\bigr) \le
  \|f\|^{2}_{\fcs(G,\B)} \rho\bigl(\trips<\xi,\xi>\bigr).
\end{equation*}
Since this holds for all $\rho$, we have established \eqref{eq:83}.
This completes the proof of Theorem~\ref{thm-yam-equi} with the
exception of the proof of the existence of approximate identities of
the required type.

\section{Proof of Proposition~\ref{prop-ai}}
\label{sec:proof-prop-refpr}
Recall that we write $A$ for the \emph{\cs-algebra} $\sa_{0}(\go;\B)$.
Abusing notation a bit, we'll let $\sa_{c}(\go;\B)^{+}$ denote the
positive elements in $A$ with compact support.

\begin{lemma}
  \label{lem-A-approx-id}
  Let $\Lambdasubc:=\set{a\in \sa_{c}(\go;\B)^{+}:\|a\|\le 1}$.  Then
  $\Lambdasubc$ is a net directed by itself such that for all
  $\xi\in\sa_{c}(T;\E)$, we have $a\cdot \xi\to \xi$ uniformly, where
  $a\cdot \xi(t):=a\bigl(r(t)\bigr) \cdot \xi(t)$.
\end{lemma}
\begin{proof}
  Note that if $a\ge b\ge 0$ in $A^{+}$, then for each $u\in\go$, we
  have $a(u)\ge b(u)\ge0$.  Also, given $b\in A(u)^{+}$, there is an
  $a\in \Lambdasubc$ such that $a(u)=b$.  It follows that, for any $t\in T$,
  $a\bigl(r(t)\bigr) \cdot \xi(t)\to \xi(t)$ as $a\nearrow 1$ (since
  $\xi(t)\in E(t)$, and $E(t)$ is a left Hilbert
  $A\bigl(r(t)\bigr)$-module).

  If $a\cdot g$ does not converge to $g$ uniformly, then there is an
  $\epsilon>0$, a subnet $\set{a_{i}}$ and $t_{i}\in \supp \xi$ such
  that
  \begin{equation}
    \label{eq:66}
    \|a\bigl(r(t_{i})\bigr)\cdot \xi(t_{i})-\xi(t_{i})\|\ge\epsilon.
  \end{equation}
  Since $\supp \xi$ is compact, we can pass to a subnet, relabel, and
  assume that $t_{i}\to t$.  Then
  \begin{multline*}
    \|a\bigl(r(t_{i})\bigr)\cdot \xi(t_{i})-\xi(t_{i}) \| \le
    \|a\bigl(r(t_{i})\bigr)\cdot \xi(t_{i}) - a\bigl(r(t)\bigr)\cdot \xi(t)\| \\
    + \|a\bigl(r(t)\bigr) \cdot \xi(t)-\xi(t)\| + \|\xi(t)-\xi(t_{i})\|.
  \end{multline*}
  Since the terms on the right-hand side all tend to zero with $i$, we
  eventually contradict \eqref{eq:66}.  This completes the proof.
\end{proof}

\begin{lemma}
  \label{lem-approx-a}
  Suppose that $a\in \Lambdasubc:=\sa_{c}(\go;\B)^{+}$, that
  $\epsilon>0$ and that $K\subset T$ is compact.  Then there are
  $\xi_{1},\dots,\xi_{n}\in\sa_{c}(T;\E)$ such that
  \begin{equation}
    \label{eq:72}
    \|a\bigl(r(t)\bigr) - \sum_{i=1}^{n}\;
    \blip\B<\xi_{i}(t),\xi_{i}(t)>\|<\epsilon\quad\text{for all $t\in K$.}
  \end{equation}
\end{lemma}
\begin{proof}
  Since $\E$ is an equivalence, $E(t)$ is a full left Hilbert
  $B\bigl(r(t)\bigr)$-module.  It follows from
  \cite{muhwil:nyjm08}*{Lemma~6.3} that for each $t\in T$, there are
  $\zeta_{1}^{t},\dots,\zeta_{n_{t}}^{t}\in \sa_{c}(T;\E)$ such that
  \eqref{eq:72} holds at $t$.  Since $K$ is compact and since
  $b\mapsto \|b\|$ is upper semicontinuous, there is a finite cover
  $\set{V_{1},\dots,V_{m}}$ of $K$ and sections
  $\zeta_{1}^{j},\dots,\zeta_{n_{j}}^{j}$ such that
  \begin{equation*}
    \|a\bigl(r(t)\bigr)
    -\sum_{i=1}^{n_{j}}\lip\B<\zeta_{i}^{j}(t),\zeta_{i}^{j}(t)> \|
    <\epsilon\quad\text{for all $t\in V_{j}$.}
  \end{equation*}
  Let $\set{\phi_{j}}$ be a partition of unity subordinate to the
  $\set{V_{j}}$ so that each $\phi_{j}\in C_{c}^{+}(T)$ with
  $\supp\phi_{j} \subset V_{j}$, $\sum_{j}\phi_{j}(t)=1$ if $t\in K$
  and the sum is less than or equal to $1$ otherwise.  Then
  \begin{equation*}
    \|a\bigl(r(t)\bigr) -
    \sum_{j=1}^{m}\phi_{j}(t)\sum_{n=1}^{n_{j}}\lip\B
    <\zeta_{i}^{j}(t),\zeta_{i}^{j}(t) >\|<\epsilon\quad\text{for all
      $t\in K$.}
  \end{equation*}
  Now we can let $\xi_{ij}(t):=\phi_{j}(t)^{\frac12}\zeta_{i}^{j}(t)$.
  This suffices.
\end{proof}

\begin{remark}
  \label{rem-myf}
  To ease the notation a bit, we are going to let
  \begin{equation}\label{eq:74}
    \myf(x,t)=\sum_{i=1}^{n}\;\lip\B<\xi_{i}(x\cdot t),\xi_{i}(t)>
  \end{equation}
  so that $\myf\bigl(r(t),t)\bigr)$ is the sum appearing in
  \eqref{eq:72}.  Notice that $\myf(x,t)\in B(x)$.
\end{remark}

\begin{proof}[Proof of Proposition~\ref{prop-ai}]
  In view of Examples \ref{ex-self-equi}~and \ref{ex-self-dual-equi},
  it suffices to treat only the case of $\gcb$ acting on
  $\sa_{c}(T;\E)$.

  As a first step, for each $a\in\Lambdasubc$ and finite set $F\subset
  \sa_{c}(T;\E)$, we will produce a net $\set{e_{j}^{a,F}}_{j\in J}$
  of the required form such that, as $j$ increases,
  $e_{j}^{a,F}*\xi\to a\cdot \xi$ in the inductive limit topology for
  each $\xi\in F$.  These nets are to be indexed by pairs
  $j=(V,\epsilon)$ where $\epsilon>0$ and $V$ is a conditionally
  compact neighborhood of $\go$ all contained in a fixed conditionally
  compact neighborhood $V_{0}$.  Since we will arrange that $\supp
  e_{V,\epsilon}^{a,F}\subset V$, it will follow that
  \begin{equation*}
    \supp e_{V,\epsilon}^{a,F}*\xi\subset V\cdot \supp\xi\subset
    V_{0}\cdot \supp \xi.
  \end{equation*}
  Since $V_{0}$ is conditionally compact, there is a compact set
  $K_{F}$, depending only on $F$, such that
  \begin{equation*}
    \supp  e_{V,\epsilon}^{a,F}*\xi \subset V_{0}\cdot \supp \xi \subset
    V_{0}\cap r_{G}^{-1}(r_{T}\bigl(\supp \xi)\bigr)\cdot \supp \xi
    \subset K_{F}. 
  \end{equation*}
  Consequently, we just have to show that $e_{j}^{a,F}*\xi\to a\cdot
  \xi$ uniformly.

  Therefore, we fix $a\in\Lambdasubc$ and $F\subset \sa_{c}(T;\E)$.
  We also fix $(V,\epsilon)\in J$.  Let $D\subset T$ be a compact set
  such that $V_{0}\cdot \supp\xi\subset D$ for all $\xi\in F$.  Given
  $\epsilon>0$, Lemma~\ref{lem-approx-a} implies that there are
  $\xi_{1},\dots,\xi_{n}\in\sa_{c}(T;\E)$ such that
  \begin{equation}
    \label{eq:73}
    \|a\bigl(r(t)\bigr) -\myf\bigl(r(t),t\bigr)\|<\epsilon \quad\text{for
      all $t\in D$,}
  \end{equation}
  where $\myf$ is defined as in \eqref{eq:74}.
  \begin{remark}
    \label{rem-norm-myf}
    Since $\|a(u)\|\le1$ for all $u$, we can assume that
    $\|\myf\bigl(r(t),t\bigr)\|\le 2$ for all $t\in D$.
  \end{remark}

  As in \cite{muhwil:nyjm08}*{\S6}, we can find functions
  $\phi_{1},\dots,\phi_{m}\in C_{c}^{+}(T)$ such that if
  \begin{equation*}
    F(x,t):=\sum_{j=1}^{m}\phi_{j}(t)\phi_{j}(x^{-1}\cdot t),
  \end{equation*}
  then
  \begin{gather}
    F(x,t)=0\quad\text{if $x\notin V$ or $t\notin D$,}\label{eq:75}\\
    \Bigl| \int_{G}\int_{H}F(x,t\cdot h)\,d\hh^{s(t)}(h)\,d\hg^{r(t)}(x) -
    1
    \Bigr| < \epsilon\quad\text{provided $t\in D$ and} \label{eq:76}\\
    \int_{G}\int_{H}F(x,t\cdot h)\,d\hh^{s(t)}(h)\,d\hg^{r(t)}(x) \le 2
    \quad\text{for all $t$.}\label{eq:77}
  \end{gather}
  Then we define$
    e = e^{a,F}_{V,\epsilon}$
  by
  \begin{equation*}
    e(x)=\sum_{ij} \tlips<\omega_{ij},\omega_{ij}>,
  \end{equation*}
  where $\omega_{ij}(t):=\phi_{i}(t)\xi_{j}(t)$.  Thus
  \begin{equation*}
    e(x)=\int_{H} F(x,t\cdot h)\myf(x,t\cdot h)\,\hh^{s(t)}(h).
  \end{equation*}

\begin{claim}
  There is a conditionally compact neighborhood $V$ of $\go$ such that
  $y\in V$ and $(s,t)\in D\starr D$ implies that
  \begin{equation}\label{eq:79}
    \|\myf(y,t)\xi(y^{-1}\cdot s)-\myf\bigl(r(t),t\bigr)
    \xi(s)\|<\epsilon \quad\text{for all $\xi\in F$.}
  \end{equation}
\end{claim}
\begin{proof}
  [Proof of Claim] It suffice to produce $V$ such that \eqref{eq:79}
  holds for a fixed $\xi\in F$.  If the claim were false, then for
  some $\epsilon_{0}>0$ and every neighborhood $V$ of $\go$ inside
  some fixed conditionally compact neighborhood $V_{0}$, we could find
  $(s_{V},t_{V})\in D\starr D$ and $y_{V}\in V\cap
  r_{G}^{-1}\bigl(r_{X}(D)\bigr)$ such that
  \begin{equation}
    \label{eq:78}
    \|\myf(y_{V},t_{V})\xi(y_{V}^{-1}\cdot s_{V})
    -\myf\bigl(r(t_{V}),t_{V}\bigr) \xi(s_{V})\|\ge\epsilon_{0}.
  \end{equation}
  Since $V_{0}$ is conditionally compact, $V_{0}\cap
  r_{G}^{-1}\bigl(r_{X}(D)\bigr)$ has compact closure.  Therefore, we
  can pass to a subnet, relabel, and assume that $(s_{V},t_{V})\to
  (s,t)$ while $y_{V}\to r(t)=r(s)$.  Since $b\mapsto \|b\|$ is upper
  semicontinuous on $\B$, this eventually contradicts \eqref{eq:78}.
  This completes the proof of the claim.
\end{proof}
Now we compute as follows:
\begin{align}
  \|e*\xi(t)&{}-a\cdot\xi(t)\| \notag \\
  &= \Bigl\| \int_{G}\int_{H} F(y,t\cdot h)\myf(y,t\cdot h)
  \xi(y^{-1}\cdot t)\,\hh^{s(t)}(h)\,\hg^{r(t)}(y) -
  a\bigl(r(t)\bigr)\cdot
  g(t)\Bigr\| \notag\\
  &=  \int_{G}\int_{H} F(y,t\cdot h)
  \| \myf(y,t\cdot h) \xi(y^{-1} \cdot t)-
  \myf\bigl(r(t),t\cdot h\bigr) \cdot \xi(t)\|\label{eq:80}\\
  &\hskip3in d \hh^{s(t)}(h)\,d\hg^{r(t)}(y)
   \notag \\
  & \qquad + \int_{G}\int_{H} F(y,t\cdot h)\|\bigl(\myf\bigl(r(t),t\cdot
  h)-a\bigl(r(t)\bigr)\bigr)
  \cdot \xi(t)\bigr)\| \label{eq:81} \\
&\hskip3in d\hh^{s(t)}(h)\,d\hg^{r(t)}(y) \notag \\
  & \qquad \qquad + \Bigl| \int_{G}\int_{H} F(y,t\cdot h)
  \,d\hh^{s(t)}(h)\,d\hg^{r(t)}(y)- 1 \Bigr| \|a\bigl(r(t)\bigr)\cdot
  \xi(t) \|. \label{eq:82}
\end{align}
Since $\supp e*\xi\subset D$, $\|e*\xi(t)-a\cdot \xi(t)\|=0$ if
$t\notin D$.  Therefore the integrand in \eqref{eq:80} is nonzero only
if $y\in V$ and $(t,t\cdot h)\in D\starr D$.  Thus \eqref{eq:80} is
bounded by $2\epsilon$ by our choice of $V$ and by \eqref{eq:77}.
Since \eqref{eq:81} also vanishes if $t\cdot h\notin D$, \eqref{eq:81}
is bounded by $2\epsilon\|\xi\|_{\infty}$ by \eqref{eq:73} and
\eqref{eq:77}.  Equation~\eqref{eq:82} is bounded by
$\epsilon\|\xi\|_{\infty}$ by \eqref{eq:76} (since $t\in D$) and since
$\|a(u)\|\le 1$ for all $u$.  Thus $\|e*\xi(t)-a\cdot\xi(t)\|\le
2\epsilon+2\epsilon\|\xi\|_{\infty} + \epsilon\|\xi\|_{\infty}$.  It
follows that $e_{V,\epsilon}^{a,F}*\xi\to a\cdot \xi$ uniformly with
$(V,\epsilon)$ for all $\xi\in F$.  This completes the ``first step''.

For the next step, view $\set{e^{a,F}_{V,\epsilon}}$ as net indexed by
increasing $a$ and $F$ and decreasing $V$ and $\epsilon$.  Since we
have already dealt with the supports, we will complete the proof by
showing that $\set{e^{a,F}_{V,\epsilon}}$ has a subnet
$\set{d_{\lambda}}_{\lambda\in\Lambda}$ such that
$d_{\lambda}*\xi\to\xi$ uniformly for all $\xi\in\sa_{c}(T;\E)$.

Let $\Lambda$ be the collection of $5$-tuples $(a,F,V,\epsilon,n)$
where $a\in\Lambdasubc$, $F$ is a finite subset of $\sa_{c}(T;\E)$,
$V$ is a conditionally compact neighborhood of $\go$ contained in
$V_{0}$, $\epsilon>0$ and $n\in\Z^{+}$ is such that
\begin{equation*}
  \|e_{V,\epsilon}^{a,F}*\xi-\xi\|_{\infty}<\frac 1n\quad\text{for
    all $\xi\in F$.}
\end{equation*}
To see that $\Lambda$ is directed, suppose that
$(a_{i},F_{i},V_{i},\epsilon_{i},n_{i})$ is an element of $\Lambda$
for $i=1,2$.  Let $F=F_{1}\cup F_{2}$.  Using
Lemma~\ref{lem-A-approx-id}, there is a $b\in\Lambdasubc$ such that
\begin{equation*}
  \|b\cdot\xi-\xi\|_{\infty}<\frac1{2(n_{1}+n_{2})}\quad\text{for all
    $\xi\in F$.}
\end{equation*}
By the first part of this proof, we can find $(V,\epsilon)$ dominating
$(V_{i},\epsilon_{i})$ for $i=1,2$ such that
\begin{equation*}
  \| e^{b,F}_{V,\epsilon}*\xi-b\cdot\xi\|_{\infty}<
  \frac1{2(n_{1}+n_{2})}\quad \text{for all $\xi\in F$.}
\end{equation*}
Then $(b,F,V,\epsilon,n_{1}+n_{2})\in\Lambda$ and $\Lambda$ is
directed.

We then get a subnet $\set{d_{\lambda}}_{\lambda\in\Lambda}$ of
$\set{e^{a,F}_{V,\epsilon}}$ by letting
$d_{a,F,V,\epsilon,n}=e^{a,F}_{V,\epsilon}$.  To see that
$d_{\lambda}*\xi\to\xi$ in the inductive limit topology, it suffices,
since $V\subset V_{0}$ gives control of the supports, to show that
$d_{\lambda}*\xi\to\xi$ uniformly.  If $\delta>0$, then there is a
$n_{0}$ such that $\frac1{n_{0}}<\delta$, and we can find $a_{0}$ such
that implies that
\begin{equation*}
  \|a_{0}\cdot\xi-\xi\|_{\infty}< \frac1{2n_{0}}.
\end{equation*}
If $F:=\set{\xi}$, then we can find $(V_{0},\epsilon_{0})$ such that
\begin{equation*}
  \|e_{V_{0},\epsilon_{0}}^{a_{0},F_{0}}*\xi-a_{0}\cdot\xi\|_{\infty}<
  \frac1{2n_{0}}. 
\end{equation*}
Therefore $(a_{0},F_{0},V_{0},\epsilon_{0},n_{0})\in\Lambda$ and if
$(a,F,V,\epsilon,n)\ge (a_{0},F_{0},V_{0},\epsilon_{0},n_{0})$, then
we have
\begin{equation*}
  \|e^{a,F}_{V,\epsilon}*\xi-\xi\|_{\infty}<\delta.
\end{equation*}
This completes the proof.
\end{proof}

\appendix

\section{Upper Semicontinuous Banach Bundles}
\label{sec:usc-banach-bundles}

We are interested in fibred \cs-algebras as a groupoid $G$ must act on
the sections of a bundle that is fibred over the unit space (or over
some $G$-space).  In \cite{ren:jot87} and in \cite{kmrw:ajm98}, it was
assumed that the algebra $A$ was the section algebra of a \cs-bundle
as defined, for example, by Fell in \cite{fd:representations1}.
However recent work has made it clear that the notion of a \cs-bundle,
or for that matter a Banach bundle, as defined in this way is
unnecessarily restrictive, and that it is sufficient to assume only
that $A$ is a $C_{0}(\go)$-algebra
\cites{gal:94,gal:kt99,khoska:jot04,khoska:jram02}.  However, our
approach here, as in \cite{kmrw:ajm98} (and in \cite{ren:jot87}),
makes substantial use of the total space of the underlying bundle.
Although it predates the term ``\cox-algebra'', the existence of a
Banach 
bundle whose section algebra is a given Banach \cox-module goes back to
\cites{hof:74,hofkei:lnm79,hof:lnm77,dg:banach}.  We
give some of the basic definitions and properties here for the sake of
completeness.

The basics on \cox-algebras\index{cox-algebras@\cox-algebras} is available from numerous sources.  A
summary can be found in \cite{wil:crossed}*{Appendix~C.1}.  A
discussion of \usc-\cs-bundles and their connection to \cox-algebras
is laid out in \cite{wil:crossed}*{Appendix~C.2}.  Here we need a bit
more as the definition of Fell bundles requires \usc-Banach bundles
(as opposed to \usc-\cs-bundles).

This definition is a minor variation on
\cite{dg:banach}*{Definition~1.1} and should be compared with
\cite{wil:crossed}*{Definition~C.16}.
\begin{definition}
  \label{def-usc-bundle}
  A \emph{\usc-Banach bundle}\index{upper semicontinuous-Banach
    bundle@\usc-Banach bundle} over a topological space $X$ is a
  topological space $\A$ together with a continuous, open surjection
  $p=p_{\A}:\A\to X$ and complex Banach space structures on each fibre
  $\A_{x}:= p^{-1}(\set x)$ satisfying the following axioms.
  \begin{itemize}
  \item [B1:] The map $a\mapsto\|a\|$ is upper semicontinuous from
    $\A$ to $\R^{+}$.  (That is, for all $\epsilon>0$,
    $\set{a\in\A:\|a\|\ge\epsilon}$ is closed.)
  \item [B2:] If $\A*\A:=\set{(a,b)\in\A\times\A:p(a)=p(b)}$, then
    $(a,b)\mapsto a+b$ is continuous from $\A*\A$ to $\A$.
  \item [B3:] For each $\lambda\in\C$, $a\mapsto \lambda a$ is
    continuous from $\A$ to $\A$.
  \item [B4:] If $\set{a_{i}}$ is a net in $\A$ such that $p(a_{i})\to
    x$ and such that $\|a_{i}\|\to 0$, then $a_{i}\to0_{x}$ (where
    $0_{x}$ is the zero element in $\A_{x}$).
  \end{itemize}
\end{definition}

Since $\set{a\in\A:\|a\|<\epsilon}$ is open for all $\epsilon>0$, it
follows that whenever $a_{i}\to0_{x}$ in $\A$, then $\|a_{i}\|\to 0$.
Therefore the proof of
\cite{fd:representations1}*{Proposition~II.13.10} implies that
\begin{itemize}
\item [B3$'$:] The map $(\lambda,a)\to\lambda a$ is continuous from
  $\C\times\A$ to $\A$.
\end{itemize}
\begin{definition}
  \label{def-cs-bundle}
  An \usc-\cs-bundle\index{upper
    semicontinuous-\cs-bundle@\usc-\cs-bundle} is an \usc-Banach
  bundle $:\A\to X$ 
  such that each fibre is a \cs-algebra such that
  \begin{itemize}
  \item [B5:] The map $(a,b)\mapsto ab$ is continuous from $\A*\A$ to
    $\A$.
  \item [B6:] The map $a\mapsto a^{*}$ is continuous from $\A$ to
    $\A$.
  \end{itemize}
\end{definition}

If axiom B1 is replaced by
\begin{itemize}
\item [B1$'$:] The map $a\mapsto \|a\|$ is continuous,
\end{itemize}
then $p:\A\to X$ is called a Banach bundle (or a \cs-bundle).  Banach
bundles are studied in considerable detail in \S\S13--14 of Chapter~II
of \cite{fd:representations1}.

If $p:\A\to X$ is an \usc-Banach bundle, then a continuous function
$f:X\to \A$ such that $p\circ f=\id_{X}$ is called a \emph{section}.
The set of sections is denoted by $\sa_{}(X;\A)$.  We say that
$f\in\sa_{}(X;\A)$ \emph{vanishes at infinity} if the the closed set
$\set{x\in X:|f(x)|\ge\epsilon}$ is compact for all $\epsilon>0$.  The
set of sections which vanish at infinity is denoted by
$\sa_{0}(X;\A)$, and the latter is a Banach space
with respect to the supremum norm: $\|f\|=\sup_{x\in X}\|f(x)\|$ (cf.
\cite{dg:banach}*{p.~10} or \cite{wil:crossed}*{Proposition~C.23}); in
fact, $\sa_{0}(X;\A)$ is a Banach \cox-module for the natural
\cox-action on sections. (In particular, the uniform
limit of sections is a section.)  We also use $\sa_{c}(X;\A)$ for the vector
space of sections with compact support (i.e., $\set{x\in
  X:f(x)\not=0_{x}}$ has compact closure).    Moreover, if $p:\A\to X$ is an
\usc-\cs-bundle, then the set of sections is clearly a $*$-algebra
with respect to the usual pointwise operations, and $\sa_{0}(X;\A)$
becomes a \cox-algebra with the obvious $\cox$-action.  However, for
arbitrary $X$, there is no reason to expect that there are any
non-zero sections --- let alone non-zero sections vanishing at
infinity or which are compactly supported.  A \usc-Banach bundle is
said to have \emph{enough sections} if given $x\in X$ and $a\in\A_{x}$
there is a section $f$ such that $f(x)=a$.  If $X$ is a Hausdorff
locally compact space and if $p:\A\to X$ is a Banach bundle, then a
result of Douady and Soglio-H\'erault implies there are enough
sections \cite{fd:representations1}*{Appendix~C}.  Hofmann has noted
that the same is true for \usc-Banach bundles over Hausdorff locally
compact spaces \cite{hof:lnm77} (although the details remain
unpublished \cite{hof:74}).  In this article, we are assuming all our
\usc-Banach bundles have enough sections.

The following lemma is useful as it shows the topology on $\A$ is tied
to the continuous sections.
\begin{lemma}
  \label{lem-top-sections}
  Suppose that $p:\A\to X$ is an \usc-Banach-bundle.  Suppose that
  $\set{a_{i}}$ is a net in $\A$, that $a\in A$ and that
  $f\in\sa_{0}(X;\A)$ is such that $f\bigl(p(a)\bigr)=a$.  If
  $p(a_{i})\to p(a)$ and if $\|a_{i}-f\bigl(p(a_{i})\bigr) \|\to 0$,
  then $a_{i}\to a$ in $\A$.
\end{lemma}
\begin{proof}
  We have $a_{i}-f\bigl(p(a_{i})\bigr)\to 0_{p(a)}$ by axiom~B4.
  Hence
  \begin{equation*}
    a_{i}=(a_{i}-f\bigl(p(a_{i})\bigr) + f\bigl(p(a_{i})\bigr)\to
    0_{p(a)} + a =a.\qed
  \end{equation*}
  \renewcommand{\qed}{}
\end{proof}

A slightly more general result is
\cite{wil:crossed}*{Proposition~C.20}.  Results such as these can be
used to show that the section algebra $\sa_{0}(X;\A)$ is complete ---
see for example the proof of \cite{wil:crossed}*{Proposition~C.23}.

We will want to make repeated use of the following.  It has a
straightforward proof similar to that given in
\cite{wil:crossed}*{Proposition~C.24}.

\begin{lemma}
  \label{lem-dense-sections}
  Suppose that $p:\A\to X$ is an \usc-Banach bundle over a locally
  compact Hausdorff space $X$, and that $B$ is a subspace of
  $A=\sa_{0}(X;\A)$ which is closed under multiplication by functions
  in $\cox$ and such that $ \set{f(x):f\in B} $ is dense in $A(x)$ for
  all $x\in X$.  Then $B$ is dense in~$A$.
\end{lemma}

Now we come to our vector-valued Tietze Extension Theorem.  The proof
is lifted from \cite{fd:representations1}*{Theorem~II.14.8}.  However,
we have to make accommodations for the lack of continuity of the norm
function on an \usc-Banach bundle.
\begin{prop}
  \label{prop-tietze}\index{Tietze Extension Theorem}\index{upper
    semicontinuous-Banach bundle@\usc-Banach bundle!Tietze Extension
    Theorem} 
  Suppose that $p:G\to\B$ is an \usc-Banach bundle and that $Y$ is a
  closed subset of $G$.  If $g\in\sa_{c}(Y;\B\restr Y)$, then there is
  a $f\in\gcb$ such that $f(y)=g(y)$ for all $y\in Y$.
\end{prop}
\begin{proof}
  Let $C:=\supp g$ and let $U$ be a pre-compact open neighborhood of
  $C$ in $G$.  Let
  \begin{equation*}
    \AA:=\set{f\restr Y:f\in\gcb}\subset \sa_{c}(Y;\B\restr Y).
  \end{equation*}
  Since every $\psi\in C_{c}(Y)$ is the restriction of some
  $\phi\in\ccg$ (by the scalar-valued Tietze theorem
  \cite{wil:crossed}*{Lemma~1.42}), $\AA$ is dense in
  $\sa_{c}(Y;\B\restr Y)$ in the \ilt{} by
  Lemma~\ref{lem-dense-sections}.  Hence there is $\set{f_{n}}\subset
  \gcb$ such that
  \begin{equation*}
    f_{h}\restr Y\to g\quad\text{uniformly on $\overline{U}\cap Y$.}
  \end{equation*}
  By multiplying by a function which is $1$ on $C$ and vanishes off
  $U$, we can assume that each $f_{n}$ vanishes off $U$.  Passing to a
  subsequence and relabeling, we can assume that
  \begin{equation*}
    \sup\set{\|f_{n}(y)-f_{n-1}(y)\|:y\in \overline{U}\cap Y} <
    \frac1{2^{n}} \quad(n\ge2).
  \end{equation*}
  Let $h_{n}':=f_{n}-f_{n-1}$.  Since $x\mapsto \|h_{n}'(x)\|$ is
  upper semicontinuous,
  \begin{equation*}
    A_{n}:=\set{x\in G:\|h_{n}'(x)\|\ge \frac 1{2^{n}}}
  \end{equation*}
  is closed and disjoint from the closed set $\overline{U}\cap Y$.
  Therefore there is a $\phi_{n}\in\ccpg$ such that $0\le
  \phi_{n}(x)\le 1$ for all $x$, $\phi_{n}(y)=1$ if $y\in
  \overline{U}\cap Y$ and $\phi_{n}(x)=0$ if $x\in A_{n}$.  Then if we
  let $h_{n}=\phi_{n}\cdot h_{n}'$, we have arranged that
  $\|h_{n}(x)\|\le 2^{-n}$ for all $x\in G$, and that
  $h_{n}(y)=h_{n}'(y)$ if $y\in \overline{U}\cap Y$.  Since $B(x)$ is
  complete, we can define a section $f:G\to \B$ by
  \begin{equation*}
    f(x):=f_{1}(x)+\sum_{n=2}^{\infty} h_{n}(x).
  \end{equation*}
  Clearly $f$ vanishes off $U$ and as it is the uniform limit of
  elements of $\gcb$, it too is in $\gcb$.

  On one hand, if $y\in \overline{U}\cap Y$, then
  \begin{align*}
    f(y)&=f_{1}(y)+\sum_{n=2}^{\infty} h_{n}'(y) \\
    &= \lim_{n}f_{n}(y) \\
    &=g(y).
  \end{align*}
  On the other hand, if $y\notin \overline{U}\cap Y$, then both $g(y)$
  and $f(y)$ are zero.  Thus $g=f\restr Y$ as required.
\end{proof}

\section{An Example: The Scalar Case}
\label{sec:an-example:-scalar}

The most basic example of a Fell bundle over $G$ is the trivial bundle
$\B:=G\times\C$.  Then we can identify $\gcb$ with $\ccg$.  Hence,
we can talk about
representations and \prerep s of $\ccg$ (see 
Definition~\ref{def-pre-rep} and Definition~\ref{def-repn}).
In this section, we want to review the disintegration theorem in the
scalar case as is appears in the literature.  Then we want to obtain
that formulation as a Corollary to our 
Theorem~\ref{thm-fell-disintegration}.  Let's take our time and recall
the basic definitions.

\begin{definition}
  \label{def-iso-bundle}
  If $X*\HH$ is a Borel Hilbert Bundle, then its \emph{isomorphism
    groupoid}\index{isomorphism
    groupoid}\index{Iso@$\operatorname{Iso}(X*\HH)$} is the groupoid 
  \begin{equation*}
    \operatorname{Iso}(X*\HH) := \set{(u,V,v):\text{$V:\H(v)\to\H(u)$ is
        a unitary.}}
  \end{equation*}
  with the weakest Borel structure such that
  \begin{equation*}
    (u,V,v)\mapsto\bip (Vf_{n}(v)|f_{m}(u))
  \end{equation*}
  is Borel for each $n$ and $m$ with $\set{f_{n}}$ a fundamental
  system for $X*\HH$.
\end{definition}

\begin{definition}
  \label{def-uni-repn}
  A unitary representation\index{unitary representation} of a groupoid
  $G$ with Haar system 
  $\set{\lambda^{u}}_{u\in\go}$ is a triple $(\mu,\go*\HH,L)$
  consisting of a quasi-invariant measure $\mu$ on $\go$, a Borel
  Hilbert bundle $\go*\HH$ over $\go$ and a Borel
  homomorphism\footnote{The natural thing is really a ``almost everywhere
    representation''\index{almost everywhere representation} (in the sense of
    \cite{muh:cbms}*{Remark~3.23}).  But we can produce an strict
    homomorphism in the disintegration theorem, so that is what we
    work with now.} $\hat L:G\to\operatorname{Iso}(\go*\HH)$ such that
  \begin{equation}\label{eq:64}
    \hat L(x)=\bigl(r(x),L_{x},s(x)\bigr).
  \end{equation}
\end{definition}

\begin{remark}
  \label{rem-right-form}
  In Definition~\ref{def-uni-repn}, it is important to note that
the groupoid homomorphism $\hat L:G\to\operatorname{Iso}(\go*\HH)$ has
the form specified in \eqref{eq:64}.  In general, a groupoid
homomorphism $L':G\to\operatorname{Iso}(\go*\HH)$ need only satisfy
$L'(x)=\bigl(\rho(x),L_{x}',\sigma(x)\bigr)$ for appropriate maps
$\rho,\sigma$ of $G$ into $\go$.  This makes passing from ``almost
everywhere'' homomorphisms to everywhere homomorphisms via Ramsay's
results a bit more problematic than indicated in the literature.  We
will pay attention to this detail below.
\end{remark}

Then, in analogy with Proposition~\ref{prop-strict-rep}, we have the
following.

\begin{prop}
  \label{prop-integrate}
  If $(\mu,\go*\HH,L)$ is a unitary representation of groupoid $G$,
  then we obtain a $\|\cdot\|_{I}$-norm bounded representation of
  $\ccg$ on
  \begin{equation*}
    \H:=L^{2}(\go*\HH,\mu),
  \end{equation*}
  called the integrated form of $(\mu,\go*\HH,L)$, determined by
  \begin{equation*}
    \bip(L(f)h|k)=\int_{G} f(x)
    \bip(L_{x}\bigl(h(s(x)\bigr)|
    {k\bigl(r(x)\bigr)})\Delta(x)^{-\half}\, d\nu(x).
  \end{equation*}
\end{prop}

Then the classical form of the disintegration in the scalar case is
given as follows.

\begin{thm}[Renault's Proposition 4.2]
  \label{thm-ren-4.2}\index{disintegration theorem!scalar case}
  Suppose that $L:\ccg\to\Lin(\H_{0})$ is a \prerep{} of $G\times\C$
  on $\H_{0}\subset \H$.  Then $L$ is bounded for the
  $\|\cdot\|_{I}$-norm and extends to a \emph{bona fide}
  representation of $\ccg$ on $\H$ which is equivalent to the
  integrated form of a unitary representation
  $(\mu,\go*\HH,\hat\sigma)$ of $G$.
\end{thm}

\begin{proof}
Notice that
  Theorem~\ref{thm-fell-disintegration} implies that $L$ is bounded
  and is equivalent to a \emph{bona fide} representation, still called
  $L$, on $L^{2}(\go*\HH,\mu)$ which is the integrated form of a
  strict representation $(\mu,\go*\HH,\pi)$ for a Borel $*$-functor
  $\pi:G\times\C\to \operatorname{End}(\go*\HH)$ (as in
  Proposition~\ref{prop-strict-rep}).

  For each $x\in G$, let $\sigma(x):=\pi(x,1)$.  Note that for all
  $\tau\in C$, we have $\pi(x,\tau)=\tau\sigma(x)$.  Then,
  \begin{align*}
    x\mapsto&\bip(\pi\bigl(f(x)\bigr)\xi\bigl(s(x)\bigr)|
    {\eta\bigl(r(x)\bigr)} ) \\
    &= f(x)\bip(\sigma(x)\xi\bigl(s(x)\bigr)| {\eta\bigl(r(x)\bigr)})
  \end{align*}
  is Borel for all $f\in \ccg$.  Since we can find $f_{n}\in\ccg$ such
  that $f_{n}(x)\nearrow1$ for all $x\in G$, it follows that
  \begin{equation*}
    x\mapsto \bip(\sigma(x)\xi\bigl(s(x)\bigr)|{\eta\bigl(r(x)\bigr)}
  \end{equation*}
  is Borel.  In particular, we can define a decomposable operator $P$
  in $B\bigl(L^{2}(\go*\HH,\mu)\bigr)$ by
  \begin{equation*}
    (P\xi)(u)=\sigma(u)\xi(u)\quad\text{for all $u\in\go$.}
  \end{equation*}
  Since $u\in\go$ implies that
  $\sigma(u)=\sigma(u)^{2}=\sigma(u)^{*}$, $\sigma(u)$ is a projection
  for all $u$.  Hence $P$ is a projection as well.  Furthermore,
  \begin{align}
    L(f)\xi(u)&=\int_{G}\pi\bigl(f(x)\bigr) \xi\bigl(s(x)\bigr)
    \Delta(x)^{-\half} \,d\lambda^{u}(x)\notag \\
    &= \int_{G}f(x)\sigma(x)\xi\bigl(s(x)\bigr) \Delta(x)^{-\half}
    \,d\lambda^{u}(x)\label{eq:65} \\
    &=\sigma(u)L(f)\xi(u).\notag
  \end{align}
  Therefore $L(f)=PL(f)$.  Since $L$ is nondegenerate, we must have
  $P=I$, and $\sigma(u)=I_{\H(u)}$ for $\mu$-almost all $u\in\go$.
  Therefore we can replace $\H(u)$ by $\H'(u):=\sigma(u)\H(u)$ and
  assume from here on that
  $\sigma(u)=I_{\H(u)}$.\footnote{Technically, $L^{2}(\go*\HH,\mu)$ is
    unitarily isomorphic to $L^{2}(\go*\HH',\mu)$ and we are replacing
    $L$ by its counterpart on the latter space.}  Having done so, we
  note that each $\sigma(x)$ is unitary:
  \begin{equation*}
    \sigma(x)^{*}\sigma(x)=\sigma\bigl(s(x)\bigr) = I_{\H(s(x))}
    \quad\text{and} \quad \sigma(x)\sigma(x)^{*}=\sigma\bigl(r(x)\bigr)
    = I_{\H(r(x))}.
  \end{equation*}
  Now the result follows from \eqref{eq:65} after defining
  $\hat\sigma(x):= \bigl(r(x),\sigma(x),s(x)\bigr)$.\footnote{It is
    possible that some of the $\sigma(u)$, for $u\in\go$ are zero.
    Since it may seem inappropriate to call the zero operator on the
    zero space ``unitary'', we can proceed as follows.  The set
    $Q\subset \go$ such that $\sigma(u)=0$ is Borel and clearly
    saturated.  If $F:=\go\setminus Q$, then $G=G\restr F \cup G\restr
    Q$.  Since $P=I$, $Q$ is $\mu$-null and $G\restr Q$ is
    $\nu$-null. For $u\in Q$, we may simply redefine $\H(u)$ to be
    $\C$ and for $x\in G\restr Q$ let $\hat\sigma(x):=
    \bigl(r(x),1,s(x)\bigr)$.}
\end{proof}

\printindex

%\bibliographystyle{amsxport} 
%\bibliography{references-nov01}

\begin{bibdiv}
\begin{biblist}

\bib{con:course}{book}{
      author={Conway, John~B.},
       title={A course in functional analysis},
      series={Graduate texts in mathematics},
   publisher={Springer-Verlag},
     address={New York},
        date={1985},
      volume={96},
}

\bib{dkr:ms08}{article}{
      author={Deaconu, Valentin},
      author={Kumjian, Alex},
      author={Ramazan, Birant},
       title={Fell bundles and groupoid morphisms},
        date={2008},
     journal={Math. Scan.},
      volume={103},
      number={2},
       pages={305\ndash 319},
}

\bib{din:integration74}{book}{
      author={Dinculeanu, N.},
       title={Integration on locally compact spaces},
   publisher={Noordhoff International Publishing},
     address={Leiden},
        date={1974},
        note={Translated from the Romanian},
      review={\MR{MR0360981 (50 \#13428)}},
}

\bib{dg:banach}{book}{
      author={Dupr{\'e}, Maurice~J.},
      author={Gillette, Richard~M.},
       title={Banach bundles, {B}anach modules and automorphisms of
  ${C}^*$-algebras},
   publisher={Pitman (Advanced Publishing Program)},
     address={Boston, MA},
        date={1983},
      volume={92},
        ISBN={0-273-08626-X},
      review={\MR{85j:46127}},
}

\bib{fd:representations1}{book}{
      author={Fell, James M.~G.},
      author={Doran, Robert~S.},
       title={Representations of {$*$}-algebras, locally compact groups, and
  {B}anach {$*$}-algebraic bundles. {V}ol. 1},
      series={Pure and Applied Mathematics},
   publisher={Academic Press Inc.},
     address={Boston, MA},
        date={1988},
      volume={125},
        ISBN={0-12-252721-6},
        note={Basic representation theory of groups and algebras},
      review={\MR{90c:46001}},
}

\bib{fd:representations2}{book}{
      author={Fell, James M.~G.},
      author={Doran, Robert~S.},
       title={Representations of {$*$}-algebras, locally compact groups, and
  {B}anach {$*$}-algebraic bundles. {V}ol. 2},
      series={Pure and Applied Mathematics},
   publisher={Academic Press Inc.},
     address={Boston, MA},
        date={1988},
      volume={126},
        ISBN={0-12-252722-4},
        note={Banach $*$-algebraic bundles, induced representations, and the
  generalized Mackey analysis},
      review={\MR{90c:46002}},
}

\bib{fol:real}{book}{
      author={Folland, Gerald~B.},
       title={Real analysis},
     edition={Second},
   publisher={John Wiley \& Sons Inc.},
     address={New York},
        date={1999},
        ISBN={0-471-31716-0},
        note={Modern techniques and their applications, A Wiley-Interscience
  Publication},
      review={\MR{2000c:00001}},
}

\bib{hah:tams78}{article}{
      author={Hahn, Peter},
       title={Haar measure for measure groupoids},
        date={1978},
        ISSN={0002-9947},
     journal={Trans. Amer. Math. Soc.},
      volume={242},
       pages={1\ndash 33},
      review={\MR{MR496796 (82a:28012)}},
}

\bib{hr:abstract}{book}{
      author={Hewitt, Edwin},
      author={Ross, Kenneth~A.},
       title={Abstract harmonic analysis. {V}ol. {I}: {S}tructure of
  topological groups. {I}ntegration theory, group representations},
      series={Die Grundlehren der mathematischen Wissenschaften, Band 115},
   publisher={Springer-Verlag},
     address={New York},
        date={1963},
      review={\MR{MR0156915 (28 \#158)}},
}

\bib{hof:74}{unpublished}{
      author={Hofmann, Karl~Heinrich},
       title={Banach bundles},
        date={1974},
        note={Darmstadt Notes},
}

\bib{hof:lnm77}{incollection}{
      author={Hofmann, Karl~Heinrich},
       title={Bundles and sheaves are equivalent in the category of {B}anach
  spaces},
        date={1977},
   booktitle={{$K$}-theory and operator algebras ({P}roc. {C}onf., {U}niv.
  {G}eorgia, {A}thens, {G}a., 1975)},
      series={Lecture Notes in Math},
      volume={575},
   publisher={Springer},
     address={Berlin},
       pages={53\ndash 69},
      review={\MR{58 \#7117}},
}

\bib{hofkei:lnm79}{incollection}{
      author={Hofmann, Karl~Heinrich},
      author={Keimel, Klaus},
       title={Sheaf-theoretical concepts in analysis: bundles and sheaves of
  {B}anach spaces, {B}anach {$C(X)$}-modules},
        date={1979},
   booktitle={Applications of sheaves ({P}roc. {R}es. {S}ympos. {A}ppl. {S}heaf
  {T}heory to {L}ogic, {A}lgebra and {A}nal., {U}niv. {D}urham, {D}urham,
  1977)},
      series={Lecture Notes in Math.},
      volume={753},
   publisher={Springer},
     address={Berlin},
       pages={415\ndash 441},
      review={\MR{MR555553 (81f:46085)}},
}

\bib{khoska:jram02}{article}{
      author={Khoshkam, Mahmood},
      author={Skandalis, Georges},
       title={Regular representation of groupoid {$C\sp *$}-algebras and
  applications to inverse semigroups},
        date={2002},
        ISSN={0075-4102},
     journal={J. Reine Angew. Math.},
      volume={546},
       pages={47\ndash 72},
      review={\MR{MR1900993 (2003f:46084)}},
}

\bib{khoska:jot04}{article}{
      author={Khoshkam, Mahmood},
      author={Skandalis, Georges},
       title={Crossed products of {$C\sp *$}-algebras by groupoids and inverse
  semigroups},
        date={2004},
        ISSN={0379-4024},
     journal={J. Operator Theory},
      volume={51},
      number={2},
       pages={255\ndash 279},
      review={\MR{MR2074181 (2005f:46122)}},
}

\bib{kmrw:ajm98}{article}{
      author={Kumjian, Alexander},
      author={Muhly, Paul~S.},
      author={Renault, Jean~N.},
      author={Williams, Dana~P.},
       title={The {B}rauer group of a locally compact groupoid},
        date={1998},
        ISSN={0002-9327},
     journal={Amer. J. Math.},
      volume={120},
      number={5},
       pages={901\ndash 954},
      review={\MR{2000b:46122}},
}

\bib{gal:94}{thesis}{
      author={Le~Gall, Pierre-Yves},
       title={Th\'eorie de {K}asparov \'equivariante et groupo\"\i des},
        type={Th\`ese de Doctorat},
        date={1994},
}

\bib{gal:kt99}{article}{
      author={Le~Gall, Pierre-Yves},
       title={Th\'eorie de {K}asparov \'equivariante et groupo\"\i des. {I}},
        date={1999},
        ISSN={0920-3036},
     journal={$K$-Theory},
      volume={16},
      number={4},
       pages={361\ndash 390},
      review={\MR{2000f:19006}},
}

\bib{muh:cbms}{techreport}{
      author={Muhly, Paul~S.},
       title={Coordinates in operator algebra},
 institution={CMBS Conference Lecture Notes (Texas Christian University 1990)},
        date={1999},
        note={In continuous preparation},
}

\bib{muh:cm01}{incollection}{
      author={Muhly, Paul~S.},
       title={Bundles over groupoids},
        date={2001},
   booktitle={Groupoids in analysis, geometry, and physics ({B}oulder, {CO},
  1999)},
      series={Contemp. Math.},
      volume={282},
   publisher={Amer. Math. Soc.},
     address={Providence, RI},
       pages={67\ndash 82},
      review={\MR{MR1855243 (2003a:46085)}},
}

\bib{mrw:jot87}{article}{
      author={Muhly, Paul~S.},
      author={Renault, Jean~N.},
      author={Williams, Dana~P.},
       title={Equivalence and isomorphism for groupoid {$C^*$}-algebras},
        date={1987},
        ISSN={0379-4024},
     journal={J. Operator Theory},
      volume={17},
      number={1},
       pages={3\ndash 22},
      review={\MR{88h:46123}},
}

\bib{muhwil:plms395}{article}{
      author={Muhly, Paul~S.},
      author={Williams, Dana~P.},
       title={Groupoid cohomology and the {D}ixmier-{D}ouady class},
        date={1995},
     journal={Proc. London Math. Soc. (3)},
       pages={109\ndash 134},
}

\bib{muhwil:nyjm08}{book}{
      author={Muhly, Paul~S.},
      author={Williams, Dana~P.},
       title={Renault's equivalence theorem for groupoid crossed products},
      series={NYJM Monographs},
   publisher={State University of New York University at Albany},
     address={Albany, NY},
        date={2008},
      volume={3},
        note={Available at http://nyjm.albany.edu:8000/m/2008/3.htm},
}

\bib{rae:ma88}{article}{
      author={Raeburn, Iain},
       title={Induced {$C\sp *$}-algebras and a symmetric imprimitivity
  theorem},
        date={1988},
        ISSN={0025-5831},
     journal={Math. Ann.},
      volume={280},
      number={3},
       pages={369\ndash 387},
      review={\MR{90k:46144}},
}

\bib{rw:morita}{book}{
      author={Raeburn, Iain},
      author={Williams, Dana~P.},
       title={Morita equivalence and continuous-trace {$C^*$}-algebras},
      series={Mathematical Surveys and Monographs},
   publisher={American Mathematical Society},
     address={Providence, RI},
        date={1998},
      volume={60},
        ISBN={0-8218-0860-5},
      review={\MR{2000c:46108}},
}

\bib{ram:am71}{article}{
      author={Ramsay, Arlan},
       title={Virtual groups and group actions},
        date={1971},
     journal={Advances in Math.},
      volume={6},
       pages={253\ndash 322 (1971)},
      review={\MR{43 \#7590}},
}

\bib{ram:jfa82}{article}{
      author={Ramsay, Arlan},
       title={Topologies on measured groupoids},
        date={1982},
     journal={J. Funct. Anal.},
      volume={47},
       pages={314\ndash 343},
}

\bib{ren:jot87}{article}{
      author={Renault, Jean},
       title={Repr\'esentations des produits crois\'es d'alg\`ebres de
  groupo\"\i des},
        date={1987},
     journal={J. Operator Theory},
      volume={18},
       pages={67\ndash 97},
}

\bib{ren:jot91}{article}{
      author={Renault, Jean},
       title={The ideal structure of groupoid crossed product \cs-algebras},
        date={1991},
     journal={J. Operator Theory},
      volume={25},
       pages={3\ndash 36},
}

\bib{rud:real}{book}{
      author={Rudin, Walter},
       title={Real and complex analysis},
   publisher={McGraw-Hill},
     address={New York},
        date={1987},
}

\bib{wil:crossed}{book}{
      author={Williams, Dana~P.},
       title={Crossed products of {$C{\sp \ast}$}-algebras},
      series={Mathematical Surveys and Monographs},
   publisher={American Mathematical Society},
     address={Providence, RI},
        date={2007},
      volume={134},
        ISBN={978-0-8218-4242-3; 0-8218-4242-0},
      review={\MR{MR2288954 (2007m:46003)}},
}

\bib{yam:xx87}{unpublished}{
      author={Yamagami, Shigeri},
       title={On the ideal structure of {$C^*$}-algebras over locally compact
  groupoids},
        date={unpublished manuscript, 1987},
}

\end{biblist}
\end{bibdiv}

\def\noopsort#1{}\def\cprime{$'$} \def\sp{^}
% \bib, bibdiv, biblist are defined by the amsrefs package.

\end{document}